\documentclass[10pt]{article}

\usepackage{latexsym,amssymb,amsmath,theorem}
\input diagrams

\numberwithin{equation}{section}

\newlength{\dinwidth}
\newlength{\dinmargin}
\def\endexem{\hfill{$\Box$}\medskip}

\theoremheaderfont{\scshape}
\newtheorem{thm}{Theorem}[section]
\newtheorem{prop}[thm]{Proposition}
\newtheorem{lem}[thm]{Lemma}
\newtheorem{cor}[thm]{Corollary}
\newtheorem{defin}[thm]{Definition}
\newtheorem{defprop}[thm]{Definition/Proposition}
\theorembodyfont{\rmfamily}
\newtheorem{example}[thm]{Example}
\newtheorem{rema}[thm]{Remark}

\newcommand{\bdefin}{\begin{defin}}
\newcommand{\blemma}{\begin{lem}}
\newcommand{\bprop}{\begin{prop}}
\newcommand{\btheor}{\begin{thm}}
\newcommand{\bcoro}{\begin{cor}}
\newcommand{\bconj}{\begin{conj}}
\newcommand{\bdefprop}{\begin{defprop}}
\newcommand{\bexam}{\begin{example}}
\newcommand{\edefin}{\end{defin}}
\newcommand{\elemma}{\end{lem}}
\newcommand{\eprop}{\end{prop}}
\newcommand{\etheor}{\end{thm}}
\newcommand{\ecoro}{\end{cor}}
\newcommand{\econj}{\end{conj}}
\newcommand{\brem}{\begin{rema}}
\newcommand{\erem}{\endexem\end{rema}}
\newcommand{\edefprop}{\end{defprop}}
\newcommand{\eexam}{\endexem\end{example}}

\def\1#1{{\bf #1}}
\def\2#1{{\cal #1}}
\def\3#1{{\sl #1}}
\def\4#1{{\tt #1}}
\def\5#1{{\sf #1}}
\def\6#1{{\mathfrak #1}}
\def\7#1{{\mathbb #1}}

\def\qed{\hfill{$\blacksquare$}\medskip}

 \def\f{\varphi}  
    
\def\om{\omega}   \def\D{\Delta}

\def\id{\mathrm{id}}

\newcommand{\Hom}{\mathrm{Hom}}
\newcommand{\Mor}{\mathrm{Mor}}
\newcommand{\End}{\mathrm{End}}
\newcommand{\Aut}{\mathrm{Aut}}
\newcommand{\Rep}{\mathrm{Rep}}
\newcommand{\Corep}{\mathrm{Corep}}
\newcommand{\Mod}{\mathrm{Mod}}
\newcommand{\Comod}{\mathrm{Comod}}
\newcommand{\Vect}{\mathrm{Vect}}
\newcommand{\Nat}{\mathrm{Nat}}
\newcommand{\obj}{\mathrm{Obj}}

\newcommand{\mcirc}{\,\circ\,}

\def\prf{\noindent \emph{Proof.\ }}

\def\dar{\Delta_{\mathrm{r}}}

\newcommand{\ahh}{\hat{A} \hspace{-.55ex}\hat{\rule{0ex}{2.0ex}}\hspace{.55ex}}
\newcommand{\dhh}{\hat{\D} \hspace{-.95ex}\hat{\rule{0ex}{2.05ex}}\hspace{.95ex}}

\newcommand{\aqgs}{aqg}

\newcommand{\restr}{\upharpoonright}
\newcommand{\rarr}{\rightarrow}
\newcommand{\lrarr}{\longrightarrow}
\newcommand{\impl}{\Rightarrow}
\newcommand{\oli}{{\overline{\imath}}}
\newcommand{\ve}{\varepsilon}
\newcommand{\op}{{\mbox{\scriptsize op}}}

\newcommand{\DS}{\displaystyle}

\newarrow{Congruent} 33333

\begin{document}

\title{Representations of Algebraic Quantum Groups and \\ Reconstruction Theorems for
Tensor Categories}

\author{M.\ M\"uger, J.\ E.\ Roberts and L.\ Tuset}
\date{\today}
\maketitle

\begin{abstract}
We give a pedagogical survey of those aspects of the abstract representation theory of quantum
groups which are related to the Tannaka-Krein reconstruction problem. We show that 
every concrete semisimple tensor $*$-category with conjugates is equivalent to the
category of finite dimensional non-degenerate $*$-representations of a discrete algebraic
quantum group. Working in the self-dual framework of algebraic quantum groups, we then
relate this to earlier results of S.\ L.\ Woronowicz and S.\ Yamagami. We
establish the relation between braidings and $R$-matrices in this context. Our approach
emphasizes the role of the natural transformations of the embedding functor. Thanks to the
semisimplicity of our categories and the emphasis on representations rather than
corepresentations, our proof is more direct and conceptual than previous
reconstructions. As a special case, we reprove the classical Tannaka-Krein result for
compact groups. It is only here that analytic aspects enter, otherwise we proceed in a 
purely algebraic way. In particular, the existence of a Haar functional is reduced to a
well known general result concerning discrete multiplier Hopf $*$-algebras.
\end{abstract}

\section{Introduction}
Pontryagin's duality theory for locally compact abelian groups and the Tannaka-Krein
theory for compact groups are two major results in the theory of harmonic analysis on
topological groups \cite{HR}. The Pontryagin duality theorem is a statement concerning
characters, whereas the Tannaka-Krein theorem is a statement involving irreducible unitary
representations. These two notions coincide whenever the group is both abelian and
compact.

Pontryagin's theorem can be stated more generally for locally compact quantum groups
\cite{KV1}, a notion evolving out of Kac algebras \cite{ES} and compact matrix
pseudogroups \cite{woro1}. 

The theory of representations is most naturally developed in the language of tensor
categories \cite{cwm}.  The category of finite dimensional representations of a compact
group is a symmetric tensor $*$-category with conjugates \cite{DR}, and the Tannaka-Krein
theorem tells us how to reconstruct the group from the latter.  In 1988, starting from a
tensor $*$-category \cite{glr} with conjugates and admitting a generator, but assuming
neither a symmetry nor a braiding, S.\ L.\ Woronowicz reconstructed \cite{woro2} a compact
matrix pseudogroup \cite{woro1} having the given category as its category of finite
dimensional unitary corepresentations. No general definition of a compact quantum group
existed at the time. Once it did, his proof generalizes to categories without generator. 
Alternatively, one can use the fact that every $*$-category and every compact quantum
group are inductive limits of categories with generators and of compact matrix
pseudogroups, respectively.

For the Tannaka-Krein results mentioned so far, the starting point is a concrete category, i.e.\ a
(non-full) subcategory of the tensor category of Hilbert spaces. It is conceptually more
satisfactory to start from an abstract tensor category together with a faithful tensor functor into
the category of Hilbert spaces.  In the work of N. Saavedra Rivano \cite{SR} the group associated
with a concrete symmetric tensor category was identified as the group of natural monoidal
automorphisms of the embedding functor as Tannaka had effectively done  \cite{tannaka} before the
advent of category theory. The idea of considering the natural transformations of the embedding
functor was generalized in the work of K.-H. Ulbrich \cite{ulb}, where a given concrete tensor
category is identified as the category of comodules over a Hopf algebra, cf.\ also \cite{yet}.

While the role of the natural transformations of the embedding functor is obscure in
Woronowicz's approach, they appear at least implicitly in the work of S.\ Yamagami
\cite{yama1}, who considered representations of discrete quantum groups. His approach
has two drawbacks. On the one hand, he assumes that the category is equipped with
an additional piece of structure, an `$\varepsilon$-structure'. A very similar notion
naturally arises in recent axiomatizations of tensor categories with two-sided duals
\cite{bw}, but it is quite superfluous if one works with $*$-categories. This will 
become clear in our treatment. (Yamagami has recently proved a result \cite[Theorem 3.6]{yama2}
implying that, passing if necessary to an equivalent tensor $*$-category, there is an essentially
unique $\varepsilon$-structure.) On the other hand, the von Neumann algebraic formulation of
discrete quantum groups used in \cite{yama1} is considerably more involved than current definitions,
in that it unnecessarily mixes analytic and algebraic structures.

In neither of \cite{woro1,yama1} does one find a complete proof of equivalence between the
given category and the representation category of the derived quantum group. Also
questions of uniqueness and braidings have not been addressed.

It is well known that the self-duality of the category of finite dimensional Hopf algebras
breaks down in the infinite dimensional case. Motivated by the desire to find a purely
algebraic framework for quantum groups admitting a version of Pontryagin duality, A. Van
Daele developed the theory of algebraic quantum groups (\aqgs) \cite{VD}. This is achieved
by admitting non-unital algebras and requiring the existence of a positive and
left-invariant functional, the Haar functional. All compact and all discrete quantum
groups are \aqgs, and every aqg has a unique analytic extension to a locally compact
quantum group \cite{KD}, having an equivalent von Neumann algebraic version \cite{KV2}. 
A discrete multiplier Hopf $*$-algebra \cite{VD2} can be shown to have a Haar functional
rendering it a discrete aqg. 

The purpose of the present paper is to give a coherent and reasonably complete survey of
the Tannaka-Krein theory of quantum groups. The only other review we are aware of is
\cite{JS1}, which appeared ten years ago. At the time, no appropriate self-dual category
of quantum groups existed. (In the same year, P. Podl\'{e}s and Woronowicz, in defining a
discrete quantum group, took the first step in this direction \cite{woro3}.) The approach
to Tannakian categories motivated by algebraic geometry is well reviewed in \cite{breen},
where, however, only symmetric tensor categories are considered.

Our approach to generalized Tannaka-Krein theory adopts the philosophy on quantum groups
in \cite{woro2,yama1}, meaning a self-dual category with emphasis on $*$-categories. In
contradistinction to these authors we wish to distinguish categorical from quantum group
aspects as well as algebraic from analytic aspects as far as possible. Following
\cite{SR,DM,ulb,yet}, we emphasize natural transformations of the embedding functor.
Yet, our use of natural transformations is more direct and we work with representations
rather than corepresentations.  The algebra product is just the composition of natural
transformations and the coproduct is defined directly in terms of the tensor structure,
the tensor unit giving rise to the counit. This yields a quantum semigroup, whose
discreteness is an immediate consequence of the semisimplicity of the category. (The
semisimplicity of $*$-categories avoids appealing to the theorem of Barr-Beck
in \cite{SR,DM} and to proceed in a pedestrian way, using only the definition of
natural transformations.) The coinverse now arises from the conjugation in the
category. The result is a discrete multiplier Hopf $*$-algebra. Thus our reconstruction of
the discrete aqg is purely algebraic, the existence of the Haar functional following by
quantum group theory. 

The selfduality of the category of \aqgs\ and the existence of analytic extensions allow us
to to relate our reconstruction result to those of Woronowicz and Yamagami on one hand and
the purely algebraic ones on the other \cite{ulb,yet}. In particular, making use of the
universal unitary corepresentation of an \aqgs\ introduced by J. Kustermans \cite{ku-e}, we
prove that the tensor $*$-category of $*$-representations 
(or modules) of a discrete aqg $(A,\D)$ is equivalent to the tensor $*$-category of
unitary corepresentations (or comodules) of the dual compact aqg $(\hat{A},\hat{\D})$. We
provide these categories with conjugates (for representations of discrete \aqgs\ this
has not been worked out before). We also show that the tensor $*$-category of pointwise
continuous finite dimensional $*$-representations of a sub Hopf $*$-algebra of the maximal
dual (or Sweedler dual) Hopf $*$-algebra of a compact quantum group is equivalent to any
of the tensor $*$-categories mentioned above, whenever the Hopf $*$-algebra separates the
regular functions associated with the compact quantum group. This is a useful result since
it applies to the deformed universal enveloping Lie algebras $U_q (\6g)$ of M. Jimbo and
V.G. Drinfeld. It therefore links these axiomatic quantum group results to the more
familiar context of quantum groups given by deformations of semisimple Lie algebras $\6g$.
Indeed, this is how the latter can be shown to produce compact (or discrete) quantum
groups and how Woronowicz \cite{woro2} constructed the compact matrix pseudogroup 
$SU_q(N)$. 

The correspondence between infinite dimensional representations and corepresentations was
already established in \cite{woro3} and \cite{ku-e}, but only for the objects of the
respective categories, morphisms and the tensor structure were not considered. In
\cite{VDZ} tensor structure and braiding were taken into account in a purely algebraic
context and, while studying the amenability of quantum groups, the correspondence between
tensor C*-categories of infinite dimensional representations and corepresentations was
established in \cite{BCT,BT}, for the various analytic versions of aqg and lcqg. The
latter results rely on the theory of infinite dimensional representations and
corepresentations and of the construction of the universal corepresentation for lcqg
developed in \cite{ku-u}.  But none of this work touched on Tannaka-Krein
reconstruction, since conjugates do not exist in the infinite dimensional
case. Another type of reconstruction result for tensor C*-categories involving infinite
dimensional objects was undertaken in \cite{DPR} from the point of view of multiplicative
unitaries and the regular corepresentation. 

For a discrete aqg $(A,\D)$ we establish a bijection between braidings of the category
$\Rep_f(A,\D)$ and $R$-matrices in the multiplier algebra $M(A\otimes A)$. To the best of
our knowledge this is the first such result rigorously proven in an axiomatic framework
for quantum groups. If the category is symmetric and the embedding functor maps the
braiding into the canonical braiding of the category of Hilbert spaces, the discrete aqg
is cocommutative.

For any discrete aqg there is a compact group $G$, the intrinsic group, and for
cocommutative $(A,\D)$ we prove an equivalence $\Rep_f(A,\D)\simeq\Rep_f\,G$ of tensor
$*$-categories. This is the only point in our approach were analysis plays a role, in that
we use the theorems of Gelfand, Krein-Milman, Stone-Weierstrass, etc.  In the theory of
locally compact quantum groups it is well known that commutative and cocommutative quantum
groups are Kac algebras. Furthermore, commutative (resp.\ cocommutative) lcqg are of the
form $(C(G),\D)$ (resp.\ $(C^*_r (G),\D)$) for a locally compact group $G$.  Finally, a
commutative compact aqg is the algebra of regular functions on a compact group $G$, but
cocommutative discrete \aqgs\ are inconvenient to characterize. This is why, alternatively,
we give the more instructive and direct proof of the above equivalence of categories. In
passing we give a description of a cocommutative aqg in terms of the intrinsic group $G$.

It cannot be emphasized strongly enough that all Tannaka-Krein type results discussed so
far depend on the tensor $*$-category being concrete, i.e.\ coming with a faithful tensor
functor into the category of Hilbert spaces. There are applications in pure mathematics
and in quantum field theory, where such a functor is not given a priori.  For symmetric
categories, it was first shown by S.\ Doplicher and J.\ E.\ Roberts \cite{DR} that such a
functor always exists. An alternative approach in a more algebraic setting was given by
P. Deligne \cite{del}.  The questions of existence and uniqueness of an embedding functor
will be addressed anew in a sequel to this paper.

The above discussion did not follow the order of presentation. Let us therefore give a
brief overview of the organization of this paper. In the next section we provide some
preliminaries on tensor $*$-categories and \aqgs.  Concerning the former we are quite
brief, since much of this material is almost universally known. Concerning the latter we
focus in particular on the discrete and compact cases and discuss the examples related to
groups.  In Section 3 we treat the representation and corepresentation theory of discrete
and compact \aqgs, respectively, from a tensor $*$-category point of view.  To this
end, we recall the universal corepresentation due to Kustermans, and we discuss
conjugates in these categories.  The special case of cocommutative discrete \aqgs\ is
considered in Section 4. Section 5 is the heart of this paper.  There we construct a
discrete quantum semigroup from a tensor $*$-category, deriving the coinverse from the
conjugation and leading to a discrete aqg. The precise statement of the
generalized Tannaka-Krein theorem for quantum groups is then made in Theorem
\ref{t-main-equiv}.  In the final Section 6 we establish the bijection between braidings
and $R$-matrices. In the case of a symmetric tensor $*$-category with symmetric embedding
functor we recover the classical Tannaka-Krein theorem for compact groups.


\section{Preliminaries}
\subsection{Tensor Categories} \label{ss-categ}
For the definitions of categories, functors and natural transformations we refer, e.g., to
\cite{cwm}.  In this subsection we briefly recall some of the less standard notions of category
theory which will be needed here, others will be introduced as we proceed.  We may occasionally say
`arrows' instead of `morphisms' in order to avoid confusion with algebra homomorphisms. All
categories which we consider are essentially small, i.e.\ equivalent to a small category.  We mostly 
speak of `tensor categories' rather than `monoidal categories' but cannot avoid the adjective
`monoidal'. In view of the coherence theorems for (braided, symmetric) tensor categories 
we may assume all tensor categories to be strict, satisfying $X\otimes(Y\otimes Z)=(X\otimes
Y)\otimes Z$ for all $X,Y,Z$. Following widespread use, we also consider the tensor categories of
vector spaces as strict, appealing to the canonical isomorphisms to identify $X\otimes(Y\otimes Z)$
with $(X\otimes Y)\otimes Z$. We often write $XY$ in the place of $X\otimes Y$.  By $\2H$ we mean
the strictification of $\mathrm{Hilb}$. (We also suppress all other canonical isomorphisms in
$\mathrm{Hilb}$, identifying $B(H\otimes K)=B(H)\otimes B(K)$.) Note, however, that one cannot
assume all tensor functors to be strict without losing generality. Thus we need the following.

\bdefin \label{d-mon-func}
Let $\2C, \2D$ be (strict) tensor categories. A functor $F: \2C\rarr\2D$ is weakly monoidal if 
there exist morphisms $d_{X,Y}: F(X)\otimes F(Y)\rarr F(X\otimes Y)$ for all $X,Y\in\2C$ 
and a morphism $e: \11_\2D\rarr F(\11_\2C)$ such that
\begin{enumerate}
\item The morphisms $d_{X,Y}$ are natural w.r.t.\ both arguments.
\item For all $X,Y,Z\in\2C$ the following diagram commutes:
\begin{diagram}
F(X)\otimes F(Y)\otimes F(Z) & \rTo^{d_{X,Y}\otimes\id_{F(Z)}} & F(X\otimes Y)\otimes F(Z) \\
\dTo^{\id_{F(X)}\otimes d_{Y,Z}} & & \dTo_{d_{X\otimes Y,Z}} \\
F(X)\otimes F(Y\otimes Z) & \rTo_{d_{X, Y\otimes Z}} & F(X\otimes Y\otimes Z)
\end{diagram}
\item The following compositions are the identity morphisms of $F(X)$
\begin{diagram}
F(X) \equiv F(X)\otimes \11_\2D & \rTo^{\id_{F(X)}\otimes e} & F(X)\otimes F(\11_\2C) & \rTo^{d_{X,\11}} &
   F(X\otimes\11_\2C) \equiv F(X)
\end{diagram}
\begin{diagram}
F(X) \equiv \11_\2D\otimes F(X) & \rTo^{e\otimes\id_{F(X)}} & F(\11_\2C)\otimes F(X) & \rTo^{d_{\11,X}} &
   F(\11_\2C\otimes X) \equiv F(X)
\end{diagram}
for all $X\in\2C$.
\end{enumerate}
The functor $F$ is monoidal if $e$ and all the $d_{X,Y},\ X,Y\in\2C$ are isomorphisms. If
$\2C,\2D$ are tensor $*$-categories and $F$ is $*$-preserving, the isomorphisms $e, d_{X,Y}$ are
required to be unitary. 
\edefin

We will mainly be interested in linear categories over the complex field. Viz., for all $X,Y$,
$\Mor(X,Y)$ is a complex vector space, and the compositions $\circ,\otimes$ of morphisms are
bilinear. All functors are supposed $\7C$-linear. A {\it $*$-operation} on a $\7C$-linear (tensor)
category $\2C$ is a map which assigns to an arrow $s\in\Mor(X, Y)$ another arrow
$s^*\in\Mor(Y,X)$. This map has to be antilinear, involutive ($s^{**}=s$), contravariant 
($(s\,\circ\, t)^*=t^*\mcirc  s^*$) and, if $\2C$ is monoidal, monoidal ($(s\times t)^*=s^*\times t^*$). 
A $*$-operation is {\it positive} iff $s^*\mcirc  s=0$ implies $s=0$. A (tensor) $*$-category is a 
$\7C$-linear (tensor) category with a positive $*$-operation. A morphism $s:X\rarr Y$ is an isometry
iff it satisfies $s^*\circ s=\id_X$. A morphism $s$ is unitary iff $s$ and $s^*$ are isometries. A
functor $F$ between $*$-categories is $*$-preserving if $F(s^*)=F(s)^*$ for every morphism $s$.

An object $Z$ is a direct sum of $X_1,X_2$, denoted $Z\cong X_1\oplus X_2$, if there are isometries 
$s_i: X_i\rarr Z, i=1,2$ such that $s_1\circ s_1^*+s_2\circ s_2^*=\id_Z$. A $*$-category `has direct
sums' if for every pair $X_1, X_2$ there exists a direct sum. An object $Y$ is a subobject of $X$ if
there exists an isometry $s:Y\rarr X$. A $*$-category `has subobjects' if for every 
$p=p\circ p=p^*\in\End(X)$ there exist $Y$ and an isometry $s:Y\rarr X$ such that $s\circ s^*=p$. 
An object $X$ in a $\7C$-linear category is irreducible if it is non-zero and $\End\,X=\7C\,\id_X$. 

\bdefin A $*$-category is semisimple if it has finite dimensional spaces of morphisms, a zero object, 
direct sums and subobjects. A tensor $*$-category is semisimple if, in addition, the tensor unit
$\11$ is irreducible. 
\edefin
In a semisimple $*$-category $\2C$, $\End(X)$ is a finite dimensional $C^*$-algebra for every $X$,  
and every object is a finite direct sum of irreducible objects. (It is well-known that a category
which is semisimple in our sense is abelian and semisimple in the usual sense, i.e.\ all exact
sequences split.) 

Let $\2C$ be a tensor $*$-category and $X\in\2C$. A `solution of the conjugate equations' is a
triple $(\overline{X},r,\overline{r})$, where $\overline{X}\in\2C$ and 
$r:\11\rarr \overline{X}\otimes X,\ \overline{r}:\11\rarr X\otimes\overline{X}$ satisfy
\begin{eqnarray*}
\overline{r}^*\otimes\id_X\mcirc\id_X\otimes r &=&\id_X, \\
   r^*\otimes\id_{\overline{X}}\mcirc\id_{\overline{X}}\otimes\overline{r} &=&\id_{\overline{X}}.
\end{eqnarray*}
A tensor $*$-category $\2C$ `has conjugates' if there is a solution of the conjugate equations for
every $X\in\2C$. A solution $(X,r,\overline{r})$ is normalized iff 
$r^*\circ r=\overline{r}^*\circ\overline{r}$. 
It is a standard solution iff there are irreducible objects $X_i, i\in I_X$, solutions 
$(\overline{X_i},r_i,\overline{r}_i)$ of the conjugate equations and isometries
$v_i: X_i\rarr X, w_i: \overline{X_i}\rarr\overline{X}, i\in I_X$ satisfying
\[ v_i^*\mcirc v_j=\delta_{ij}\,\id_{X_i}, \quad\quad\sum_i v_i\mcirc v_i^*=\id_X, \]
similar equations for $w_i$, and we have
\[ r=\sum_i w_i\otimes v_i\mcirc r_i, \quad \overline{r}=\sum_i v_i\otimes w_i\mcirc\overline{r_i}. \]
We define the (intrinsic or categorical) dimension $d(X)\in\7R_+$ of $X$ by $r^*\circ r=d(X)\id_\11$
where $(\overline{X},r,\overline{r})$ is a normalized standard solution. 
One can prove the following facts, cf.\ \cite{LR}. The dimension is additive under direct sums and
multiplicative under tensor products. It takes values in the set 
$\{ 2\cos\frac{\pi}{n}, n=3,4,\ldots \}\cup[2,\infty)$, in particular $d(X)\ge 1$ with
$d(X)=1$ iff $X\otimes\overline{X}\cong\11$ iff $X$ is invertible, i.e.\ there exists $Y$ such
that $X\otimes Y\cong\11$. If $d(X)=1$ then $X$ is irreducible.
In the category $\2H$ of Hilbert spaces we have $d_\2H(H)=\dim_\7C H$. 

We briefly comment on a somewhat more general setting. A $C^*$-(tensor) category is a (tensor)
$*$-category where $\Mor(X,Y)$ is a Banach space for every pair $(X, Y)$ of objects and the norms
satisfy $\|X^*\mcirc  X\|=\|X\|^2$. A $W^*$-category is a $C^*$-category where every $\Mor(X,Y)$ is
the dual of a Banach space for every pair $(X,Y)$. In a $C^*$-category with conjugates and
irreducible tensor unit all spaces of morphisms are finite dimensional. This is useful in
applications where this finite dimensionality is not known a priori, like in quantum field
theory. Conversely, every $*$-category which is semisimple in our sense is a $W^*$-category, cf.\
\cite{mue06}. In a $W^*$-category, every morphism $s: X\rarr Y$ has a polar decomposition $s=pu$,
where $p$ is positive and $u$ a partial isometry. As a consequence, whenever $\Mor(X,Y)$ contains  
a split monic (or isomorphism), it also contains an isometry (respectively, unitary). This shows
that most of the above definitions, e.g., of direct sums, are equivalent to the the usual ones as
given, e.g., in \cite{cwm}.

If $\2C$ is a semisimple tensor $*$-category we denote the set of isomorphism classes of irreducible
objects by $I_\2C$. Let $(X_i, i\in I_\2C)$ be a complete set of irreducible objects and write
$d_i=d(X_i)$. We have $0\in I_\2C$ such that $X_0\cong\11$. If $\2C$ has conjugates, $I_\2C$ comes
with an involution $i\mapsto\bar{\imath}$ such that 
$X_{\bar{\imath}}$ is a conjugate of $X_i$. For $i,j,k\in I_\2C$ we define 
$N_{ij}^k=\dim\Mor(X_k,X_i\otimes X_j)$. These numbers satisfy the following properties.
\begin{enumerate}
\item For every pair $(i,j)$ there are only finitely many $k\in I_\2C$ such that $N_{ij}^k\ne 0$. We
have $X_i\otimes X_j\cong\bigoplus_{k\in I_\2C} N_{ij}^k\,X_k$ and thus
$d_id_j=\sum_k N_{ij}^k d_k$.
\item $\DS \sum_l N_{ij}^l N_{lk}^m = \sum_n N_{in}^m N_{jk}^n$ for all $i,j,k,m\in I_\2C$.
\item 
$\DS N_{ij}^k=N_{k\overline{\jmath}}^i=N_{ \overline{\imath}k}^j=N_{\overline{k}i}^{\overline{\jmath}}=N_{j\overline{k}}^{\overline{\imath}}=N_{\overline{\jmath},
\overline{\imath}}^{\overline{k}}$ for all $i,j,k\in I_\2C$.
\item $\DS N_{ij}^0=\delta_{i,\overline{\jmath}}$. 
\item If $\2C$ is braided, cf.\ Section \ref{s-braid}, then $N_{ij}^k=N_{ji}^k$ for all $i,j,k\in
I_\2C$. 
\end{enumerate}


\subsection{Algebraic Quantum Groups}
In this subsection we briefly outline those aspects of the theory of \aqgs\ that will be needed in
the sequel. For the details and proofs, see the original references \cite{VD2,VD}.

Every algebra will be a (not necessarily unital) associative
algebra over the complex field~$\7C$. The identity map on a
set $V$ will be denoted by $\iota$.
If $V$ and $W$ are linear spaces, $V'$ denotes the linear
space of linear functionals on $V$ and  $V \otimes W$ denotes
the linear space tensor product of $V$ and~$W$. The {\em flip
map} $\sigma$ from ${V \otimes W}$ to ${W \otimes V}$ is the
linear map sending ${v\otimes w}$ onto ${w\otimes v}$, for all
$v\in V$ and $w\in W$. If $V$ and $W$ are Hilbert spaces,
${V\otimes W}$ denotes their Hilbert space tensor product; we
denote by $B(V)$ and $B_0(V)$ the C*-algebras of bounded
linear operators and compact operators on~$V$, respectively.
If $V$ and $W$ are algebras, ${V \otimes W}$ denotes their
algebra tensor product. If $V$ and $W$ are C*-algebras, then
${V\otimes W}$ will denote their C*-tensor product with
respect to the minimal C*-norm.

An algebra $A$ is {\it non-degenerate} if for any $a\in A$ such that
$ab=0$ for all $b\in A$ or $ba=0$ for all $b\in A$, we have $a=0$.
Obviously, all unital algebras are non-degenerate. If $A$ and $B$ are non-degenerate algebras, so is
$A\otimes B$. From now on {\it all} algebras are assumed to be non-degenerate.

Let $A$ be a $*$-algebra and denote by $\End\,A$ the unital algebra of linear maps from
$A$ to itself. Let
\[ M(A)= \{ x\in\End\,A \ | \ \exists y\in\End\,A\ \mbox{such that}\ 
     x(a)^* b = a^* y(b) \quad \forall a,b \in A \}. \]
Then $M(A)$ is a unital subalgebra of $\End\,A$. The linear map $y$ associated to a given 
$x\in M(A)$ is uniquely determined by non-degeneracy and we denote it by $x^*$. The unital
algebra $M(A)$ becomes a $*$-algebra when endowed with the involution ${x \mapsto x^*}$. 
This unital $*$-algebra is called the {\em multiplier algebra} of $A$.

Suppose that $A$ is an ideal in a $*$-algebra~$B$. For $b \in B$, define
$L_b \in M(A)$ by $L_b(a)=b a$, for all $a \in A$. Then the map 
$L\colon B \rightarrow M(A)$, $b \mapsto L_b$, is a homomorphism. If $A$ is an 
{\em essential} ideal in $B$ in the sense that an element $b$ of~$B$ is
necessarily  equal to zero if $b a=0$, for all $a\in A$, or $a
b=0$, for all $a\in A$, then $L$ is injective. In particular,
$A$ is an essential ideal in itself (by
non-degeneracy) and therefore we have an injective homomorphism 
$L\colon A \rightarrow M(A)$. We identify the 
image of $A$ under $L$ with $A$. Then $A$ is an essential ideal of
$M(A)$. (In fact, $M(A)$ is the `largest' algebra containing $A$ as an essential ideal.) 
Obviously, $M(A) = A$ iff $A$ is unital. 

If $A$ and $B$ are $*$-algebras, then it is
easily verified that  $A\otimes B$ is an essential
self-adjoint ideal in $M(A)\otimes M(B)$. Hence, by the
preceding remarks, there exists a canonical injective
$*$-homomorphism from $M(A)\otimes M(B)$ into $M(A\otimes B)$.
We use this to identify $M(A)\otimes M(B)$ as a unital
$*$-subalgebra of $M(A\otimes B)$. In general, these algebras
are not equal.
A linear map $\pi \colon A \to B$ is said to be  {\em non-degenerate} if $\pi(A)B=B$ and $B\pi(A)=B$. 
Here, as elsewhere, $\pi(A)B$ denotes the linear span of ${\{\pi(a)b\mid a\in A,\  b\in B\}}$.
Whenever $\pi$ is non-degenerate and multiplicative (resp.\ non-degenerate and antimultiplicative), 
there exists a unique extension to a unital homomorphism (resp.\ antihomomorphism) 
${\overline{\pi} \colon M (A) \rightarrow M (B)}$. We shall henceforth
use the same symbol $\pi$ to denote the original map and its
extension~$\overline{\pi}$. A {\em representation} of a $*$-algebra is a non-degenerate homomorphism
$\pi:A\rarr B(K)$, where $K$ is a Hilbert space. 

If $\omega$ is a linear functional on~$A$ and $x \in M(A)$, we
define the linear functionals $x \omega$ and $\omega x$ on $A$
by setting $(x \omega)(a) = \omega(a x)$ and $(\omega x)(a) =
\omega(x a)$, for all $a \in A$.
We need the leg numbering notation.  
Take three $*$-algebras $A$, $B$, $C$. It can be shown that there exists a 
non-degenerate $*$-homomorphism 
$\theta_{13} : A \otimes C \rightarrow M(A \otimes B \otimes C)$ such that 
$\theta_{13}(a \otimes c) = a \otimes 1 \otimes c$, for all $a \in A$, $c \in C$. Thus, it
has a unique extension to  
$M(A \otimes C)$. Set $x_{13} = \theta_{13}(x)$, for all $x \in M(A \otimes C)$.
The other variants of the leg numbering notation are defined similarly.

The triple product $A\otimes A\otimes A$ is an essential ideal of both $M(A\otimes A)\otimes A$
and $A\otimes M(A\otimes A)$, thus $M(M(A\otimes A)\otimes A)\subset M(A\otimes A\otimes A)$ and
$M(A\otimes M(A\otimes A))\subset M(A\otimes A\otimes A)$. 

\begin{defin} \label{def-mult-hopf-alg}
A multiplier Hopf $*$-algebra $(A, \D)$  consists of a $*$-algebra $A$ and a
$*$-homo\-mor\-phism $\D$ from $A$ into $M(A\otimes A)$ 
such that
\begin{enumerate}
\item $(\D \otimes \iota)\D = (\iota \otimes \D)\D$.
\item The linear mappings $T_1,  T_2$ from $A \otimes A$ into $M(A \otimes A)$ such that 
\begin{eqnarray*}
\ T_1(a \otimes b) &=& \D(a)(b \otimes 1), \\
  T_2(a \otimes b) &=& \D(a)(1 \otimes b) 
\end{eqnarray*}
for all $a, b \in A$, are bijections from $A \otimes A$ to $A \otimes A$.
\end{enumerate}
\end{defin}

\noindent Here the (unique) extension of $\D \otimes \iota: A\otimes A\rarr M(A\otimes A\otimes A)$ 
to $M(A\otimes A)$ is understood and similarly for $\iota \otimes \D$. Condition (ii)
implies that $\D(a)(b\otimes 1), \D(a)(1\otimes b)\in M(A\otimes A)$ lie in fact in
$A\otimes A$.

We say that two multiplier Hopf $*$-algebras $(A_1, \D_1 )$ and $(A_2, \D_2 )$
are {\em isomorphic}, and write $(A_1, \D_1 )\cong (A_2, \D_2 )$, if there exists 
a bijective $*$-homomorphism $\theta : A_1\rarr A_2$ such that $(\theta\otimes\theta )\D_1 =\D_2\theta$. 

The following result shows that multiplier Hopf $*$-algebras share most properties with
Hopf algebras. Let $m: A\otimes A\rarr A$ denote the linear extension of the multiplication
map. Note that $m$ is a $*$-homomorphism iff $A$ is commutative. 

\bprop 
\label{p-converse}
Let $(A,\D)$ be a multiplier Hopf $*$-algebra. There exist a unique $*$-homo\-mor\-phism 
$\ve: A\rarr\7C$ and a unique invertible antimultiplicative $S\in\End\,A$ such that 
\begin{enumerate}
\item $\DS (\ve\otimes\iota)(\D(a)(1\otimes b)) = ab$,
\item $\DS (\iota\otimes\ve)(\D(a)(b\otimes 1)) = ab$,
\item $\DS m(S\otimes\iota)(\D(a)(1\otimes b)) = \ve(a)b$,
\item $\DS m(\iota\otimes S)(\D(a)(b\otimes 1)) = \ve(a)b$
\end{enumerate}
for all $a,b\in A$. We call $\ve$ the counit and $S$ the coinverse. Conversely, if 
$(A\otimes 1)\D (A)\subset A\otimes A$ and $(A\otimes 1)\D (A)\subset A\otimes A$ and there exist
linear maps $\varepsilon :A\rarr\7C$ and $S:A\rightarrow A$ satisfying properties 1-4 with $S$
invertible, then $(A,\D)$ is a multiplier Hopf $*$-algebra. 
\eprop

We denote  $m_\op=m\sigma$ and $\D_\op=\sigma\D$ where 
$\sigma(a\otimes b)=b\otimes a$ is the flip $*$-automorphism of $A\otimes A$ and 
$M(A\otimes A)$. $A_\op$ denotes the vector space $A$ with multiplication $m_\op$.
The following result follows easily from the uniqueness of the counit and coinverse.

\bprop The antipode $S$ is invertible and
\begin{enumerate}
\item $\DS S^{-1}(a)=S(a^*)^*\quad\forall a\in A$.
\item $\DS(S\otimes S)(\D(a)(1\otimes b)) = (1\otimes S(b))\D_\op S(a) \quad\forall a,b\in A$.
\item $\DS \ve S=\ve$.
\item $(A_\op, \D)$ and $(A,\D_\op)$ are multiplier Hopf $*$-algebras with  
counit $\ve$ and coinverse $S^{-1}$.
\end{enumerate}
\eprop

\bdefin A linear functional $\omega$ on a multiplier Hopf $*$-algebra $(A,\D)$ is
called  
\begin{enumerate} 
\item positive iff $\omega(a^*a)\geq 0\quad\forall a\in A$,
\item faithful iff it is positive and $\omega(a^*a) = 0 \Rightarrow a =0 \quad \forall a\in A$,
\item left-invariant iff
$(\iota\otimes\omega )(\D(a)(b\otimes 1)) = \omega(a)b \ \ \forall a,b \in A$,
\item right-invariant iff
$(\omega\otimes\iota)((1\otimes b)\D(a)) = \omega(a)b \ \ \forall a,b \in A$.
\end{enumerate}
\edefin

Thus $\omega$ is left-invariant iff $\omega S$ is right-invariant.

\begin{defin}
An aqg is a multiplier Hopf $*$-algebra $(A, \D)$ which admits 
a non-zero left-invariant positive linear functional $\varphi$.
\end{defin}

\btheor \label{thm-faithful}
Let $(A, \D)$ be an aqg with positive left-invariant functional
$\varphi$. Then 
\begin{enumerate}
\item $\varphi$ is faithful.
\item If $\omega$ is a left-invariant functional, there exists $c\in\7C$ such that $\omega=c\varphi$. 
\item If $\omega$ is a right-invariant functional, there exists $c\in\7C$ such that 
$\omega=c\varphi S$. 
\end{enumerate}
\etheor

An immediate consequence is that the $*$-operation is positive, viz.\ 
$a^*a=0\ \impl\ a=0$, which clearly is stronger than non-degeneracy. Furthermore, there
exists a unique complex number $\mu$ such that $\varphi S^2 = \mu \varphi$. It can be
proved that $|\mu| = 1$, but it has recently been discovered that $\mu\neq 1$ for
the quantum group version of $ax+b$ \cite{VD3}. Every aqg $(A,\D )$ has the property that
to any $a\in A$, there exists $c\in A$ such that $ac=ca=a$.
Two \aqgs\ are said to be {\it isomorphic} if they are isomorphic as multiplier Hopf $*$-algebras
and we use the same notation $\cong$ to denote this.

Set $\psi = \varphi S$. Then $\psi$ is a non-zero right-invariant linear functional on
$A$. However, in general, $\psi$ will not be positive. It is known that there
exists a non-zero positive right-invariant linear functional on $A$. 
To $\varphi$ there existence a unique automorphism $\rho$ on $A$ such that $\varphi(ab) =
\varphi(b\rho (a))$, for every $a, b \in A$. We refer to this as the weak KMS-property of
$\varphi$. Moreover, we have $\D \rho = (S^2 \otimes \rho)\D$ and 
$\rho(\rho(a^*)^*) = a$, for every $a \in A$. 
Also there exists an automorphism $\rho'$ for the right-invariant functional $\psi$, that is,
$\rho'$ is an automorphism on $A$ such that $\psi(ab) = \psi(b\rho '(a))$, for every 
$a, b \in A$. 

It is possible to introduce a modular function for \aqgs. It is an
invertible element $\delta$ in $M(A)$ such that 
$(\varphi \otimes \iota)(\D (a)(1\otimes b)) = \varphi(a)\delta b$, for every 
$a, b \in A$. This modular function satisfies

\[ \D (\delta) = \delta \otimes \delta, \quad \varepsilon (\delta) = 1, \quad S(\delta) = \delta^{-1}. \]

\noindent As in the classical case we have
\[\varphi(S(a)) = \varphi(a\delta)=\mu\varphi (\delta a) \quad \forall a\in A. \]

We now discuss duality within the category of \aqgs.
Define a subspace $\hat A$ of $A'$ by
\[\hat A = \{\varphi a \ | \ a \in A\} = \{a \varphi\  | \ a \in A \}. \]
Then $\hat A$ is a non-degenerate $*$-algebra under the definitions:
\begin{enumerate}
\item For every $\omega_1, \omega_2 \in \hat A$ and $a \in A$, we have 
$(\omega_1 \omega_2)(a) = (\omega_1 \otimes \omega_2)(\D(a))$.
\item For every $\omega \in \hat A$ and $a \in A$, define 
$\omega^*(a) = \overline {\omega(S(a)^*)}$.
\end{enumerate}

The comultiplication $\hat \D$ is defined on $\hat A$ by 
$\hat \D(\omega)(x \otimes y) = \omega(xy)$, for every $\omega \in \hat A$ and 
$x, y \in A$. For this to make sense, $M(A \otimes A)$ should be embedded in $(A \otimes A)'$
in a proper way. A definition of the comultiplication $\hat \D$ that does not use such
an embedding can be found in \cite{ku-e}. Then $(\hat A, \hat \D)$
is a multiplier Hopf $*$-algebra with counit $\hat \varepsilon$ and the antipode $\hat S$
given by: 
\begin{enumerate}
\item $\hat \varepsilon(\omega) = \omega(1)$, for every $\omega \in \hat A$.
\item $\hat S(\omega)(a) = \omega(S(a))$, for every $\omega \in \hat A$ and $a \in A$. 
\end{enumerate}

Define $\hat a = a\varphi \in \hat A$, for $a \in A$. The map $A \rightarrow \hat A$ 
sending $a$ to $\hat a$ is a bijection, and is referred to as the Fourier transform. Next define the
linear functional $\hat{\psi}$ on $\hat A$ by setting 
$\hat{\psi}(\hat a) = \varepsilon(a)$, for every 
$a \in A$. It is possible to prove that $\hat\psi$ is right-invariant, and that
$\hat \psi(\hat a ^* \hat a)=\varphi (a^* a)$, for every $a \in A$. 
It follows that $\hat \psi$ is a non-zero positive linear
functional on $\hat A$, and that the dual $(\hat A, \hat \D)$ is an
aqg. Let $\hat{\rho}$ denote the automorphism of $\hat{A}$ such that 
$\hat{\psi}(\hat a\hat b)=\hat{\psi}(\hat b\hat{\rho}(\hat a ))$, for all $a,b\in A$.
The aqg version of Pontryagin's duality theorem takes the following
form:  

\btheor
\label{dual}
Let $(A, \D)$ be an aqg. Then the double dual 
$(\ahh, \dhh)$ is an aqg isomorphic to 
$(A, \D)$. More precisely, there exists a canonical $*$-isomorphism  
$\theta : A \rightarrow \hat {\hat A}$ such that $\theta(a)(\omega) = \omega(a)$, 
for all $a \in \hat A$, and that satisfies the equation 
$(\theta \otimes \theta)\D = \hat {\hat \D} \theta$. 
\etheor

\noindent Note that an aqg $(A,\D )$ is commutative iff $(\hat{A},\hat{\D})$ is cocommutative. 

We shall need to consider an object associated to an aqg called its analytic extension. See
\cite{KD} for full details.  If $A$ is a C*-algebra, then $M(A)$ denotes the usual multiplier
algebra of $A$, and $*$-homomorphism between C*-algebras are assumed to be non-degenerate in the
usual operator algebraic sense.  We need first to recall the concept of a GNS~pair. Suppose given a
positive linear functional $\omega$ on a $*$-algebra~$A$. Let $H$ be a Hilbert space, and let
${\Lambda\colon A\to H}$ be a linear map with dense range for which
$(\Lambda(a),\Lambda(b))=\om(b^*a)$, for all $a,b\in A$.  Then we call $(H,\Lambda)$ a {\em
GNS~pair} associated to~$\om$.  As is well known, such a pair always exists and is essentially
unique. 

If $\f$ is a left Haar integral on an aqg $(A,\D)$, and $(H,\Lambda)$ is an associated GNS~pair,
then it can be shown that there is a unique $*$-homomorphism ${\pi\colon A\to B(H)}$ such that 
$\pi(a)\Lambda(b)=\Lambda(a b)$, for all ${a,b\in A}$.  Moreover, $\pi$ is faithful and
non-degenerate. We let $A_r$ denote the norm closure of $\pi(A)$ in $B(H)$. Thus, 
$A_r$ is a non-degenerate C*-subalgebra of $B(H)$.  There exists a unique non-degenerate
$*$-homomorphism ${\dar : A_r \rightarrow M(A_r\otimes A_r)}$ such
that, for all $a \in A$ and all $x \in A \otimes A$, we have 
${\dar (\pi(a)) (\pi \otimes \pi) (x)} = {(\pi \otimes \pi) (\Delta (a)x)}$ and ${(\pi \otimes \pi)
(x)\dar (\pi(a))} = {(\pi \otimes \pi) (x \Delta (a))}$.

Define a linear map
$\hat{\Lambda} :\hat{A}\rightarrow H$ by setting
$\hat{\Lambda}(\hat {a})=\Lambda (a)$, for all $a\in\2A$. Since
$\hat{\psi} (\hat{b}^*\hat{a})
=\varphi (b^*a)=(\Lambda(a),\Lambda(b))$, for all $a,b\in A$, it follows that
$(H,\hat{\Lambda})$ is a GNS-pair associated to $\hat{\psi}$.
It can be shown that it is unitarily equivalent to the
GNS-pair for a  left Haar integral $\hat{\varphi}$ of $(\hat
A,\hat \D)$. Hence, we can use $(H,\hat\Lambda)$ to define a
representation of the analytic extension $(\hat{A}_r, \hat \dar)$ of $(\hat A, \hat \D)$
on the space $H$. There is a unique 
$*$-homomorphism $\hat{\pi}_r :\hat{A}\rightarrow B(H)$ such
that  $\hat{\pi}_r (a )\hat{\Lambda}(b )=\hat{\Lambda}(a b)$,
for all $a,b\in\hat{A}$. Moreover, $\hat\pi$ is faithful and
non-degenerate. Let $\hat A_r$ be the norm closure
of $\hat{\pi}_r (A)$ in $B(H)$, so $\hat A_r$ is a
non-degenerate C*-subalgebra of $B(H)$. There exist a unique $*$-homomorphism 
$\hat{\D}_r: \hat{A}_r \rightarrow M(\hat{A}_r\otimes\hat{A}_r)$ such that, for
all $a \in \hat{A}$ and $x \in \hat{A} \otimes \hat{A}$,
$\hat{\D}_r(\hat\pi(a)) (\hat\pi \otimes \hat\pi) (x) =
(\hat\pi \otimes \hat\pi) (\hat{\D} (a)x)$ and
$(\hat\pi \otimes \hat\pi) (x)\hat \D_r (\hat\pi(a))=(\hat\pi\otimes\hat\pi)(x\hat{\D}(a))$.


\subsection{Discrete and Compact Case}

The following proposition is well known and easy to prove.

\begin{prop} \label{Berkeley}
Let $I$ be a set and $H_i, i\in I$ finite dimensional Hilbert spaces. Consider the
algebraic direct sum 
\[ A=\bigoplus_{i\in I} \ B(H_i) \]
\begin{enumerate} 
\item $A$ is a $*$-algebra with the usual $*$-operation on $B(H_i)$.
\item $A$ is unital iff $I$ is finite. 
\item The canonical embedding maps $\iota_i: B(H_i)\hookrightarrow A$ and projections 
$p_i: A\rarr B(H_i)$ are $*$-homomorphisms. 
\item Let $M(A)$ be the multiplier algebra and $L: A\rarr M(A)$ the canonical embedding map. 
There exists a unique isomorphism
\[ P: \ M(A)\rarr\prod_{i\in I} \ B(H_i) \]
such that $PL=\id$, where we use the obvious identification of the direct sum with a subalgebra of 
the direct product. 
\item If $\pi$ is an irreducible ($*$)-representation of $A$ there is $i\in I$ and an
invertible (unitary) $U\in B(H,H_i)$ unique up to a scalar such that 
\[ \pi(\cdot)=\pi\iota_i p_i(\cdot)=U p_i(\cdot) U^{-1}. \]
\item Any ($*$-)representation of $A$ is a direct sum of irreducible
($*$-)representations.
\end{enumerate}
\end{prop}

From now on the isomorphism $P$ will be used to make an identification 
$M(A)\equiv\prod_{i\in I} \ B(H_i)$ and therefore suppressed.

\bdefin \label{def-comp-disc}
A multiplier Hopf $*$-algebra $(A, \D)$ is called discrete if $A=\oplus_i B(H_i)$.
An aqg is called discrete if it is discrete as a multiplier Hopf $*$-algebra.
\edefin

The proof of the following result is included, since it is instructive for our purposes.

\btheor \label{t-disc2}
Let $(A,\Delta)$ be a discrete multiplier Hopf $*$-algebra. Then there is a unique
$f\in M(A)$ satisfying 
\begin{enumerate}
\item $f$ is positive and invertible,
\item $\DS S^2(a)=faf^{-1}\ \forall a\in A$,
\item $\DS Tr_i\, p_i(f)=Tr_i\, p_i(f^{-1})$, where $Tr_i$ is the usual trace on
$B(H_i)$. 
\end{enumerate}
It also satisfies $\DS S(f)=f^{-1}$ (where the unique extension of $S$ to $M(A)$ is
understood). 
\etheor
\prf To simplify the formulae we introduce the antilinear map $\chi(a)=S(a^*)$ which
satisfies $\chi^2=\id_A$. For every $i\in I$ we pick an antiunitary $J_i: H_i\rarr H_i$. 
Now, $J_i p_i(\cdot)J_i^{-1}$ is a (linear) irreducible representation of $A$ on
$H_i$ for any $i\in I$. By Prop.\ \ref{Berkeley} there is a unique $c(i)\in I$ and an
invertible $V_i: H_i\rarr H_{c(i)}$ such that
\begin{equation} \label{equ}
V_iJ_i \, p_i\chi(\cdot)J_i^{-1}V_i^{-1}= p_{c(i)}. 
\end{equation}
Now $U_i\equiv V_iJ_i$ is unique up to a scalar. We calculate
\[ U_ip_i(a)U_i^{-1}=U_ip_i(\chi\chi(a))U_i^{-1}=p_{c(i)}(\chi(a))=
   U_{c(i)}^{-1}p_{c(c(i))} U_{c(i)}. \]
In the first step we have used $\chi^2=\id$ and then we have used (\ref{equ}) twice. 
This shows that the irreducible representations $p_i$ and $p_{c(c(i))}$ are equivalent and
therefore $c$ is involutive: $c(c(i))=i$. From now on we write $\overline{\imath}=c(i)$. We have
\[ U_{\overline{\imath}}U_ip_i(\cdot)U_i^{-1}U_{\overline{\imath}}=p_i(\cdot), \]
and by irreducibility $U_{\overline{\imath}}U_i=k_i 1\in B(H_i)$. 
Using the freedom (up to a scalar) in the choice of the $U_i$ we 
can assume that $k_i=1\ \forall i\in I$. (If $i=\overline{\imath}$ we can absorb $k_i$ in
$U_i$. If $i\ne \overline{\imath}$, we can achieve $U_{\overline{\imath}}U_i=1$, and since the $U_i$
are invertible this implies that automatically $k_{\overline{\imath}}=1$.)

From (\ref{equ}) we obtain $p_iS(a)=U_i^{-1}p_\oli(a^*)U_i$ and applying this twice yields 
\[  p_i S^2(a) = U_i^{-1}p_\oli S(a)^* U_i 
   = U_i^{-1}[U_\oli^{-1}p_{\overline{\oli}}(a^*)U_\oli]^*U_i
   =  U_i^{-1}U_\oli^* p_i(a)U_\oli^{-1 *}U_i.
\]
Thus $p_i S^2(a)=G_i p_i(a) G_i^{-1}$ with $G_i=U_i^{-1}U_\oli^*=k_i^{-1}U_\oli U_\oli^*$ 
positive and invertible for all $i\in I$. Defining
\[ F_i=\sqrt{\frac{Tr_i\,G_i^{-1}}{Tr_i\,G_i}}\,G_i, \]
we have
\[ Tr_i\, F_i=\sqrt{Tr_i\,G_i^{-1}\,Tr_i\,G_i}=Tr_i\,F_i^{-1} \]
as required. Obviously, also the $F_i$ are positive and invertible and
$p_iS^2(a)=F_ip_i(a)F_i^{-1}$. The uniqueness is obvious in view of the normalizations.
It remains to prove $S(f)=f^{-1}$, which is clearly equivalent to
$p_iS(f)=p_i(f^{-1})\ \forall i$. As above, with (\ref{equ}) we have 
$p_iS(f)=U_i^{-1}p_\oli(f^*)U_i=U_i^{-1}F_\oli^*U_i$. On the other hand,
$p_i(f^{-1})=F_i^{-1}$. Thus we have to show
\[ U_i^{-1}F_\oli^*U_i = F_i^{-1}\quad\forall i\in I. \]
Inserting the definition of the $F_i$ this is equivalent to
\[ U_i^{-1}\sqrt{\frac{Tr_\oli\,G_\oli^{-1}}{Tr_\oli\,G_\oli}}U_iU_\oli^{-1*}U_i=
   \sqrt{\frac{Tr_i\,G_i}{Tr_i\,G_i^{-1}}} U_\oli^{*-1}U_i\quad\forall i\in I. \]
This is clearly equivalent to
\[ Tr_\oli\,G_\oli^{-1}\ Tr_i\,G_i^{-1}=Tr_i\,G_i\ Tr_\oli\, G_\oli\quad\forall i\in I. \]
Plugging in $G_i=U_i^{-1}U_\oli^*$ this is an easy computation which we omit. (One uses
$U_\oli=U_i^{-1}$ and cyclicity of the traces.)
\qed

\btheor \label{thm-haar}
A discrete multiplier Hopf $*$-algebra $(A,\D)$ admits a non-zero
left-invariant positive linear functional $\varphi$, thus is an aqg. 
\etheor

\bprop \label{p-disc}
An aqg $(A,\D)$ is discrete iff there exists $h\in A$ such that 
$ah=\ve(a) h$, for all $a\in A$. 
\eprop

\bdefin
An aqg $(A,\D)$ is called compact if $A$ is unital.
\edefin

Note that an aqg $(A,\Delta)$ is compact iff it is a Hopf $*$-algebra. Hence the analytical
extension of a compact aqg is a compact quantum group in the sense of S.L. Woronowicz. 
A multiplier Hopf $*$-algebra $(A,\D )$ with $A$ unital is not necessarily an aqg, but this is the
case whenever $A$ has a C*-algebra envelope. It is easily seen that $(A,\D )$ is compact and
discrete iff $A$ is finite dimensional.

\btheor 
Let $(A,\D )$ be an aqg. Then $(A,\D )$ is compact iff $(\hat{A},\hat{\D})$ is discrete.
\etheor


\subsection{Examples}

The examples discussed below are treated in depth in \cite{ES}. The translation into the framework
of \aqgs is straightforward and is recommended as an exercise.

Throughout this section let $\Gamma$ denote a discrete group with unit $e$. 
Consider the group algebra $\7C\Gamma$ and the algebra $C_c (\Gamma )$ of finitely
supported complex valued functions with the usual pointwise operations.
They are both non-degenerate $*$-algebras. The group algebra has unit $\delta_e$, whereas 
$C_c (\Gamma )$ has a unit iff $\Gamma$ is finite. 

In the group algebra case $(\7C\Gamma, \D)$ the comultiplication is given by
$\D(\delta_g )=\delta_g\otimes\delta_g$ and yields a Hopf $*$-algebra with counit 
$\varepsilon (\delta_g )=1$ and coinverse $S(\delta_g )=\delta_{g^{-1}}$, for all $g\in\Gamma$.
The functional defined by $\varphi (\delta_g )=\delta_{g,e}$, for $g\in\Gamma$, is unital, positive
and both left and right invariant, so $(\7C\Gamma, \D)$ is a cocommutative compact aqg.  

Note that $M(C_c (\Gamma ))\cong C(\Gamma )$, where $C(\Gamma )$ is the unital $*$-algebra of all
functions on $\Gamma$, and that $C_c (\Gamma\times\Gamma )\cong C_c (\Gamma )\otimes C_c (\Gamma )$
(as $*$-algebras). Thus the algebra $C_c (\Gamma )$ has comultiplication 
$\D:C_c (\Gamma )\rightarrow M(C_c (\Gamma )\otimes C_c (\Gamma ))$
given by $\D(a)(g,h)=a(gh)$ for $a\in C_c (\Gamma )$ and $g,h\in\Gamma$. It is easy to see that the
support of $\D (a)$ for $a\in C_c (\Gamma )$ is infinite whenever $\Gamma$ is. This shows that the
Hopf algebra framework is too restrictive to cover this example. 
For $(C_c(\Gamma), \D )$ we have $\ve(a)=a(e)$ and $S(a)(g)=a(g^{-1})$
where $a\in C_c(\Gamma\times\Gamma)$ and $g\in\Gamma$.
The integral $\varphi(a)=\sum_{g\in\Gamma} a(g)$, $a\in\7C\Gamma$, w.r.t. the counting measure is
positive and both left and right invariant, so $(C_c (\Gamma ), \D)$ is a commutative discrete aqg.  

The following proposition shows that these two examples exhaust the cocommutative compact \aqgs\ and 
the commutative discrete \aqgs.

\bprop 
Let $(A,\D)$ be an aqg. Then
\begin{enumerate}
\item $(A,\D)$ is discrete and commutative iff there exists a discrete group
$\Gamma$ such that $(A,\D)\cong (C_c (\Gamma),\D)$.
\item $(A,\D)$ is compact and co-commutative iff there exists a discrete group
$\Gamma$ such that $(A,\D)\cong (\7C\Gamma, \D )$.
\end{enumerate}
\eprop

\noindent It is also easy to check that $(\hat{\7C\Gamma },\hat{\D})\cong (C_c (\Gamma),\D)$ and 
$(\hat{C_c (\Gamma)},\hat{\D})\cong (\7C\Gamma, \D )$, so they are dual to each other.

Giving a characterization of the commutative compact and the cocommutative discrete algebraic
quantum groups on the algebraic level, is less immediate and requires some (co-) representation
theory. For the moment we mention that the commutative-compact case can be described in terms of the
algebra of regular functions on a compact group, and we make here an analytic statement concerning
both cases.

\bprop
Let $(A,\D)$ be an aqg with analytic extension $(A_r, \D_r )$. Then
\begin{enumerate}
\item $(A,\D)$ is compact and commutative iff there exists a compact group
$G$ such that $(A_r,\D_r )\cong (C(G),\D_r )$.
\item $(A,\D)$ is discrete and cocommutative iff there exists a compact group
$G$ such that $(A_r, \D_r )\cong (C_r^* (G),\D_r )$.
\end{enumerate}
\eprop

Again these two cases are dual to each other. Here $C(G)$ denotes the C*-algebra of continuous
functions on $G$ with uniform norm, whereas $C_r^* (G)$ denotes the reduced group C*-algebra of
$G$. Since a compact group $G$ is amenable, $C_r^* (G)$ coincides with the universal group
C*-algebra $C^* (G)$, and therefore, as we shall see very explicitly in the next section, the
representation theory of a discrete cocommutative aqg coincides with that of $G$. We shall also give
a simple algebraic description of the associated aqg in this case.


\section{Representation Theory} \label{S-repres}
\subsection{The Universal Corepresentation}
Throughout this section $(A,\D)$ stands for an arbitrary aqg and $B$ denotes a $*$-algebra.
We recall the universal corepresentation $U$ of $(A,\D)$ and discuss its various properties \cite{ku-e}.

\begin{defin} A corepresentation of $(A,\D )$ on $B$ is an
element $V\in M(A\otimes B)$ such that $(\D\otimes \iota )V=V_{13}V_{23}$. 
\end{defin}

\begin{prop}
Suppose $V$ is a corepresentation of $(A,\D)$ on $B$. Then the following are equivalent:
\begin{enumerate}
\item $V$ is invertible in $M(A\otimes B)$ with $V^{-1}=(S\otimes\iota )V$.
\item $(\varepsilon\otimes\iota )V=1$. 
\item $V(A\otimes B)=A\otimes B =(A\otimes B)V$.
\end{enumerate}
A corepresentation $V$ satisfying these equivalent conditions is said to be non-degenerate. 
\end{prop}

\begin{defin} A unitary corepresentation $V$ of $(A,\D )$ on $B$ is
a corepresentation on $B$ which is a unitary element of $M(A\otimes B)$. 
\end{defin}

The following result is fundamental for corepresentations of \aqgs.
We regard elements of $A\otimes\hat{A}$ as endomorphisms of $A$.

\begin{thm} \label{t-univ-corep}
There exists a unique element $U\in M(A\otimes\hat{A})$ such that 
\[ [U(x\otimes\omega )](y)=(\iota\otimes\omega )(\D (y)(x\otimes 1)) \]
\[ [(x\otimes\omega )U](y)=(\omega\otimes\iota )((1\otimes x)\D (y)) \]
for all $x,y\in A$ and $\omega\in\hat{A}$. 
Moreover, we have the following properties for $U$:
\begin{enumerate}
\item $U$ is unitary in $M(A\otimes\hat{A})$.
\item $(\D\otimes\iota )U=U_{13}U_{23}$.
\item $(\iota\otimes\hat{\D})U=U_{12}U_{13}$.
\item $(\omega\otimes\iota )U=\omega$ for all $\omega\in\hat{A}$.
\item $(\iota\otimes a)U=a$ for all $a\in A$, where $a$ acts on $\hat{A}$ by the identification $A=\ahh$.
\end{enumerate}
\end{thm}

Note that claims 1 and 2 say that $U$ is a unitary corepresentation of $(A,\D )$ on $\hat{A}$, 
whereas claims 1 and 3 say that $\sigma(U)$ is a unitary corepresentation of
$(\hat{A},\hat{\D})$ on $A$. Below we will explain why $U$ is called the universal corepresentation
of $(A,\D)$. It follows that $\sigma (U)$ is the universal corepresentation of $(\hat{A},\hat{\D})$.  

Let $\Hom(A,B)$ denote the set of $*$-homomorphisms from $A$ to $M(B)$ satisfying
$\theta (A)B=B$.

\begin{prop}
\label{p-functP}
Let $V\in M(A\otimes B)$ and consider the linear map
\[ \pi_V:\ \hat{A}\rightarrow M(B), \quad \omega \mapsto (\omega\otimes\iota) V. \] 
Then the following equivalences hold:
\begin{enumerate}
\item $V$ is a corepresentation (of $(A,\D )$ on $B$) iff $\pi_V$ is multiplicative. 
\item $V$ is a non-degenerate corepresentation iff $\pi_V$ is multiplicative and non-degenerate.
\item $V$ is a unitary corepresentation iff $\pi_V\in\Hom(\hat{A},B)$. 
\end{enumerate}
\end{prop}

The universality of $U$ can now be formulated as follows.

\begin{thm}
\label{t-univ-corr}
For any unitary corepresentation $V$ of $(A,\D)$ on $B$ and any
$\theta\in\Hom(\hat{A},B)$, we have:
\begin{enumerate}
\item $(\iota\otimes \pi_V )U=V$.
\item $\pi_{(\iota\otimes\theta )U}=\theta$. 
\end{enumerate}
\end{thm}

The result remains valid for the non-degenerate situation without $*$-operation.

\begin{defin} 
\label{d-tensor}
Let $B'$ be non-degenerate $*$-algebra. For $\theta\in\Hom(A,B)$ and 
$\theta '\in \Hom(A,B')$ we define $\theta\times\theta'\in \Hom(A,B\otimes B' )$
by $\theta\times\theta' =(\theta\otimes\theta' )\D$. If $V$ and $V'$ are unitary corepresentations of 
$(A,\D)$ on $B$ and $B'$, respectively, we define a unitary corepresentation $V\times V'$ of $(A,\D)$ on 
$B\otimes B'$ by $V\times V' =V_{12}V'_{13}$.  
\end{defin}

\begin{thm} 
\label{t-bij}
The map sending unitary corepresentations $V$ of $(A,\D)$ on $B$ to $\pi_V\in\Hom(\hat{A},B)$,
is a bijection. Let $V'$ be a unitary corepresentation of $(A,\D)$ on $B'$. Then
$\pi_{V\times V'} =\pi_V\times\pi_{V'}$. For the universal corepresentation $U$ of $(A,\D)$ on
$\hat{A}$ we have $\pi_U =\iota\in\Hom(\hat{A},\hat{A})$, and finally for the trivial unitary
corepresentation $1\otimes 1\in A\otimes B(\7C )$ of $(A,\D)$ we have 
$\pi_{1\otimes 1}=\hat{\varepsilon}$, where  
$\hat{\varepsilon}$ is the counit of $(\hat{A},\hat{\D})$.  
\end{thm}

\prf The only part which is not immediate from the results above, is the 
identity $\pi_{V\times V'} =\pi_V\times\pi_{V'}$. To show this, first note that
\[ V\times V'=V_{12}V'_{13}=[(\iota\otimes\pi_V )U]_{12}[(\iota\otimes\pi_{V'} )U]_{13}  
=(\iota\otimes\pi_V\otimes\pi_{V'})(U_{12}U_{13}) \]
\[ =(\iota\otimes\pi_V\otimes\pi_{V'})(\iota\otimes\hat{\D})U
=(\iota\otimes(\pi_V\otimes\pi_{V'})\hat{\D})U
=(\iota\otimes\pi_{V}\times\pi_{V'})U. \] 
Thus we may conclude that
\[ \pi_{V\times V'}(\omega )=(\omega\otimes\iota\otimes\iota )(V\times V' )
=(\omega\otimes\iota\otimes\iota )(\iota\otimes\pi_{V}\times\pi_{V'})U
=(\pi_{V}\times\pi_{V'})(\omega\otimes\iota )U =(\pi_{V}\times\pi_{V'})(\omega ), \]
for all $\omega\in\hat{A}$. 
\qed

\brem \label{r-cat}
The correspondence established in Theorem \ref{t-bij} is in fact a functor if we define arrows for
representations and corepresentations on $*$-algebras in the following way. For any objects
$\pi\in\Hom(A,B)$ and $\pi'\in \Hom(A,B')$ define an arrow 
$f:\pi\rightarrow\pi'$ to be $f\in \Hom(B,B')$ such that $f\pi =\pi'$. 
Similarly, if $V$ and $V'$ are unitary corepresentations of $(A,\D)$ on $B$ and $B'$, respectively, 
we define an arrow $f:V\rightarrow V'$ between these objects to be $f\in \Hom(B,B')$ such
that $V' =(\iota\otimes f)V$. Then the correspondence $V\mapsto\pi_V$ in Theorem \ref{t-bij} is an
equivalence of tensor categories (they are not $*$-categories). 
\erem

\brem \label{r-actcoact}
We remark on the tensor categories of actions and coactions of an aqg $(A,\D )$ on $*$-algebras.  

By an action $\gamma$ of $(A,\D)$ on a $*$-algebra $B$, we mean a surjective linear map 
$\gamma :A\otimes B\rarr B$ such that $\gamma (m\otimes\iota )=\gamma (\iota\otimes\gamma )$. An
arrow $f:\gamma\rightarrow \gamma'$ is an element $f\in\Hom(B,B')$ such that 
$f\gamma=\gamma' f$. We do not define tensor products in the general case. 

Given $\pi\in\Hom(A,B)$, we define the action $\gamma_{\pi}$ of $(A,\D)$ on $B$ by
$\gamma_{\pi}(a\otimes b)=\pi (a)b$, for all $a\in A$ and $b\in B$. Using non-degeneracy of $f$, we
see that for any $f\in \Hom(B,B')$, we have $f:\pi\rarr \pi'$ iff $f:\gamma_{\pi}\rarr\gamma_{\pi'}$.   
Thus $\pi\mapsto\gamma_{\pi}$ is an equivalence of the tensor category of non-degenerate
$*$-representations of $(A,\D)$ on $*$-algebras to the tensor category of actions of $(A,\D )$ on
$*$-algebras contained in the image of  
the functor and with tensor product $\gamma_{\pi}\times \gamma_{\pi '}$ simply defined to be the
action $\gamma_{\pi\times\pi'}$. We have not incorporated the $*$-preserving property of $\pi$
in $\gamma_{\pi}$ and the axioms for actions.  

A coaction $\delta$ of $(A,\D )$ on a $*$-algebra $B$ is a $\delta\in\Hom(B,A\otimes B)$ 
such that $(\iota\otimes\delta)\delta =(\D\otimes\iota)\delta$ and 
$(\varepsilon\otimes\iota)\delta=\iota$. 
Such coactions form a tensor category, an arrow $f:\delta\rightarrow\delta '$ being an element
$f\in\Hom(B,B')$ such that $(\iota\otimes f)\delta =\delta 'f$ and the 
tensor product $\delta\times\delta '$ being the comodule of $(A,\D )$ on $B\otimes B'$ 
given by $(\delta\times\delta ')(b\otimes b')=\delta (b)_{12}\delta' (b')_{13}$, for all
$b\in B$ and $b'\in B'$. 

Let $V$ be a unitary corepresentation of $(A,\D )$ on $B$. Then the linear map 
$\delta_V :B\rightarrow M(A\otimes B)$ given by $\delta_V (b)=V (1\otimes b)V^*$ for all $b\in B$,
is indeed a coaction of $(A,\D^{op})$ on $B$. Non-degeneracy of $\delta_V$ follows from 
$V(A\otimes B)=A\otimes B$, and the formula
$(\iota\otimes\delta_V)\delta_V =(\D^{op}\otimes\iota)\delta_V$ follows from the calculation
\begin{eqnarray*} \lefteqn{  (\iota\otimes\delta_V )\delta_V (b)
   =(\iota\otimes\delta_V )V(1\otimes\delta_V (b))((\iota\otimes\delta_V )V)^*
   =V_{23}V_{13}V_{23}^* (1\otimes V(1\otimes b)V^* )(V_{23}V_{13}V_{23}^* )^* } \\
  && =V_{23}V_{13}(1\otimes 1\otimes b)V_{13}^* V_{23}^*
=(\D^{op}\otimes \iota )V(1\otimes 1\otimes b)((\D^{op}\otimes\iota )V)^*
=(\D^{op}\otimes\iota )(V(1\otimes b)V^* ) \\
 &&=(\D^{op}\otimes\iota )\delta_V (b)
\end{eqnarray*}
for all $b\in B$. Note that any arrow $f:V\rightarrow V'$ will be an arrow 
$f:\delta_V\rightarrow \delta_{V'}$, but the converse is not true (consider $B'$ commutative). So we
have a tensor functor $V\mapsto\delta_V$ from the tensor category of unitary corepresentations of
$(A,\D)$ on $*$-algebras to the tensor category of coactions of $(A,\D^{op})$ on $*$-algebras, which
in general, is not an equivalence.
\erem


\subsection{Representations vs.\ Corepresentations as Tensor $*$-Categories}

We now proceed to establish the correspondence between representations of an aqg $(A,\D)$ and
corepresentations of $(\hat{A},\hat{\D})$, restricting ourselves to finite dimensional Hilbert
spaces. For homomorphisms $\pi: A\rarr\End K$ there are two different notions of non-degeneracy
which, fortunately, are equivalent for $K$ finite dimensional.

\blemma \label{l-non-deg}
Let $A$ be an algebra and $K$ a finite dimensional vector space. For a homomorphism 
$\pi: A\rarr\End\,K$ the conditions $\pi(A)\End\,K=\End\,K$ and $\pi(A)K=K$ are equivalent. 
\elemma
\prf For $K$ finite dimensional we clearly have $(\End\,K)K=K$. Assuming $\pi(A)\End K=\End K$ we
compute 
\[ \pi(A)K=\pi(A)((\End K)K)=(\pi(A)\End K)K=(\End K)K=K. \]
Conversely, assume $\pi(A)K=K$. With the isomorphism 
$\alpha: K\otimes K^*\rarr \End\,K$ determined by $\alpha(x\otimes\phi)(z)=x\phi(z)$
for $x,z\in K, \phi\in K^*$ we have
$(s\alpha(x\otimes\phi))(z)=sx\phi(z)=\alpha(sx\otimes\phi)(z)$ and therefore
\[ \pi(A)\End\,K=\pi(A)\alpha(K\otimes K^*)=\alpha((\pi(A)K)\otimes K^*)
   =\alpha(K\otimes K^*)=\End\,K, \]
as desired.
\qed

\begin{defin}
\label{d-repf} Let $(A,\D )$ be an aqg.
Let $\Rep_f (A,\D )$ denote the class of (algebraically) non-degenerate $*$-representations
of $A$ on finite dimensional Hilbert spaces including the zero representation $0$. 
Let $\pi_1, \pi_2\in \Rep_f (A,\D )$.
Then $\pi_i\in\Hom(A,B(K_i ))$, so we define
$\pi_1\times\pi_2 \in\Hom(A,B(K_1 )\otimes B(K_2 ))
=\Hom(A,B(K_1\otimes K_2 ))\subset\Rep_f (A,\D )$ 
according to Definition~\ref{d-tensor}.
\end{defin}

We regard $\Rep_f (A,\D )$ as a tensor $*$-category with representations as objects and 
intertwiners as arrows. (Recall that a zero representation $0$ is not regarded as an irreducible 
object.) The tensor product is the one given in the definition above.
Clearly, $\Rep_f (A,\D )$ is a tensor $*$-category with the counit $\varepsilon$ as 
the irreducible unit. Of course, the $*$-operation in the category is the usual Hilbert space adjoint
of operators.

\begin{defin}
\label{d-corepf}
We say $V$ is a finite dimensional unitary corepresentation of an aqg $(A,\D )$ on a Hilbert space
$K$ if $K$ is finite dimensional and $V$ is a unitary corepresentation of $(A,\D )$ on the unital
$*$-algebra $B(K)$ of linear operators on $K$. Let $\Corep_f (A,\D)$ denote the class of
all such unitary corepresentations including the zero corepresentation $0$. Let 
$V, V'\in \Corep_f (A,\D )$. Then we define the unitary corepresentation $V\times V'$ of
$(A,\D )$ on $B(K)\otimes B(K')=B(K\otimes K' )$ according to Definition~\ref{d-tensor}. Clearly, we
have $V\times V'\in \Corep_f (A,\D )$. 
\end{defin}
 
We regard $\Corep_f (A,\D )$ as a tensor $*$-category with corepresentations as objects and
intertwiners as arrows. Recall that $T\in B(K,K' )$ is an intertwiner between two corepresentations
$V,V'\in \Corep_f (A,\D )$ iff $T(\omega\otimes\iota )V =(\omega\otimes\iota )V' T$, for all
$\omega\in\hat{A}$. The tensor product is the one given in the definition above.
Clearly, $\Corep_f (A,\D )$ is a tensor $*$-category with the unit $1\otimes 1$ as the irreducible
unit. Again, the $*$-operation in the category is the usual Hilbert space adjoint of operators.   

\medskip

The two tensor categories are related in the following way. 

\begin{thm} \label{t-dualcor}
Let notation be as above. The correspondence
\[ P: \ \Corep_f (A,\D) \rarr \Rep_f (\hat{A},\hat{\D} ), \quad V \mapsto \pi_V \]
provided by Proposition \ref{p-functP} together with the identity map on morphisms 
gives rise to an isomorphism of tensor $*$-categories.
\end{thm}

\prf The proof is immediate from Theorem \ref{t-bij}.
\qed

\begin{rema}
A similar result can be obtained for the infinite dimensional case. One needs then to talk about
$*$-representations which are non-degenerate in the $C^*$-algebra sense, so
$\pi\in \Rep (A,\D )$ iff the exists a (possibly infinite dimensional) Hilbert space $K$ and a 
$*$-homomorphism $\pi :A\rightarrow B(K)$ such that the vector space $\pi (A)K$ is dense in $K$.  
The tensor $*$-category $\Corep (A_r, \D_r )$ of unitary corepresentations on Hilbert spaces
consists of Hilbert spaces $K$ and unitaries $V\in M(A_r\otimes B_0 (K))$ such that 
$(\D_r\otimes\iota )V =V_{13}V_{23}$.
For the exact correspondence thus established, see \cite{BMT}.
\end{rema}

We will look at another way of obtaining a tensor $*$-category $\Rep_f (A^o_s, \hat{\D})$
from a compact aqg $(A,\D)$ which is equivalent to $\Corep_f (A,\D )$.
This can sometimes come in very handy, especially when dealing with representations of quantized 
universal enveloping Lie algebras $U_q (\6g)$.  

First recall that $(A,\D )$ is a Hopf $*$-algebra with counit $\varepsilon$, coinverse $S$ and unit
$I$, so the vector space $A'$ of all linear functionals on $A$ is a unital $*$-algebra with unit
$\varepsilon$, product $\omega\eta =(\omega\otimes\eta)\D\in A'$ and $*$-operation $\omega^*\in A'$ 
given by $\omega^* (a)=\overline{\omega (S(a)^* )}$, for all $a\in A$. Define a linear map $\hat{\D}
:A'\rightarrow (A\otimes A)'$ by $\hat{\D}(\omega )(a\otimes b)=\omega (ab)$, 
for all $\omega\in A'$ and $a,b\in A$. Consider the subspace $A^o$ of $A'$ given by 
\[A^o =\{\omega\in A'\ |\ \hat{\D}(\omega )\in A'\otimes A'\},\]
where we understand the embedding $A'\otimes A'\subset (A\otimes A)'$. Then
$(A^o, \hat{\D} )$ is called the {\it maximal dual} Hopf $*$-algebra of $(A,\D )$. 
Note that for $\hat{f}\in M(\hat{A})$ as in Theorem \ref{t-disc2}, we have $\hat{f}\in A^o$ and 
$\hat{\D}(\hat{f })=\hat{f}\otimes\hat{f}$.
Regard $A^o$ as a locally convex topological vector space with pointwise convergence,
so $\omega_{\lambda}\rightarrow\omega$ in $A^o$ if $\omega_{\lambda}(a)\rightarrow\omega (a)$,
for all $a\in A$. The continuous unital $*$-representations 
of $(A^o, \hat{\D} )$ on finite dimensional Hilbert spaces
clearly form a tensor $*$-category with conjugates given by the same formulas as for 
$\Rep_f (\hat{A},\hat{\D})$. Note that in general, the maximal dual Hopf $*$-algebra can be very small!

\btheor Let $(A,\D)$ be a compact aqg and let $(A^o_s, \hat{\D})$ be a sub Hopf $*$-algebra
of the maximal dual Hopf $*$-algebra of $(A,\D )$ which separates the elements of $A$.
(Thus $A^o_s$ is a unital $*$-subalgebra of $A^o$ and $\hat{\D}$ is given by restriction). 
Consider the tensor $*$-category $\Rep_f (A^o_s, \hat{\D})$ of continuous unital $*$-representations
of $(A^o_s, \hat{\D} )$ on finite dimensional Hilbert spaces with arrows, $*$-operation and tensor
product defined in the obvious way. Then the assignment
\[N: \ \Corep_f (A,\D )\rarr \Rep_f (A^o_s, \hat{\D}), \quad V\mapsto \pi_V \ \mbox{with} \ 
   \pi_V (\omega )=(\omega\otimes\iota )V, \ \forall\omega\in A^o_s, \]
together with the identity map on arrows, is an equivalence of tensor $*$-categories.
\etheor

\prf 
The only part which requires a proof is that any $\pi\in \Rep_f (A^o_s, \hat{\D})$ is of the form 
$\pi_{V_{\pi}}$ for some $V_{\pi}\in \Corep_f (A,\D )$. Thus let 
$\pi :A^o_s\rightarrow B(K)$ be a continuous unital $*$-representation of $A^o_s$ on a finite
dimensional Hilbert space $K$. We will provide a $V_{\pi}\in A\otimes B(K)$ such that
$\pi_{V_{\pi}}=\pi$. 
Using the fact that $A^o_s$ separates the elements of $A$, it is then clear that $V_{\pi}$ will be
unitary, that $(\D\otimes\iota )V_{\pi} =(V_{\pi})_{13}(V_{\pi})_{23}$ and that the $*$-preserving
tensor functor $N$ is an equivalence. 

Now pick an orthonormal basis $(e_i )$ for $K$ and consider the system $m_{ij}$ of matrix units for
$B(K)$ defined as $m_{ij}(e_k )=\delta_{jk}e_i$, for all $i,j,k$. Consider, for fixed $i$ and $j$,
the linear functional on $A^o_s$ given by $\omega\mapsto (\pi (\omega )e_j, e_i )$ for all
$\omega\in A^o_s$. Since $\pi$ is continuous w.r.t pointwise topology on $A^o_s$, it follows from  
\cite[Proposition 2.4.4]{ped}, that there exists a unique $a_{ij}\in A$ such that 
$\omega (a_{ij})=(\pi (\omega )e_j, e_i )$, for all $\omega\in A^o_s$.
Define now $V_{\pi}\in A\otimes B(K)$ by $V_{\pi}=\sum_{ij} a_{ij}\otimes m_{ij}$.
By construction $\pi_{V_{\pi}}=\pi$ and we are done.
\qed

\brem
It is known that for the usual quantized compact \aqgs\ $(A_q, \D )$, the associated quantized
universal enveloping Lie algebras $U_q (\6g )$ are Hopf $*$-algebras with $A^o_s =U_q (\6g )$ which
separate the elements of $A_q$ \cite[Sect.\ 7.1.5]{schm}. Thus the tensor $*$-category of continuous
unital $*$-representations of $(U_q (\6g ),\hat{\D})$ on finite dimensional Hilbert spaces is
equivalent to $\Corep_f (A_q, \D )$. Also $\hat{f}\in U_q (\6g)$, so the intrinsic dimension can be
read off conveniently for such representations \cite{RT}. We return to this issue in the next
subsection. 
\erem

We now recast our results in the language of modules and comodules. 

\bprop
Any non-degenerate $*$-representation $\pi$ of a discrete aqg $(A,\D )$ on a finite dimensional
Hilbert space $K$ gives rise to an $A$-module $K$. Namely, define a linear map
$\alpha_\pi: A\otimes K\rightarrow K$ by $\alpha_\pi (a\otimes \xi )=\pi (a)\xi$,
for all $a\in A$ and $\xi\in K$. Then the following properties hold:
\begin{enumerate}
\item $\alpha_{\pi}(m\otimes\iota )=\alpha_{\pi}(\iota\otimes\alpha_{\pi} )$.
\item $\alpha_{\pi} (A\otimes K)=K$.
\item $(\alpha_{\pi}(a^*\otimes\xi ),\eta )=(\xi, \alpha_{\pi}(a\otimes \eta)) \ \ 
\forall a\in A\ \ \forall \xi, \eta\in K$.
\end{enumerate}
Denote by $\Mod_f (A,\D )$ the tensor $*$-category of linear maps 
$\alpha :A\otimes K\rightarrow K$ satisfying properties 1-3 above, with arrows 
$t:\alpha\rightarrow\alpha '$ being linear maps $t:K\rightarrow K'$ satisfying 
$\alpha' =t\alpha$ and tensor product $\alpha\times\alpha ':A\otimes K\otimes K'\rightarrow K\otimes K'$ 
uniquely determined by 
\[(\alpha\times\alpha ')(\iota\otimes\alpha\otimes\alpha )(a\otimes b\otimes\xi\otimes c\otimes\xi')
   =\sum_k \alpha (a_k \otimes\xi )\otimes \alpha' (b_k\otimes \xi '), \]
for all $a,b\in A$, $\xi\in K$ and $\xi'\in K'$ and where $\sum_k a_k\otimes b_k =\D(a)(b\otimes c)$. 
Then the assignment $\pi\mapsto\alpha_{\pi}$ is a monoidal $*$-preserving equivalence of 
$\Rep_f (A,\D )$ and $\Mod_f (A,\D )$.  
\eprop

\prf
Straightforward, once one shows that $\alpha\times\alpha '$ is well-defined. This is most easily 
done by considering the corresponding well-defined tensor product $\pi\times\pi'$ 
of the associated non-degenerate $*$-representations, i.e, of $\pi$ and $\pi'$ such that
$\alpha =\alpha_{\pi}$ and $\alpha'=\alpha_{\pi'}$, and then
noting that $\alpha\times\alpha '=\alpha_{\pi\times\pi'}$.
\qed

\bprop
Any unitary corepresentation $V$ of a compact aqg $(A,\D)$ on a finite dimensional
Hilbert space $K$ gives rise to an $A$-comodule $K$. Namely, define a linear map
$\beta_V :K\rarr A\otimes K$ by $\beta_V (\xi )=\sum_k a_k\otimes x_k\xi$,
for all $\xi\in K$ and where $V=\sum_k a_k\otimes x_k\in A\otimes B(K)$. 
Then the following properties hold:
\begin{enumerate}
\item $(\D^{op}\otimes\iota )\beta_V =(\iota\otimes\beta_V )\beta_V$.
\item $(S\otimes(\cdot,\eta))\beta_V (\xi )=(*\otimes (\xi,\cdot ))\beta_V (\eta )\ \ \ \ \forall \xi, \eta\in K$,
\end{enumerate}
where $\varepsilon$ and $S$ are the counit and coinverse of $(A,\D)$, respectively. Let
$\Comod_f (A,\D^{op})$ denote the tensor $*$-category of linear maps 
$\beta :K\rarr A\otimes K$ satisfying properties 1-2 above, with arrows $t:\beta\rightarrow\beta'$
being linear maps $t:K\rightarrow K'$ satisfying $\beta' t=(\iota\otimes t)\beta$ and tensor product
$\beta\times\beta ':K\otimes K'\rightarrow A\otimes K\otimes K'$ given by 
\[(\beta\times\beta ')(\xi\otimes \xi')
=\beta (\xi )_{12}\beta'(\xi ')_{13},\]
for all $\xi\in K$ and $\xi'\in K'$.
Then the assignment $V\mapsto\beta_V$ is a monoidal $*$-preserving equivalence of 
$\Corep_f (A,\D )$ and $\Comod_f (A,\D^{op})$.  
\eprop

\prf
Straightforward. How to produce a unitary corepresentation $V\in A\otimes B(K)$ of $(A,\D)$ from a
prescribed comodule $\beta :K\rarr A\otimes K$ of $(A,\D^{op})$ perhaps calls for some explanation. 
Pick an orthonormal basis $(e_k )$ for $K$ and consider the system $(m_{ij})$ of matrix units of
$B(K)$ given by $m_{ij}(e_k )=\delta_{jk}\ e_i$ for all $i,j,k$. Define $V\in A\otimes B(K)$ by
$V=\sum_{ij} a_{ij}\otimes m_{ij}$, where the elements $a_{ij}\in A$ are given by
\[a_{ij}=(\iota\otimes (\cdot,e_i ))\beta (e_j ),\]
for all $i,j$. It is now easily checked that $\beta_V =\beta$. 
\qed

We thus get the following restatement of Theorem \ref{t-dualcor}.

\btheor \label{t-dual-comod}
Let $(A,\D)$ be a discrete aqg. Then $\Mod_f (A,\D )$ and $\Comod_f (\hat{A},\hat{\D}^{op})$ 
are equivalent as semisimple tensor $*$-categories. 
\etheor


\subsection{Conjugates for Representations and Corepresentations}
Suppose that $(A,\D )$ is a discrete aqg with counit $\varepsilon$
and coinverse $S$. We now show that the semisimple tensor $*$-categories $\Rep_f (A,\D )$ and 
$\Corep_f (\hat{A},\hat{\D} )$ have conjugates, and provide a conjugate
object for every non-zero object. 

Suppose $\pi$ is a non-zero non-degenerate $*$-representation on the finite dimensional Hilbert
space $K$. Let $f\in M(A)$ be as in Theorem \ref{t-disc2} and consider the positive operator 
$\pi(f^{-1})\in B(K)$. Pick any Hilbert space $\overline{K}$ and any invertible antilinear operator
$J:K\rarr\overline{K}$ such that $J^* J=\pi (f^{-1})$. This can clearly be done. Now pick an
orthonormal basis $(e_i )$ for $K$ and define linear maps $r: \7C\rightarrow \overline{K}\otimes K$
and $\overline{r} :\7C\rightarrow K\otimes\overline{K}$ by 
\[ r(1)=\sum_i Je_i\otimes e_i \quad \mbox{and} \quad  \overline{r}(1)=\sum_i e_i\otimes J^{*-1}e_i,  \]
respectively. Next define 
\[ \overline{\pi}:\ A\rarr B(\overline{K}), \quad a \mapsto J\pi S(a^* )J^{-1}. \]
As $J$ and $J^{-1}$ both are antilinear, $\overline{\pi}(a)$ is linear and thus
$\overline{\pi}(a)\in B(\overline{K})$ for all $a\in A$. By the following proposition,
$\overline{\pi}$ is indeed a conjugate to $\pi$ and $r$ and $\overline{r}$ are normalized solutions
of the conjugate equations for $\pi$ and $\overline{\pi}$. 

\bprop
\label{p-conj}
Let $(A,\D )$ be a discrete aqg and $\pi$ a non-zero object of $\Rep_f (A,\D )$. 
Then $\overline{\pi}$ is a non-zero object of $\Rep_f (A,\D )$. 
Moreover, the following hold:
\begin{enumerate} 
\item $r\in\Hom (\varepsilon, \overline{\pi}\times\pi )$,
\item $\overline{r}\in\Hom (\varepsilon, \pi\times\overline{\pi})$,
\item $\overline{r}^*\otimes\id_{\pi}\mcirc\id_{\pi}\otimes r=\id_{\pi}$,
\item $r^*\otimes\id_{\overline{\pi}}\mcirc\id_{\overline{\pi}}\otimes
    \overline{r}=\id_{\overline{\pi}}$, 
\item $r^* \mcirc  r=\overline{r}^*\mcirc \overline{r}$.
\end{enumerate}
\eprop

\prf Claims 3 and 4 hold for any invertible antilinear map $J$, as is easily verified.
Claim 5 is simply a restatement of the fact $\mathrm{Tr} \pi(f)=\mathrm{Tr} \pi(f^{-1})$ 
stated in Theorem \ref{t-disc2}. To show that $\overline{\pi}\in \Rep_f (A,\D )$,
we first note that $\overline{\pi}$ is linear (as $Ad\,J$ and $*$ are both antilinear), multiplicative
and non-degenerate, thus non-zero. To see that $\overline{\pi}$ is $*$-preserving, first
notice that $S(a^*)^* =S^{-1}(a)$ and $S^2 (a)=faf^{-1}$ for all $a\in A$, and then calculate
\[ \overline{\pi}(a)^* =J^{*-1}\pi S^{-1}(a)J^* =J\pi (f)\pi S^{-1}(a)\pi (f^{-1})J^{-1}
=J\pi S(a)J^{-1}=\overline{\pi}(a^* ),\]
for all $a\in A$.

It remains to show relations 1 and 2. We prove only the first, the second being proved
similarly. Now $r\in\Hom (\varepsilon, \overline{\pi}\times\pi )$ simply 
means that $\varepsilon (a)r(1)=(\overline{\pi}\times\pi )(a)r(1)$, for all $a\in A$.
By the non-degeneracy of $\pi$ it thus suffices to show that
\[ (\varepsilon (a)r(1),J^{*-1}e_j\otimes \pi (b^* )e_l )
=((\overline{\pi}\times\pi )(a)r(1),J^{*-1}e_j\otimes \pi (b^* )e_l ) \]
for all $a,b\in A$ and all $j,l$. On the l.h.s.\ we have
\[ (\varepsilon (a)r(1),J^{*-1}e_j\otimes \pi (b^* )e_l )=\varepsilon (a)\sum_i (Je_i, J^{*-1}e_j)
   (e_i, \pi (b^* )e_l )\] 
\[=\varepsilon (a)\sum_i \delta_{ij}(e_i, \pi (b^* )e_l )
=\varepsilon (a)(e_j, \pi (b^* )e_l )
=(\pi (\varepsilon (a)b )e_j, e_l ).\]
To see that the r.h.s.\ coincides with this expression, first write 
$(1\otimes b)\D (a)=\sum_k a_k\otimes b_k$ and notice that 
\[ \sum_k b_k S(a_k^* )^* =\sum_k b_k S^{-1}(a_k )=m(\iota\otimes S^{-1})(\sum_k b_k\otimes a_k )
=m(\iota\otimes S^{-1})((b\otimes 1)\D^{op} (a))=\varepsilon (a)b.\]
Hence 
\begin{eqnarray*} \lefteqn{ 
   ((\overline{\pi}\times\pi )(a)r(1),J^{*-1}e_j\otimes \pi (b^* )e_l )
=\sum_i ((\overline{\pi}\otimes\pi )\D(a)(Je_i\otimes e_i ),J^{*-1}e_j\otimes \pi (b^* )e_l ) } \\
 &=& \sum_i ((\overline{\pi}\otimes\pi )(1\otimes b)(\overline{\pi}\otimes\pi )\D(a)(Je_i\otimes e_i ),J^{*-1}e_j\otimes e_l ) \\ 
&=& \sum_i ((\overline{\pi}\otimes\pi )((1\otimes b)\D(a))(Je_i\otimes e_i ),J^{*-1}e_j\otimes e_l ) \\
 &=& \sum_{ik} (\overline{\pi}(a_k )Je_i \otimes\pi (b_k )e_i, J^{*-1}e_j\otimes e_l )
=\sum_{ik} (J\pi S(a_k^* )e_i \otimes\pi (b_k )e_i, J^{*-1}e_j\otimes e_l ) \\
 &=& \sum_{ik} (J\pi S(a_k^* )e_i, J^{*-1}e_j )(\pi (b_k )e_i, e_l )
=\sum_{ik} (e_j, \pi S(a_k^* )e_i )(\pi (b_k )e_i, e_l )\\
 &=& \sum_{ik} (\pi S(a_k^* )^* e_j, e_i)(\pi (b_k )e_i, e_l )
=\sum_k (\pi (b_k )\sum_i (\pi S(a_k^* )^* e_j, e_i )e_i, e_l ) \\
 &=& \sum_k (\pi (b_k )\pi S(a_k^* )^* e_j,  e_l )
=(\pi (\sum_k b_k S(a_k^* )^* )e_j,  e_l ) 
=(\pi (\varepsilon (a)b )e_j, e_l ),
\end{eqnarray*}
as desired.
\qed

Note that for a discrete aqg $(A,\D )$, the intrinsic dimension $d(\pi)$ of an irreducible object 
$\pi\in \Rep_f (A,\D )$ is then given by
\[ d(\pi) =r^* \mcirc 
  r=\overline{r}^*\mcirc \overline{r}=\mathrm{Tr}\,\pi(f)=\mathrm{Tr}\,\pi(f^{-1}).\] 
In fact, the latter two expressions can be thought of as the quantum dimension of $\pi$ \cite{RT}
and gives the intrinsic dimension for any $\pi\in\Rep_f (A,\D )$. By the Schwarz inequality for the
inner product given by $\mathrm{Tr}$ it follows that if $d(\pi )$ equals the dimension of the
Hilbert space on which $\pi$ acts, then $\pi (f)=1$. Thus the von Neumann extension of $(A,\D )$
is a Kac algebra iff the intrinsic dimension of any finite dimensional representation coincides with
the dimension of the Hilbert space on which it acts.
\medskip

Suppose now that $(A,\D)$ is a compact agq with counit $\varepsilon$ and coinverse $S$. We would
like to find an expression for the conjugate unitary corepresentation $\overline{V}$ of a non-zero
object $V\in \Corep_f (A,\D )$. To this end, we will use the correspondence between 
$\Rep_f(\hat{A},\hat{\D} )$ and $\Corep_f (A,\D)$ established in Theorem \ref{t-dualcor}. 
Recall that the dual $(\hat{A},\hat{\D})$ of $(A,\D )$ is a discrete aqg. Let $\hat{\varepsilon}$
denote its counit, $\hat{S}$ its coinverse, and let $\hat{f}\in M(\hat{A})$ be as in Theorem
\ref{t-disc2}. 

Let $V\in A\otimes B(K)$ be a non-zero unitary corepresentation of $(A, \D)$ on a finite dimensional
Hilbert space $K$. Pick a finite dimensional Hilbert space $\overline{K}$ and a antilinear map
$J:K\rarr\overline{K}$ such that $J^* J=\pi_V (\hat{f}^{-1})\in B(K)$. Given any invertible 
antilinear map $J$ we define a linear map $j :B(K)\rightarrow B(\overline{K})$ by 
$j(x)=J x^*J^{-1}$, for all $x\in B(K)$. 

\bprop Suppose $(A,\D)$ is a compact aqg.
Let $V$ be a unitary corepresentation of $(A,\D)$ on a finite dimensional Hilbert space $K$ 
with $J^* J=\pi_V (\hat{f}^{-1})\in B(K)$.
Define $\overline{V}\in A\otimes B(\overline{K})$ by $\overline{V}=(S^{-1}\otimes j)V$.
Then the following hold:
\begin{enumerate}
\item $V (1\otimes J^* J )=(1\otimes J^* J)(S^2\otimes \iota )V$,
\item $\overline{V}\in \Corep_f (A,\D )$, 
\item $\overline{\pi_V}=\pi_{\overline{V}}$. (Here it is understood that the same 
$J: K\rarr\overline{K}$ is used in both constructions.)
\end{enumerate}
\eprop

\prf
To prove 1., observe that
\[(\omega\otimes\iota )V J^* J=\pi_V (\omega )\pi_V (\hat{f}^{-1})=\pi_V (\omega\hat{f}^{-1} )
=\pi_V (\hat{f}^{-1} )\pi_V (\hat{f}\omega\hat{f}^{-1} )\]
\[=J^* J\pi_V (\hat{S}^2 (\omega )) =J^* J (\hat{S}^2 (\omega )\otimes\iota )V
=J^* J (\omega\otimes\iota )(S^2\otimes \iota )V \]
for all $\omega\in\hat{A}$.

To prove 3., note first that $\hat{S}(\omega)(a)=\omega S(a)$ and 
$\omega^* (a)=\overline{\omega (S(a)^*)}$ for all $a\in A$ and $\omega\in\hat{A}$.
Writing $V=\sum_k a_k\otimes x_k$ we compute
\begin{eqnarray*}
\sum_k \hat{S}(\omega^* )(a_k )x_k &=& \sum_k \omega^* S(a_k )x_k 
  =(\omega^*\otimes\iota )(S\otimes\iota )V =(\omega^*\otimes\iota )(V^* ) \\
  &=& \sum_k \omega^* (a_k^* )x_k^* = \sum_k \overline{\omega (S(a_k^* )^* )}x_k^*
  =\sum_k \overline{\omega S^{-1}(a_k)}x_k^*. 
\end{eqnarray*}
Since $(S\otimes\iota )V =V^*$, we thus obtain
\begin{eqnarray*}
  \overline{\pi_V}(\omega ) &=& J\pi_V \hat{S}(\omega ^* )J^{-1}
   = J(\hat{S}(\omega^* )\otimes\iota )VJ^{-1}
   =\sum_k J(\hat{S}(\omega^* )(a_k )x_k )J^{-1}
   =\sum_k J(\overline{\omega S^{-1}(a_k )}x_k^* )J^{-1} \\
  &=& \sum_k \omega S^{-1}(a_k )Jx_k^* J^{-1} 
  =(\omega\otimes\iota)(S^{-1}\otimes j)(\sum_k a_k \otimes x_k )
  =(\omega\otimes\iota)(S^{-1}\otimes j)V =\pi_{\overline{V}} (\omega )
\end{eqnarray*}
for all $\omega\in\hat{A}$, proving 3. Now 3., together with Proposition \ref{p-functP} and
Proposition \ref{p-conj}, imply 2., completing the proof.
\qed

This suggests the following proposition, which is a formulation more intrinsic to 
the category $\Corep_f (A,\D )$.

\bprop 
Let $(A,\D )$ be a compact aqg and $V\in \Corep_f (A,\D )$ be non-zero. Pick a finite dimensional
Hilbert space $\overline{K}$ and any antilinear map $J:K\rarr\overline{K}$ such that 
$V (1\otimes J^* J )=(1\otimes J^* J)(S^2\otimes \iota )V$. Then 
$\overline{V}\in A\otimes B(\overline{K})$ given by $\overline{V}=(S^{-1}\otimes j)V$ belongs to
$\Corep_f (A,\D)$. It is a conjugate of $V$ with normalized solutions 
$r\in\Hom (1\otimes 1,\overline{V}\times V )$,  
$\overline{r}\in\Hom (1\otimes 1, V\times\overline{V})$ of the conjugate equations for $V$ and
$\overline{V}$ given by
\[r(1)=\sum_i Je_i\otimes e_i\  \  \  \  \ and \  \  \  \  \overline{r}(1)=\sum_i e_i\otimes
   J^{*-1}e_i. \]
\eprop

\prf The irreducible case follows from Proposition \ref{p-conj}, Theorem \ref{t-dualcor} and the
fact that $J^* J$ is the unique (up to a positive scalar) positive operator in $B(K)$ 
with the property $V (1\otimes J^* J )=(1\otimes J^* J)(S^2\otimes \iota )V$ \cite{woro1,RT}.
We content ourselves with proving one of the least obvious parts of the proposition, namely that
$r$, defined above, belongs to $\Hom (1\otimes 1,\overline{V}\times V )$. In other words we must
show that 
$r(\omega\otimes\iota )(1\otimes 1)=(\omega\otimes\iota\otimes\iota )(\overline{V}\times V)r$ 
or $\omega (1)r(1)=(\omega\otimes\iota\otimes\iota )(\overline{V}\times V)r(1)$
for all $\omega\in\hat{A}$. Write $V=\sum_k a_k\otimes x_k \in A\otimes B(K)$, so
$\overline{V}=\sum_k S^{-1}(a_k )\otimes Jx_k^* J^{-1}$. For any $\omega\in\hat{A}$ and any
elements $e_j, e_r$ in the chosen orthonormal basis for $K$, we 
then get
\begin{eqnarray*} \lefteqn{
  ((\omega\otimes\iota\otimes\iota )(\overline{V}\times V)r(1),J^{*-1} e_j\otimes e_r ) 
  =\sum_i ((\omega\otimes\iota\otimes\iota )(\overline{V}_{12}V_{13})(Je_i\otimes e_i ),J^{*-1}
   e_j\otimes e_r ) } \\ 
&=& \sum_{ikl} \omega(S^{-1}(a_k )a_l )(Jx_k^* J^{-1} Je_i, J^{*-1}e_j )(x_l e_i, e_r )
\quad=\quad\sum_{ikl} \omega(S^{-1}(a_k )a_l )(Jx_k^* e_i, J^{*-1}e_j )(x_l e_i, e_r ) \\
&=& \sum_{ikl} \omega(S^{-1}(a_k )a_l )( e_j, x_k^* e_i )(x_l e_i, e_r )
\quad=\quad\sum_{ikl} \omega(S^{-1}(a_k )a_l )(x_k  e_j, e_i )(x_l e_i, e_r ) \\
&=& \sum_{kl} \omega(S^{-1}(a_k )a_l )(x_l \sum_i (x_k  e_j, e_i )e_i, e_r )
\quad=\quad\sum_{kl} \omega(S^{-1}(a_k )a_l )(x_l x_k  e_j, e_r ) \\
&=& \sum_{kl} \omega S^{-1}(S(a_l )a_k )(x_l x_k  e_j, e_r )
\quad=\quad\sum_{kl} \omega S^{-1}(a_l^* a_k )(x_l^* x_k  e_j, e_r ) \\
&=& (\omega S^{-1}\otimes \omega_{ e_j, e_r})(V^* V)
\quad=\quad(\omega\otimes \omega_{ e_j, e_r})(1\otimes 1 )
\quad=\quad\omega (1)( e_j, e_r ) \\
&=& \sum_i \omega (1)(Je_i, J^{*-1}e_j )(e_i, e_r ) 
\quad=\quad(\omega (1)r(1), J^{*-1}e_j \otimes e_r ),
\end{eqnarray*}
where 
\[ \omega_{ e_j, e_r}:\ B(K)\rarr\7C, \quad x\mapsto (xe_j, e_r ). \]
Thus
\[ (\omega\otimes\iota\otimes\iota )(\overline{V}\times V)r(1)=\omega (1)r(1), \]
completing the proof.
\qed

\brem 
In the last two results, the conjugates can be expressed alternatively in terms of the unitary
antipode $R$ and an antiunitary $J$ \cite{BMT}.
\erem


\section{Cocommutative Algebraic Quantum Groups vs. Compact Groups}

\begin{defin} \label{d-intrins}
Let $(A,\D )$ be an aqg. The intrinsic group $G$ of $(A,\D )$ is the following
subgroup of the unitaries in the multiplier algebra $M(A)$:
\[ G=\{g\in M(A)\ |\ \D g=g\otimes g,\ g^* g =gg^* =1\}. \]   
\end{defin} 

Since the extended comultiplication
$\D :M(A)\rightarrow M(A\otimes A)$ is a unital $*$-homomorphism and the algebra $M(A)$ is associative, 
the set $G$ is indeed a group. It is easy to see that $\varepsilon(g)=1$ and $S(g)=g^{-1}=g^*$, for
any $g\in G$.

\begin{rema} It can be shown that any bounded group-like element of $M(A)$
is automatically unitary with $S(g)=g^*$. A proof of this may be found in \cite{KD} and \cite{ES}, 
but the above definition suffices for our purposes.
\end{rema}

\begin{lem} \label{l-loc_cp} 
Let $(A,\D)$ be a discrete aqg. Equipped with the product topology, $M(A)=\prod_{i\in I} B(H_i )$ is
a complete locally convex topological vector space and $A$ is dense in $M(A)$.
\end{lem}
\prf Since the blocks $B(H_i)$ are Banach spaces, the functions $\| p_i(\cdot)\|$ form a separating 
family of seminorms which induces the product topology. Completeness is obvious by semisimplicity of
the $B(H_i),\ i\in I$. To see that $A$ is dense in $M(A)$, let $x\in M(A)$ and consider the net
$(x_{\lambda})_{\lambda\in\Lambda}$ in $A$ given by $x_{\lambda}=\oplus_{i\in\lambda}p_i (x)$, where
$\Lambda$ is the collection of finite subsets of $I$ directed by inclusions. Then clearly
$x_{\lambda}\rightarrow x$. 
\qed

\begin{prop}
If $(A,\D )$ is a discrete aqg, its intrinsic group $G$ is a compact topological
group w.r.t. the product topology on $M(A)$.
\end{prop}

\prf
Note that a net $(g_{\lambda})$ converges to $g$ in $G$ iff $p_i (g_{\lambda})\rightarrow p_i (g)$ in
norm, for all $i\in I$. By Tychonov's theorem it suffices to show that $G$ is closed in 
$\prod_{i\in I} U(B(H_i ))$. Given a net $(g_{\lambda})$ in $G$ which converges to $a\in M(A)$, we must
show that $\D a=a\otimes a$. Let $i,j\in I$ and consider the finite dimensional non-degenerate
$*$-representation $(p_i\otimes p_j )\D :A\rightarrow B(H_i\otimes H_j )$. By Proposition \ref{Berkeley},
we may write
\[ (p_i\otimes p_j )\D \simeq\bigoplus_{k\in I}\ N_{ij}^k\, p_k.  \]
Since both these expressions are non-degenerate $*$-representations, they extend to equivalent
unital $*$-representations of $M(A)$, which we now identify.
Therefore
\[ (p_i\otimes p_j )\D a =\bigoplus_{k\in I}\ N_{ij}^k\, p_k (a)
=\lim_{\lambda}\bigoplus_{k\in I}\ N_{ij}^k\, p_k (g_{\lambda})
=\lim_{\lambda}(p_i\otimes p_j )\D g_{\lambda} \]
\[ =\lim_{\lambda}p_i (g_{\lambda})\otimes p_j (g_{\lambda})
=p_i (a)\otimes p_j (a)=(p_i\otimes p_j )(a\otimes a), \]
thus $\D a =a\otimes a$ and $G$ is closed. \qed

\begin{rema}
Note furthermore that in the proof of the proposition we have shown that 
\[ \D:M(A)=\prod_{i\in I}B(H_i )\rightarrow M(A\otimes A)=\prod_{i,j\in I}B(H_i\otimes H_j ) \]
is continuous w.r.t. the product topologies.
\end{rema}

\begin{prop} \label{p-functor_D}
Suppose $(A,\D )$ is a discrete aqg with intrinsic group $G$. Let $\pi$ be a non-degenerate
$*$-representation of $A$ on a finite dimensional Hilbert space $K$. Define a map 
$u_{\pi} :G\rightarrow B(K)$ by 
\[ u_{\pi} : g\mapsto \pi (g), \]
where the extension of $\pi$ to $M(A)$ is understood. Then $u_{\pi}$ is strongly continuous and
$D: (K,\pi) \mapsto (K,u_{\pi})$ together with the identity map on the
morphisms is a faithful and tensor functor from
$\Rep_f(A,\D )$ to the category $\Rep_f\,G$ of finite dimensional continuous representations of $G$.
\end{prop}

\prf Continuity of $u_\pi$ w.r.t.\ the topology on $G$ is obvious.
As $\pi :M(A)\rightarrow B(K)$ is a unital $*$-homomorphism, clearly $(K, u_{\pi})\in\Rep\,G$.
Functoriality and faithfulness of $D$ are obvious. Monoidality follows from the calculation
\[ D(\pi_1\times\pi_2)(g)=(\pi_1\otimes\pi_2 )\D g=\pi_1 (g)\otimes\pi_2 (g)
=u_{\pi_1}(g)\otimes u_{\pi_2 }(g)=(u_{\pi_1}\times u_{\pi_2 })(g) \]
for all $g\in G$ and $\pi_i\in \Rep(A,\D )$. Every $\pi\in \Rep(A,\D )$ is equivalent to a direct
sum of the representations $p_i$ by Proposition \ref{Berkeley}.
Each $u_{p_i}$ is (strongly) continuous and therefore the direct sum is strongly continuous. \qed

Thus far we have not assumed $(A,\D)$ to be cocommutative. The following characterization of 
cocommutativity will be crucial for proving that $D$ gives rise to an equivalence of categories.

\begin{thm} \label{t-cocomm1}
A discrete aqg $(A,\D)$ is cocommutative iff 
\[ \overline{\mathrm{span}_\7C \{ g\ |\ g\in G \}}=M(A). \]
\end{thm}

\begin{rema} 
This theorem implies that 
\[ \mathrm{span}_\7C \{ p_i(g)\ |\ g\in G \}=B(H_i) \]
for every $i\in I$ and in particular that
$A=\mathrm{span}_\7C \{ gI_i \ |\ i\in I, g\in G \}$. These results are rigorous formulations of
the heuristic idea that a cocommutative aqg is `spanned by its grouplike elements'. Before we give
the proof of Theorem \ref{t-cocomm1} we show that it leads to the desired equivalence of categories. 
\end{rema}

In what follows we fix a cocommutative discrete aqg $(A,\D)$ where $A=\oplus_{i\in I} B(H_i)$ and
let $G$ denote its intrinsic group. 
For every $i\in I$, $u_i : g\mapsto p_i(g)\in B(H_i)$ is a continuous unitary representation of
$G$. By Theorem \ref{t-cocomm1}, the span of $p_i(g),\ g\in G$ is dense in $B(H_i)$, thus $u_i$ is
irreducible. This defines a map $\gamma: I\rightarrow I_G$. 

\begin{prop} \label{p-biject}
The map $\gamma$ is a bijection.
\end{prop} 
\prf 
If there is a unitary $V: H_i\rightarrow H_j$ such that $Vp_i(g)=p_j(g)V$ for all $g\in G$ then 
$Vp_i(x)=p_j(x)V$ for all $x\in \mathrm{span}_\7C \{ g\ |\ g\in G \}$. Since these $x$ are dense in
$M(A)$ by Theorem \ref{t-cocomm1} and the $p_i$ are continuous we conclude that $Vp_i(x)=p_j(x)V$
for all $x\in M(A)$, thus $i=j$ so that $\gamma: I\rightarrow I_G$ is injective. 

Obviously if $g\ne e$ there is an $i\in I$ such that $p_i(g)\ne 1_{B(H_i)}$. Since the category
$\Rep_f(A,\D)$ is monoidal and has conjugates, $\gamma(I)\subset I_G$ is closed w.r.t.\ conjugation
and tensor products and reduction. The surjectivity of $\gamma$ now follows from the following well
known group theoretical fact.
\qed

\begin{lem} \label{l-surj}
Let $G$ be a compact group and let $J\subset I_G$ be closed w.r.t.\ conjugation and tensor products
and reduction. If $J$ separates points on $G$ then $J=I_G$.
\end{lem}
\prf For every equivalence class $i\in I_G$ pick a representative $u_i$. The assumptions on $J$
imply that the span of the matrix elements $(u_j)_{nm},\ j\in J$ is a unital $*$-subalgebra of
$C(G)$. Since it separates the points of $G$ it is dense in $C(G)$ by the Stone-Weierstrass theorem
and therefore dense in $L^2(G,\mu)$, where $\mu$ is the Haar measure. If there were a
$k\in I_G\backslash J$ then by the Peter-Weyl theorem the matrix elements of $u_k$ would be
orthogonal to the dense subspace of $L^2(G,\mu)$ generated by the $(u_j)_{nm},\ j\in J$, which is a 
contradiction. 
\qed

\begin{thm} \label{t-cocomm2}
Let $(A,\D)$ be a cocommutative discrete aqg and $G$ its intrinsic group. Then the functor
$D :\Rep(A,\D )\rightarrow\Rep\,G$ induces a canonical equivalence of tensor $*$-categories: 
\begin{eqnarray*} 
  \Rep_f(A,\D) & \stackrel{\otimes}{\simeq} & \Rep_f\,G.
\end{eqnarray*}
\end{thm}
\prf In view of Proposition \ref{p-functor_D} it only remains to prove that the functor $D$ is full
and essentially surjective. By Proposition \ref{Berkeley} the category $\Rep_f(A,\D)$ is semisimple,
and for a compact group $G$ the semisimplicity of $\Rep_f G$ is well known. Recall from Section
\ref{ss-categ} that a faithful functor between semisimple categories is full if and only if it maps
simple objects to simple objects and non-isomorphic simple objects have non-isomorphic images. The
first property was used to define the map $\gamma$ of Proposition \ref{p-biject}, and the 
second is the injectivity of $\gamma$. Finally, essential surjectivity of $D$ is expressed by
surjectivity of $\gamma$. Since $D$ is monoidal it gives rise to an equivalence of tensor categories
by \cite{SR}.
\qed

Now we prove the characterization of cocommutative discrete \aqgs\ used above.\\

\noindent \emph{Proof of Theorem \ref{t-cocomm1}.}
Clearly, if the closure of $\mathrm{span}_\7C\{g\, |\, g\in G\}$ is $M(A)$, then by linearity and
continuity of $\D :M(A)\rightarrow M(A\otimes A)$, we see that $(A,\D )$ is cocommutative.
We prove the converse direction. By Lemma \ref{l-loc_cp} it suffices to show that any $x\in A$ is
the limit of linear combinations of elements of $G$. Fix $x\in A$.
For any aqg $(A,\D)$ there exists a linear inclusion $Q:A\hookrightarrow \hat{A}_r^*$ determined by
\[ Q(b)\hat{\pi}_r (\hat{a})=\hat{a}(b)=\varphi (ba) \]
for all $a,b\in A$. To see this, to any $a,b\in A$ we can choose an element $c\in A$ such that 
$\hat{c}\hat{b^*}^* =\hat{b^*}^*$. (Such $c$ can be obtained as inverse Fourier transform of a
local unit for $\hat{b^*}^*$ in $\hat{A}$.) Now we observe that
\[ Q(b)\hat{\pi}_r (\hat{a})=\varphi((b^*)^*a)
   =\hat{\psi}(\hat{b^*}^*\hat{a})=\hat{\psi}(\hat{c}\hat{b^*}^*\hat{a}) 
   =\hat{\psi}(\hat{b^*}^*\hat{a}\hat{\rho}(\hat{c}))
   =(\hat{\pi}_r (\hat{a})\hat{\Lambda}(\hat{\rho}(\hat{c})),\hat{\Lambda}(\hat{b^*} )), \]
so $Q(b)=(\cdot\hat{\Lambda}(\hat{\rho}(\hat{c})),\hat{\Lambda}(\hat{b^*} ))$.
Now $\hat{A}_r$ is a unital commutative $C^*$-algebra. Let $Y$ be the set of $*$-characters on
$\hat{A}_r$. Gelfand's theorem tells us that $Y$ is a compact Hausdorff space and that
the map ${\Gamma}:\hat{A}_r\rightarrow C(Y)$ given by ${\Gamma}(a)(y)=y(a)$, for all $a\in\hat{A}_r$
and $y\in Y$, is a unital $*$-isomorphism from $\hat{A}_r$ to the $C^*$-algebra $C(Y)$ of continuous
functions on $Y$. Thus $Q(x){\Gamma}^{-1}\in C(Y)^*$. By the Krein-Milman theorem for probability
measures \cite[Theorem 2.5.4]{ped} any element of $C(Y)^*$ is a $w^*$-limit of linear combinations
of Dirac measures $\delta_y\in C(Y)^*$. Hence we may write 
\[ Q(x){\Gamma}^{-1}=\lim_\lambda\sum_k c_{\lambda k}\delta_{y_{\lambda k}} \]
for some $c_{\lambda k}\in \7C$ and $y_{\lambda k}\in Y$. Thus we get
\begin{eqnarray*}
\hat{a}(x) &=& Q(x)\hat{\pi}_r (\hat{a})=Q(x){\Gamma}^{-1}{\Gamma}\hat{\pi}_r (\hat{a})
=\lim_\lambda\sum_k c_{\lambda k}\delta_{y_{\lambda k}}{\Gamma}\hat{\pi}_r (\hat{a}) \\
 &=& \lim_\lambda\sum_k c_{\lambda k}{\Gamma}\hat{\pi}_r (\hat{a})(y_{\lambda k})
=\lim_\lambda\sum_k c_{\lambda k}\, y_{\lambda k}\hat{\pi}_r(\hat{a})
\end{eqnarray*}
for all $a\in A$. Now the $y_{\lambda k}\hat{\pi}_r :\hat{A}\rightarrow\7C$ are unital
$*$-homomorphisms, so we may define $g_{\lambda k}\in M(A)$ by 
\[ g_{\lambda k}=(\iota\otimes y_{\lambda k}\, \hat{\pi}_r )U, \]
where $U\in M(A\otimes\hat{A})$ is the universal corepresentation of $(A,\D)$.
Since $U$ is unitary, so are the elements $g_{\lambda k}$ 
and moreover,
\begin{eqnarray*}
\D g_{\lambda k} &=&\D (\iota\otimes y_{\lambda k}\, \hat{\pi}_r )U
=(\iota\otimes\iota\otimes y_{\lambda k}\, \hat{\pi}_r)(\D\otimes\iota )U
=(\iota\otimes y_{\lambda k}\, \hat{\pi}_r)(U_{13}U_{23}) \\
 &=& (\iota\otimes y_{\lambda k}\, \hat{\pi}_r)U\otimes(\iota\otimes y_{\lambda k}\, \hat{\pi}_r)U
=g_{\lambda k}\otimes g_{\lambda k},
\end{eqnarray*}
thus all $g_{\lambda k}\in G$.

Now for any $\xi\in B(H_i)'$ and $i\in I$ observe that $\xi p_i\in \hat{A}$ since $\varphi$
is, up to a factor $d_i$, the trace on $B(H_i)$. Thus by Theorem \ref{t-univ-corep} we have 
$(\xi p_i\otimes\iota )U =\xi p_i$. Hence by the previous formulae we get
\begin{eqnarray*}
\xi p_i (x) &=& \lim_\lambda\sum_k c_{\lambda k}\, y_{\lambda k}\hat{\pi}_r(\xi p_i )
=\lim_\lambda\sum_k c_{\lambda k}\, y_{\lambda k}\hat{\pi}_r(\xi p_i\otimes\iota )U \\
 &=& \lim_\lambda\sum_k c_{\lambda k}\, \xi p_i (\iota\otimes y_{\lambda k}\, \hat{\pi}_r )U
=\lim_\lambda\sum_k c_{\lambda k}\, \xi p_i (g_{\lambda k})
=\lim_\lambda\xi p_i(\sum_k c_{\lambda k}\,  g_{\lambda k}).
\end{eqnarray*}
Since this is true for all $i\in I$ and $\xi\in B(H_i)'$, we thus get
\[ x=\lim_\lambda\sum_k c_{\lambda k}\,  g_{\lambda k} \]
as desired. 
\qed

This concludes the proof of Theorem \ref{t-cocomm1}. 

\begin{cor} \label{c-span}
For every $i\in I$ we have $\ \mathrm{span}_\7C \{ p_i(g)\ |\ g\in G \}=B(H_i)$.
\end{cor}
\prf First we show that for $0\ne a\in A$, there exist $g\in G$ such that $\varphi (ga)\neq 0$.
Thus let $0\ne a\in A$. Then there exists
$x\in A$ such that $\hat{a}(x)=\varphi(xa)\ne 0$, thus
$\hat{a}=\varphi(\cdot a)\ne 0$. Furthermore, there exists a finite set $\Lambda\subset I$ such that  
\[ \hat{a}(b)=\varphi (ba)=\sum_{i\in\Lambda} d_i Tr_i (p_i (ba)), \]
for all $b\in M(A)$. Since the product in $M(A)$ is continuous, we thus see that 
$\hat{a}: M(A)\rightarrow\7C$ is continuous. Thus by Theorem \ref{t-cocomm1} 
there exists $g\in G$ such that $\varphi(ga)\ne 0$. 

Assume now that $\mathrm{span}_\7C \{p_i(g) \ | \ g\in G\}\ne B(H_i)$. Then there exists a non-zero
$a\in B(H_i)$
such that $Tr_{B(H_i)}(p_i(g)a)=0$ for all $g\in G$. Since $\varphi\restr B(H_i)=d_i Tr_{B(H_i)}$,
this is a contradiction. 
\qed

Let $G$ be a compact group and consider a unitary representation $u:G\rightarrow B(K)$ of $G$ on some
finite dimensional Hilbert space $K$. Since $u$ is continuous (here with respect to the norm), we
have $u\in C(G,B(K))$. Since $C(G,B(K))\cong C(G)\otimes B(K)$, we may thus regard $u$ as 
a unitary element $\tilde{V}$ of the unital $*$-algebra $C(G)\otimes B(K)$.
Now $u:G\rightarrow B(K)$ being a representation simply means that 
$(\D\otimes\iota )\tilde{V} =\tilde{V}_{13}\tilde{V}_{23}$,
where $\D :C(G)\rightarrow C(G\times G)\cong C(G)\otimes C(G)$ is given by $\D (a)(g,g')=a(gg' )$,
for all $a\in C(G)$ and $g,g'\in G$. Thus $\tilde{V}$ is a unitary corepresentation of the 
commutative compact quantum group $(C(G),\D)$ in the sense of Woronowicz.
It is easy to see that one may recover $u$ from $\tilde{V}$ by
$u(g)=(\delta_g\otimes\iota )\tilde{V}$, for $g\in G$.

\btheor
Let $(A,\D)$ be a cocommutative discrete aqg, so $(\hat{A},\hat{\D})$ is a commutative compact aqg.
Let $G$ be the compact group such that $(C(G),\D)\cong (\hat{A}_r, \hat{\D}_r )$.
For any $V\in\Corep_f (\hat{A},\hat{\D})$, let
$\tilde{V}\in C(G)\otimes B(K)$ be the unitary corepresentation of $(C(G),\D)$ 
given by $\tilde{V} =(\hat{\pi}_r\otimes\iota )V$.
Then there exists an isomorphism from $G$ to the intrinsic group of $(A,\D )$ given by
\[g\mapsto\tilde{g}=(\iota\otimes\delta_g\hat{\pi}_r )U \in M(A),\]
where $U$ is the universal corepresentation of $(A,\D )$.
Moreover, the equivalence
\[ \Corep_f (\hat{A},\hat{\D})\stackrel{P}{\lrarr}\Rep_f (A,\D)\stackrel{D_f}{\lrarr}\Rep_f G\]
given by the composition
\[ V\mapsto \pi_V\mapsto u_{\pi_V}=:u_V, \]
with $P$ and $D_f$ established in Theorem \ref{t-dualcor} and Theorem \ref{t-cocomm2}, is described
more directly by
\[ u_V (\tilde{g})=(\delta_g\otimes\iota )\tilde{V},\]
for all $g\in G$.  
\etheor

\prf
From the proof of Theorem \ref{t-cocomm1}, it is easily seen that $g\mapsto\tilde{g}$ is indeed a 
homeomorphism from $G$ to the intrinsic group of $(A,\D)$ which also preserves the products.
Let $g\in G$. We proceed to show that $u_V (\tilde{g})=(\delta_g\otimes\iota )\tilde{V}$.
Since $V\in\hat{A}\otimes B(K)$ we may write $V=\sum_k \omega_k\otimes x_k$ for $\omega_k\in\hat{A}$
and $x_k\in B(K)$. Then
\begin{eqnarray*} \lefteqn{ 
u_V (\tilde{g}) = u_V((\iota\otimes\delta_g\hat{\pi}_r )U)
=\pi_V ((\iota\otimes\delta_g\hat{\pi}_r )U) } \\
&=& ([(\iota\otimes\delta_g\hat{\pi}_r )U]\otimes\iota )V
=\sum_k [(\iota\otimes\delta_g\hat{\pi}_r )U](\omega_k )x_k
=\sum_k (\omega_k\otimes\delta_g\hat{\pi}_r )U x_k \\
&=& \sum_k \delta_g\hat{\pi}_r ((\omega_k\otimes\iota )U) x_k
=\sum_k \delta_g\hat{\pi}_r (\omega_k ) x_k 
=(\delta_g\hat{\pi}_r\otimes\iota )(\sum_k \omega_k \otimes x_k ) 
=(\delta_g\hat{\pi}_r\otimes\iota )V =(\delta_g\otimes\iota )\tilde{V},
\end{eqnarray*}
as claimed. 
\qed


\section{From Concrete Tensor Categories to Discrete Quantum Groups}
The strategy we adopt in this section was been pursued by J.E. Roberts \cite{ro} in the discrete Kac
algebra case, and has then been adapted to the more algebraic setting of discrete \aqgs\ $(A, \D)$
by the authors.  

\subsection{Natural Transformations of $*$-Preserving Functors} \label{ss-nattrf}

\bprop
\label{p-defprod}
Let $\2D, \2K$ be $*$-categories and $F: \2D\rarr\2K$ a $*$-preserving functor.
Let $\Nat\,F$ be the set of all natural transformations from $F$ to itself, viz.\ 
\[ \Nat\,F= \{ b=(b_X \in \End\,F(X), \ X\in\2D ) \ \ | \ \ 
   F(s)\mcirc  b_X=b_Y \mcirc  F(s) \ \ \forall s: X\rarr Y\}.
\]   
Then with the following pointwise operations $\Nat\,F$ is a unital $*$-algebra:
\begin{eqnarray*}
   (\lambda b +\lambda' b')_X &=& \lambda b_X + \lambda' b'_X,  \\
  (b b')_X &=& b_X \mcirc  b'_X,  \\
  (b^*)_X &=& (b_X)^*,  \\
  \11_X &=& id_{F(X)},
\end{eqnarray*}
for $b, b' \in\Nat\,F,\ \lambda,\lambda'\in\7C$ and $X\in\2D$.
\eprop
\prf If $s\in \Mor(X, Y )$, then 
\[ F(s)\mcirc  (b^* )_X = (b_X \mcirc  F(s^* ))^* =(F(s^* )\mcirc  b_Y )^* =(b^* )_Y \mcirc  F(s), \]
where we have used $F(s^* )=F(s)^*$. Thus $b^*\in\Nat\,F$. Similarly, we get  
$b b', \ \lambda b + \lambda' b' \in\Nat\,F$, so with these operations $\Nat\,F$ is a
$*$-algebra. Obviously $\11$ is a unit for $\Nat\,F$.
\qed

\bdefin 
Let $B$ be a $*$-algebra and $\2K$ a $*$-category. Then the category $\Rep_\2K\,B$ of
`representations of $B$ in $\2K$' has as objects the pairs $(X, \pi_X)$, where $X\in\obj\,\2K$ and
$\pi_X$ is a non-degenerate $*$-homomorphism of $B$ into $\End\,X$. The morphisms are given by
\[ \Mor_{\Rep_\2K B} ((X,\pi_X), (Y,\pi_Y))=\{ s\in\Mor_\2K(X,Y) \ | \ s\mcirc \pi_X(b)
   =\pi_Y(b)\mcirc  s \ \ \forall b\in B\}. \]
\edefin

Note that if $\2K=\2H$, the category of finite dimensional Hilbert spaces, we have
$\End\,H=B(H)$. Thus our definition reduces to the usual notion of finite dimensional non-degenerate
$*$-representations of $B$: $\Rep_\2HB=\Rep_f\,B$.

\bprop \label{p-funct-F}
Let $\2C, \2K$ be $*$-categories and $E: \2C\rarr\2K$ a $*$-preserving functor. Define a unital
$*$-algebra $B=\Nat\,E$. There is a $*$-preserving functor $F: \2C\rarr\Rep_\2K\,B$ satisfying
$K\mcirc  F=E$, where $K: \Rep_\2K B\rarr\2K,\ (X,\pi_X)\mapsto X$ is the forgetful functor.
\eprop
\prf
We set $F(X)=(E(X), \pi_X)$ for $X\in\obj\,\2C$ and $F(s)=E(s)$ for the morphisms of $\2C$. Here
the $*$-representation $\pi_X$ is defined by
\[ \pi_X: \ B\rarr \End\,X, \quad b\mapsto b_X. \]
That $F$ really is a functor follows from
\[ F(s)\mcirc \pi_X(b)=E(s)\mcirc  b_X=b_Y\mcirc  E(s)= \pi_Y(b)\circ F(s)
   \quad\quad \forall s:X\rarr Y\ \forall b\in B, \]
where we used $b\in B=\Nat\,E$.
\qed

\begin{prop} \label{p-direct}
Let $\2D$ be a semisimple $*$-category, $\2K$ a $*$-category and $F: \2D\rarr\2K$ a $*$-preserving
functor. Then there is an isomorphism
\[ \psi_F: \ \Nat\,F \rarr \prod_{i\in I_\2D} \End\,F(X_i) \]
of $*$-algebras.
\end{prop}  
\prf
We define a map 
$\psi_F :\Nat\,F\rightarrow\prod_{i\in I_\2D} \End\,F(X_i)$ by 
\[ \psi_F (b)=\prod_{i\in I_\2D} b_{X_i} \]
for all $b\in\Nat\,F$.
Clearly $\psi_F$ is a unital $*$-algebra homomorphism. We first show that it is injective.
Suppose $\psi_F (b)=0$, so $b_{X_i}=0$ for all $i\in I_\2D$. Let $X\in\2D$. We must show that
$b_X=0$. Recall that, for any $i\in I_\2D,\ \Mor(X_i, X )$ is a Hilbert space with the inner product 
$(s,t)\id_{X_i}=t^*s$ and that every $s\in\Mor(X_i, X)$ is a multiple of an isometry.
Let $s_{i\alpha}\in \Mor(X_i, X ), \ i \in I_\2C,\ \alpha=1,\ldots, \dim \Mor(X_i, X )$
be an orthonormal basis w.r.t.\ this inner product satisfying 
$\sum_{i\alpha} s_{i\alpha} \mcirc  s_{i\alpha}^* =\id_{X}$. Hence 
\[ b_X= b_X \mcirc  F(\id_{X})=\sum_{i\alpha} b_X\mcirc  F(s_i \mcirc  s_i^* )
=\sum_{i\alpha} F(s_i )\mcirc  b_{X_i}\mcirc  F(s_i^* )=0, \]
and injectivity of $\psi_F$ follows.

Next we prove surjectivity of $\psi_F$. Given $( b_i\in\End\,F(X_i),\ i\in I_\2D )$ we need to 
construct $b\in\Nat\,F$ such that $\psi_F (b)=\prod_{i\in I_\2D} b_i$.
Thus let $X\in\2D$. As above, pick a basis of isometries $s_{i\alpha}\in \Mor(X_i, X )$, 
$i\in I_\2D$, with $\sum_{i\alpha} s_{i\alpha} \mcirc  s_{i\alpha}^* =\id_{X}$. 
Define $b_X\in\End\,F(X)$  by
\[ b_X =\sum_{i\alpha} F(s_{i\alpha})\mcirc  b_i\mcirc  F(s_{i\alpha}^* ). \]
Another choice of orthonormal bases $(s_{i\alpha})$ in $\Mor(X_i, X),\ i\in I_\2D$ does not affect
$b_X$ since any two such bases are related by an orthogonal transformation.
It remains to show that $b\in\Nat\,F$, viz.\ 
$F(s)\mcirc  b_X=b_Y\mcirc  F(s) \ \ \forall s: X\rarr Y$.
To this purpose consider $u\in \Mor(X, Y )$.
Pick isometries $t_{j\beta}\in\Mor(Y_j,Y )$,
with $\sum_{j\beta} t_{j\beta} \mcirc  t_{j\beta}^* =\id_{Y}$, and
define $b_Y\in\End\,F(Y)$ by
\[ b_Y=\sum_{j\beta} F(t_{j\beta})\mcirc  b_j\mcirc  F(t_{j\beta}^* ). \]
Then
\[ F(u)\mcirc  b_X=\sum_{i\alpha} F(u \mcirc  s_{i\alpha} )\mcirc  b_i\mcirc  F(s_{i\alpha}^* )
  =\sum_{i\alpha,j\beta} F(t_{j\beta}\mcirc  t_{j\beta}^* \mcirc  u\mcirc  s_{i\alpha} )\mcirc  b_i \mcirc  F(s_{i\alpha}^* ). \]
But $t_{j\beta}^*\mcirc  u\mcirc  s_{i\alpha}\in\Mor(X_i,X_j)$ is zero unless $i=j$ and
otherwise a scalar multiple of $\id_{X_i}$. Hence 
\[ F(u)\mcirc  b_X
=\sum_{i\alpha\beta} F(t_{i\beta}\mcirc  (t_{i\beta}^*\mcirc  u\mcirc  s_{i\alpha})) \mcirc  b_i\mcirc  F(s_{i\alpha}^*)
=\sum_{i\alpha,j\beta} F(t_{j\beta})\mcirc  b_j\mcirc  F(t_{j\beta}^*\mcirc  u\mcirc  s_{i\alpha}\mcirc  s_{i\alpha}^*)
=b_Y\mcirc  F(u). \]
The last identity in this equation follows by precisely the same sort of argument as used before.
This proves the naturality of the family $b\in\Nat\,F$.
Obviously we have $b_{X_i}=b_i$ for all $i\in I_\2D$. Thus
$\psi_F (b)=\prod_{i\in I_\2D} b_i$, which proves surjectivity of $\psi_F$.
\qed

The isomorphism $\psi_F$ depends on the chosen family $X_i, i\in I$ of irreducibles in $\2C$, but is
canonical once the latter is fixed. Thus there is no risk in suppressing it and writing
$\Nat\,F=\prod_{i\in I_\2D} \End\,F(X_i)$, as will be done in the sequel.

\bexam
\label{e-idfunc} 
The simplest application of this proposition is to the identity $\2D\rarr\2D$, when we get
\[ \Nat\,(Id)\cong\prod_{i\in I_\2D} \7C. \]
Examples of natural transformations of the identity are the `twists' appearing in the context of 
braided tensor categories. 
\eexam

\bexam \label{e-otimes}
Suppose $\2D$ is a tensor $*$-category with tensor product $\otimes$. A less trivial application of
the proposition is the case of the tensor functor $F=\otimes :\2D\times\2D\rightarrow \2D$, where
$\2D\times\2D$ is the product $*$-category. Note that the isomorphism classes of irreducible objects
in $\2D\times\2D$ are precisely $I_\2D\times I_\2D$. By semisimplicity of $\2D$ we have a finite
decomposition 
\[ X_i\otimes X_j\cong \bigoplus_{k\in I_\2D} N_{ij}^k X_k, \]
for any $i,j\in I_\2D $ and thus 
$\End(X_i\otimes X_j)\cong \bigoplus_{k\in I_\2D} M_{N_{ij}^k}(\7C)$.
We therefore obtain
\[ \Nat\,\otimes\ \cong\ \prod_{\stackrel{i,j,k\in I_\2D}{N_{ij}^k\ne 0}} M_{N_{ij}^k}(\7C). \]
\eexam

The following example is a continuation of Example \ref{e-idfunc}.
 
\bexam Suppose $\2C$ is a semisimple $*$-category.
Consider the $*$-algebra $A=\oplus_{i\in I_\2C}\7C$.
Then Proposition \ref{p-funct-F} provides an equivalence 
$F: \2C\rarr \Rep_f\,A$ of categories such that 
$F(X_i)\cong p_i$ for all $i\in I_\2C$. Let $K:\Rep_f\,A\rarr\2H$ be the forgetful functor.
Then $K\mcirc  F$ is an embedding functor for $\2C$ as a $*$-category.
By Proposition \ref{p-direct} and Example \ref{e-idfunc} we also conclude that 
\[\Nat(K\mcirc  F)\cong\Nat\,Id\cong\prod_{i\in I_\2C}\7C.\] 
(Here of course, it does not make sense to say that the functors $Id$ and $K\mcirc  F$ are
unitarily equivalent.) Suppose now that $\2C$ is in addition a tensor category. 
Then the faithful functor $K\mcirc  F$ cannot in general be monoidal. To see this
let $X_i, X_j$ with $i,j\in I_\2C$. Then $K\mcirc  F(X_i )\cong K(p_i )\cong\7C$,
so $K\mcirc  F(X_i )\otimes K\mcirc  F(X_j )\cong\7C$, whereas 
\[ K\mcirc  F (X_i\otimes X_j )\cong \bigoplus_{k\in I_\2C}N_{ij}^k\, \, \, K\mcirc  F(X_k )
\cong\bigoplus_{k\in I_\2C}N_{ij}^k\, \, \7C, \]
where we used the decomposition $X_i\otimes X_j\cong\oplus_{k\in I_\2C}N_{ij}^k\, \ X_k$.
\eexam

Given a $*$-preserving functor $E:\2C\rarr\2K$ and writing $B=\Nat\,E$ it is natural to ask under
which conditions the functor $F: \2C\rarr\Rep_\2KB$ defined in Proposition \ref{p-funct-F} is an
equivalence of categories. Instead of pursuing this line further we restrict ourselves to the case
where $\2C$ is semisimple and $\2K=\2H$. Then $\End\,E(X)$ is a finite dimensional matrix algebra
for every $X\in\2C$ and Proposition \ref{p-direct} shows that $\Nat\,E$ is a direct product of
matrix algebras. The representation theory of such algebras is quite intricate, cf.\
\cite{GT}. There it is shown, among other results, that every irreducible representation of
$B=\prod_{i\in I} B(H_i)$ is equivalent to a projection $p_i$ iff the index set $I$ has less than
measurable cardinality. In practice, this is no serious restriction but for reasons which will
become clear in the next subsection we prefer to work with direct sums instead of direct products. 

\btheor \label{t-equiv}
Let $\2C$ be a semisimple $*$-category and $E: \2C\rarr\2H$ a $*$-preserving functor.
The $*$-algebra $A=\oplus_{i\in I_\2C} B(E(X_i))$ embeds canonically into 
$\Nat\,E\cong\prod_{i\in I_\2C} B(E(X_i))$. The functor
\[ F:\ \2C\rarr\Rep_f A, \quad F(X)=(E(X),\pi_X),\ \pi_X(a)=a_X\ \forall a\in A, \quad F(s)=E(s) \]
satisfies $K\mcirc  F=E$. It gives rise to an equivalence of categories iff $E$ is faithful.
\etheor
\prf $F$ as defined above differs from the functor defined in Proposition \ref{p-funct-F} only in
the restriction of the representations $\pi_X$ to $A\subset\Nat\,E$. Thus it clearly is a
$*$-preserving functor. By definition it is clear that $K\mcirc  F=E$. Since $F$ coincides with
$E$ on morphisms, $E$ is faithful if $F$ is an equivalence. Assume now that $E$ is faithful. Let
$X\in\2C$ be irreducible. Then $F(X)$ is a representation $\pi_X: a\mapsto a_X$ of $A$ on the
Hilbert space $E(X)$. As is clear from the proof of Proposition \ref{p-direct}, for $X=X_i$ we have
$\pi_X=p_i$, which clearly is irreducible. Since $a\in A$ is a natural transformation, 
$\pi_X: a\mapsto a_X$ is an irreducible $*$-representation whenever $X$ is irreducible. 
Consider $i\ne j, i,j\in I_\2C$. By Proposition \ref{Berkeley}, the projections 
$p_i=F(X_i), p_j=F(X_j)$ are inequivalent $*$-representations. By semisimplicity of $\Rep_f A$ this
implies that $F$ is full. Since the $p_i, i\in I_\2C$ exhaust the equivalence classes of irreducible
$*$-representations of $A$ it is clear that $F$ is essentially surjective, completing the proof.
\qed

\brem 
The theorem remains true if one replaces $A$ by the direct product algebra $\Nat\,E$ provided one 
makes the additional assumption that the set $I_\2C$ of isomorphism classes of $\2C$ has less than
measurable cardinality.
\erem


\subsection{Concrete Tensor $*$-Categories and Discrete Quantum Semigroups}
From now on, $\2C$ will be a semisimple strict tensor $*$-category. (As discussed in the
Introduction, the strictness assumption does not entail a loss of generality.) We will introduce
other requirements as we proceed. By  Theorem \ref{t-equiv} we have associated a $*$-algebra $A$ to
a $*$-preserving functor $E: \2C\rarr\2H$. The aim of this subsection is to endow $A$ with the
structure of a discrete aqg under suitable conditions on $\2C$ and $E$. To bring out the special
r\^{o}le of the target category $\2H$, we work in a more general setting as far as seems
justifiable.

Let $E: \2C\rarr\2D$ be a $*$-preserving tensor functor into some tensor $*$-category $\2D$. Recall
that by Proposition \ref{p-defprod}, $\Nat\,E$ is a unital $*$-algebra for every
$*$-preserving functor $E: \2C\rarr\2D$.

\bprop 
Let $\2C, \2D$ be tensor $*$-categories, $\2C$ being strict, and let $E: \2C\rarr\2D$ be a
$*$-preserving tensor functor. Then
\[ \tilde{\Delta}: \Nat\,E\rarr\Nat\,E\otimes, \quad b=(b_X) \ \mapsto \ \tilde{\Delta}(b)\ =\ 
   (\tilde{\Delta}(b)_{Y,Z}), \quad \tilde{\Delta}(b)_{Y,Z}=b_{Y\otimes Z} \]
defines a unital $*$-homomorphism.
\eprop
\prf The morphisms in $\2C\times\2C$ are of the form $s\times t$ where 
$s: Y\rarr Y',\ t: Z\rarr Z'$. Then $s\otimes t\in\Mor(Y\otimes Z,Y'\otimes Z')$. We compute 
\[ (E\otimes) (s\times t) \, \mcirc  \tilde{\Delta}(b)_{Y,Z} = 
   E(s\otimes t) \, \mcirc  \, b_{Y\otimes Z} = b_{Y'\otimes Z'} \, \mcirc  E(s\otimes t) = 
   \tilde{\Delta}(b)_{Y',Z'} \mcirc (E\otimes) (s\times t), \]
where we used $b=(b_X)\in\Nat\,E$. We conclude that $\tilde{\D}(b)\in\Nat\,E\otimes$. That 
$b\mapsto\tilde{\D}(b)$ is a unital $*$-homomorphism follows immediately from Proposition
\ref{p-defprod}. 
\qed

We now consider the analogues of $(\Delta\otimes\iota)$, $\2C$ being strict, and $(\iota\otimes\D)$
in the present general setting. $\otimes^2$ will denote the functor 
$\otimes\mcirc (\otimes\times Id)=\otimes\mcirc (Id\times\otimes):\2C\times\2C\times\2C\rarr\2C$, 
where we have used the strictness of $\2C$.

\bprop \label{p-coassoc}
Let $\2C, \2D$ be tensor $*$-categories and let $E: \2C\rarr\2D$ be a $*$-preserving tensor
functor. There are unital $*$-homomorphisms 
$\widetilde{\D\otimes\iota}, \widetilde{\iota\otimes\D}: \ \Nat\,E\otimes\rarr\Nat\,E\otimes^2$
defined by 
\begin{eqnarray*} 
   \widetilde{\D\otimes\iota}: \ \Nat\,E\otimes\rarr\Nat\,E\otimes^2, &&  b=(b_{X,Y})\ \mapsto\      (\widetilde{\D\otimes\iota})(b)=((\widetilde{\D\otimes\iota})(b)_{U,V,W}), \quad
    (\widetilde{\D\otimes\iota})(b)_{U,V,W}=b_{U\otimes V,W}, \\
   \widetilde{\iota\otimes\D}: \ \Nat\,E\otimes\rarr\Nat\,E\otimes^2, &&
   b=(b_{X,Y})\ \mapsto\ (\widetilde{\iota\otimes\D})(b)=((\widetilde{\iota\otimes\D})(b)_{U,V,W}), \quad
    (\widetilde{\iota\otimes\D})(b)_{U,V,W}=b_{U,V\otimes W}.
\end{eqnarray*}
These morphisms satisfy $(\widetilde{\D\otimes\iota})\tilde{\D}=(\widetilde{\iota\otimes\D})\tilde{\D}$.
\eprop
\prf That $(\widetilde{\D\otimes\iota})(b)\in\Nat\,E\otimes^2$ for $b\in\Nat\,E\otimes$ follows from
the following computation:
\begin{eqnarray*} \lefteqn{(E\otimes^2) (s\times t\times u) \, \mcirc  (\widetilde{\D\otimes\iota})(b)_{Y,Z,W} = 
   E(s\otimes t\otimes u) \, \mcirc  \, b_{Y\otimes Z,W}  } \\ 
  && = b_{Y'\otimes Z',W'} \, \mcirc  E(s\otimes t\otimes u) = (\widetilde{\D\otimes\iota})(b)_{Y',Z',W'} \mcirc (E\otimes^2) (s\times t\times u), 
\end{eqnarray*}
whenever $s: U\rarr U', t: V\rarr V', u: W\rarr W'$. The argument for
$(\widetilde{\iota\otimes\D})(b)$ is completely analogous. Clearly, both maps
$(\widetilde{\D\otimes\iota}), (\widetilde{\iota\otimes\D})$ are unital $*$-homomorphisms. 
The coassociativity property follows from
\[  (\widetilde{\D\otimes\iota})\tilde{\D}(b)_{U,V,W}=\tilde{\D}(b)_{U\otimes V,W}=b_{(U\otimes V)\otimes W}
   =b_{U\otimes (V\otimes W)}=\tilde{\D}(b)_{U,V\otimes W}=(\widetilde{\D\otimes\iota})\tilde{\D}(b)_{U,V,W} \]
with $b\in\End\,E$, where we have used the strictness of $\2C$.
\qed

\bprop \label{p-ve}
There are unital $*$-homomorphisms $\tilde{\ve}: \ \Nat\,E\rarr\End(E(\11))$ and 
$\widetilde{\ve\otimes\iota}, \widetilde{\iota\otimes\ve}:\ \Nat\,E\otimes\rarr\Nat\,E$ defined by
\[
\begin{array}{lll}
  \DS \tilde{\ve}: & \Nat\,E\rarr\End(E(\11)), & b=(b_X)\ \mapsto\ b_\11, \vspace{4pt}\\
  \DS \widetilde{\ve\otimes\iota}: & \Nat\,E\otimes\rarr\Nat\,E, & b=(b_{X,Y})\ \mapsto\ (\widetilde{\ve\otimes\iota})(b)=
   ((\widetilde{\ve\otimes\iota})(b)_X), \quad (\widetilde{\ve\otimes\iota})(b)_X=b_{\11,X},   \vspace{4pt} \\
  \DS \widetilde{\iota\otimes\ve}: & \Nat\,E\otimes\rarr\Nat\,E, & b=(b_{X,Y})\ \mapsto\ (\widetilde{\iota\otimes\ve})(b)=
   ((\widetilde{\iota\otimes\ve})(b)_X), \quad (\widetilde{\iota\otimes\ve})(b)_X=b_{X,\11}. 
\end{array}
\]
They satisfy 
$(\widetilde{\ve\otimes\iota})\tilde{\D}=(\widetilde{\iota\otimes\ve})\tilde{\D}=\iota$.
\eprop
\prf Clearly, the map $\tilde{\ve}: b=(b_X)\ \mapsto\ b_\11\in\End\,E(\11)$ is a unital
$*$-homomorphism. The proofs that 
$(\widetilde{\ve\otimes\iota})(b),(\widetilde{\iota\otimes\ve})(b)$ are in fact in $\Nat\,E$ and  
that the maps $b\mapsto(\widetilde{\ve\otimes\iota})(b), \ b\mapsto(\widetilde{\iota\otimes\ve})(b)$
are $*$-homomorphisms proceed along the same lines as for coproducts and are therefore omitted. For
$b\in\Nat\,E$ the computation 
\[ (\widetilde{\ve\otimes\iota})\tilde{\D}(b)= \tilde{\D}(b)_{\11,X}=b_{\11\otimes X}=b_X
  =b_{X\otimes\11}= \tilde{\D}(b)_{X,\11} = (\widetilde{\iota\otimes\ve})\tilde{\D}(b),  \]
where we used the strictness of the unit $\11$, concludes the proof.
\qed

We have thus seen that we can define a `coproduct' $\tilde{\D}: \Nat\,E\rarr\Nat\,E\otimes$ and a
`counit' $\tilde{\ve}: \Nat\,E\rarr\End E(\11)$ satisfying analogues of the usual properties. A
priori, however, $\Nat\,E\otimes$ has nothing to do with $\Nat\,E\otimes\Nat\,E$, regardless of how
we interpret the tensor product. The following result indicates how to proceed. 
We now assume that $\2C$ is semisimple. As before, by $I_\2C$ we denote the set of isomorphism
classes of irreducible objects of $\2C$. For every $i\in I_\2C$ we pick a representative $X_i$. 

\bdefin
Given a tensor $*$-category $\2C$, an embedding functor (for $\2C$) is a faithful $*$-preserving
tensor functor $E: \2C\rarr\2H$. A concrete tensor $*$-category is a tensor $*$-category together
with an embedding functor. 
\edefin

\bprop \label{p-canon-isom}
Let $\2C$ be a semisimple tensor $*$-category and $E: \2C\rarr\2H$ an embedding
functor. With the $*$-algebra $A=\bigoplus_{i\in I_\2C} B(E(X_i))$ we have the following
canonical $*$-isomorphisms: 
\[
\begin{array}{lclcl}
  \DS \Nat\,E & \stackrel{\psi_1}{\longrightarrow}& M(A) & = & \prod_{i\in I_\2C}
     B(E(X_i)), \vspace{5pt}\\ 
  \DS \Nat\,E \otimes & \stackrel{\psi_2}{\longrightarrow}& M(A\otimes A) & = & \prod_{i,j\in I_\2C} B(E(X_i))\otimes                   B(E(X_j)), \vspace{5pt}\\
  \DS \Nat\,E\otimes^2 &\stackrel{\psi_3}{\longrightarrow} & M(A\otimes A\otimes A) & =& \prod_{i,j,k\in I_\2C}          B(E(X_i))\otimes B(E(X_j))\otimes B(E(X_k)).
\end{array}
\]
\eprop
\prf The isomorphism $\psi_1$ was established in Propositions \ref{p-direct} and \ref{Berkeley} and
suppressed subsequently. As to $\psi_2$, we have
\[\begin{array}{cclclcc} 
  \DS \Nat\,E\otimes  & = & \prod_{i,j\in I_\2C} \End\,E(X_i\otimes X_j) 
  & \cong & \prod_{i,j\in I_\2C} \End\,(E(X_i)\otimes E(X_j))  \\ 
\DS  & = & \prod_{i,j\in I_\2C} B(E(X_i)\otimes E(X_j)) &=& \prod_{i,j \in I_\2C} 
      B(E(X_i))\otimes B(E(X_j))\ &=& M(A\otimes A). 
\end{array}
\]
Here we used the identification $B(H\otimes K)\equiv B(H)\otimes B(K)$ for finite dimensional
Hilbert spaces and the identifications made after Propositions \ref{Berkeley} and
\ref{p-direct}. The isomorphism between $\End\,E(X_i\otimes X_j)$ and $\End\,E(X_i)\otimes E(X_j)$
is induced in the canonical way from the isomorphism 
$d^E_{X_i,X_j}: E(X_i)\otimes E(X_j)\rarr E(X_i\otimes X_j)$ which makes $E$ monoidal. 
Unitarity of the $d^E_{X_i,X_j}$ implies that $\psi_2$ is $*$-preserving.
The argument for $\psi_3$ is completely analogous.
\qed

\brem In the above proof it was crucial that the canonical $*$-homomorphisms 
$\otimes: \ \End\,X\otimes_\7C\End\,Y \rarr \End(X\otimes Y)$ are $*$-isomorphisms for all
$X,Y\in\2H$. A semisimple tensor $*$-category has this property iff every object is a finite direct
sum of copies of the tensor unit $\11$. [{\it Proof.} The `if' clause is easy. As to the `only if'
part, let $Z$ be irreducible and consider $X\cong Y\cong \11\oplus Z$, thus 
$X\otimes Y\cong\11\oplus Z\oplus Z\oplus Z^2$. If $Z\not\cong\11$ then 
$\End\,X\otimes_\7C\End\,Y\cong\7C^2\otimes\7C^2$ is commutative whereas $\End(X\otimes Y)$ is
non-commutative since $X\otimes Y$ contains the irreducible $Z$ with multiplicity at least two. This
contradicts $\End\,X\otimes_\7C\End\,Y \cong\End(X\otimes Y)$.] This property is known to
characterize the (semisimple) tensor $*$-category $\2H$ uniquely up to equivalence. We thus see that
with $E=Id :\2C\rightarrow\2C$, we cannot conclude that $\Nat\,E\otimes \ \cong \ M(A\otimes A)$. 
This is consistent with the computation of $\Nat\,\otimes$ in Example \ref{e-otimes}.
\erem

From now on we restrict ourselves to the case we are ultimately interested in. 

\bdefin \label{d-ve-D}
Let $\2C$ be a semisimple tensor $*$-category and let $E: \2C\rarr\2H$ be an embedding functor.
Using the preceding results we define unital $*$-homomorphisms as follows.
\[\begin{array}{ll}
  \DS \overline{\ve}: \ M(A)\rarr\7C, & \quad a\mapsto \tilde{\ve}\mcirc \psi_1^{-1}(a), \vspace{4pt}\\
  \DS \overline{\D}: \ M(A)\rarr M(A\otimes A), & \quad a\mapsto \psi_2\mcirc \tilde{\D}\mcirc \psi_1^{-1}(a).
\end{array}\]
Here we have implicitly used the irreducibility of $E(\11_\2C)\cong\11_\2H\cong\7C$, giving rise to
the isomorphism $\psi_0:\End\,E(\11)\rarr\7C,\ c\,\id_{E(\11)}\mapsto c$.
\edefin

\brem 
1. In the sequel we will usually use the isomorphisms $\psi_1, \psi_2$ to identify
$M(A)\equiv\Nat\,E,\ M(A\otimes A)\equiv\Nat\,E\otimes$ and suppress the symbols $\psi_1, \psi_2$.  
Thus for $a\in M(A), i\in I,\  a_i$ will denote both the $i$-component of $a$ in
$\prod_{i\in I} B(E(X_i))$ and the action $a_{X_i}$ on $E(X_i)$ of the natural
transformation $\psi_1^{-1}(a)\in\Nat\,E$.  

2. Let $a\in M(A)$. As $\D(a)_{ij}=a_{X_i\otimes X_j}$ and $(a\otimes a)_{ij}=a_i\otimes a_j$ we
find
\[ \D(a)=a\otimes a \quad\Leftrightarrow\quad a_{X_i\otimes X_j}=a_{X_i}\otimes a_{X_j} \quad\forall i,j\in I. \]
Thus the grouplike elements of $M(A)$ are in one-to-one correspondence with those natural
transformations of the embedding functor which respect the tensor structure.
\erem

\blemma \label{l-nondeg}
The restriction $\D=\overline{\D}\restr A$ is a non-degenerate $*$-homomorphism from $A$ to
$M(A\otimes A)$, i.e.\ $\D(A)(A\otimes A)=(A\otimes A)\D(A)=(A\otimes A)$. 
\elemma
\prf To simplify the notation we write $H_i=B(E(X_i))$ and use the isomorphisms $\psi_1, \psi_2$
to identify $M(A)$ with $\Nat\,E$ and $M(A\otimes A)$ with $\Nat\,E\otimes$. We will only prove
$\D(A)(A\otimes A)=(A\otimes A)$, the proof of the other identity being completely analogous. 
As $A=\oplus_i B(H_i)$, it suffices to show that for every $i,j\in I_\2C$ there is a $c\in A$ with 
$\D(c)(I_i\otimes I_j)=I_i\otimes I_j$, since this implies 
$a_i\otimes b_j=\D(c)(I_i\otimes I_j)(a_i\otimes b_j)\in\D(A)(A\otimes A)$ for all 
$a_i\in B(H_i), b_j\in B(H_j)$. Recall from the
Propositions \ref{p-canon-isom} and \ref{p-direct} that $M(A)$ is isomorphic to the $*$-algebra of
natural transformations of $E$, and these are uniquely determined by their actions on the
irreducibles: $M(A)\ni b\leftrightarrow (b_i)$, where $b_i\equiv b_{X_i}$. Similarly, 
$M(A\otimes A)\ni b\leftrightarrow (b_{ij})$ with $b_{ij}\equiv b_{X_i,X_j}$.  
By the definition of $\tilde{\Delta}$ we have $\tilde{\Delta}(b)_{ij}=b_{X_i\otimes X_j}$.
As usual, let $(v_{ij}^{k\alpha}, \alpha=1,\ldots,N_{ij}^k)$ be an orthonormal basis in
$\Mor(X_k,X_i\otimes X_j)$ for all $i,j,k\in I_\2C$. Since $b\in M(A)\equiv\End\,E$ is a natural
transformation we have 
\[ \tilde{\Delta}(b)_{ij}=b_{X_i\otimes X_j}= b_{X_i\otimes X_j}\mcirc \id_{E(X_i\otimes X_j)}=
   b_{X_i\otimes X_j}\mcirc E\left(\sum_{k\alpha} v_{ij}^{k\alpha}\mcirc  v_{ij}^{k\alpha*}\right)
   = \sum_{k\alpha} E(v_{ij}^{k\alpha})\mcirc b_k\mcirc E(v_{ij}^{k\alpha*}). \]
Now fix $i_0, j_0\in I_\2C$ and let $b_k=\id_{E(X_k)}$ if $X_k$ is contained in 
$X_{i_0}\otimes X_{j_0}$ and $b_k=0$ otherwise. Since $X_k\prec X_{i_0}\otimes X_{j_0}$ only for
finitely many $k$ we have $b=(b_k)\in A$ and
\[ \tilde{\Delta}(b)_{i_0j_0}=\sum_{k\alpha} E(v_{i_0j_0}^{k\alpha})\mcirc \id_{E(X_k)}\mcirc
   E(v_{i_0j_0}^{k\alpha*}) =\sum_{k\alpha} E(v_{i_0j_0}^{k\alpha}\mcirc  v_{i_0j_0}^{k\alpha*})
  = E(\id_{X_{i_0}\otimes X_{j_0}}) = \id_{E(X_{i_0}\otimes X_{j_0})}. \]
Having identified $\Nat\,E\otimes\equiv M(A\otimes A)$ this precisely means that the
$(i_0,j_0)$-component of $\D(b)\in M(A\otimes A)=\prod_{ij} B(H_i)\otimes B(H_j)$ is the identity
$I_{i_0}\otimes I_{j_0}$ and therefore $\Delta(b)(I_{i_0}\otimes I_{j_0})=(I_{i_0}\otimes I_{j_0})$
as desired.
\qed

\bcoro \label{c-extens}
Defining $*$-homomorphisms $\ve: A\rarr\7C$ and $\D: A\rarr M(A\otimes A)$ by
$\ve=\overline{\ve}\restr A$ and $\D=\overline{\D}\restr A$, $\overline{\ve}$ and $\overline{\D}$
are the unique extensions of $\ve$ and $\D$, respectively, to $M(A)$. 
\ecoro
\prf It is clear that $\ve$ is not identically zero, thus non-degenerate. Non-degeneracy of $\D$ has
been proven in Lemma \ref{l-nondeg}. Now the claim follows from the fact \cite{VD2} that a
non-degenerate $*$-homomorphism $\phi: A\rarr M(B)$ has a unique extension to a unital
$*$-homomorphism $\overline{\phi}: M(A)\rarr M(B)$.
\qed

\bprop
\label{p-cdisc}
Let $\2C$ be a semisimple tensor $*$-category and $E:\2C\rarr\2H$ an embedding functor.
Then the $*$-algebra $A=\oplus_{i\in I_\2C} B(E(X_i))$ is a discrete quantum semigroup, i.e.\ has a  
non-degenerate coproduct $\D: A\rarr M(A\otimes A)$ in the sense of \cite{VD2} and a counit 
$\ve: A\rarr\7C$.
\eprop
\prf By Lemma \ref{l-nondeg}, $\D$ is non-degenerate, thus also 
$\D\otimes\iota: A\otimes A\rarr M(A\otimes A)\otimes A\subset M(A\otimes A\otimes A)$ is
non-degenerate and extends uniquely to 
$\overline{\D\otimes\iota}: M(A\otimes A)\rarr M(A\otimes A\otimes A)$. Another unital
$*$-homomorphism $M(A\otimes A)\rarr M(A\otimes A\otimes A)$ is given by 
$\psi_3\mcirc (\widetilde{\D\otimes\iota})\mcirc \psi_2^{-1}$. The latter clearly restricts to 
$\D\otimes\iota$ on $A\otimes A$, thus it coincides with $\overline{\D\otimes\iota}$ on all of
$M(A\otimes A)$. In the same way we make sense of 
$\overline{\iota\otimes\D}: M(A\otimes A)\rarr M(A\otimes A\otimes A)$ and see that it coincides
with $\psi_3\mcirc (\widetilde{\iota\otimes\D})\mcirc \psi_2^{-1}$. Now, the coassociativity
property proved in Proposition \ref{p-coassoc} implies
\[ \overline{\D\otimes\iota}\mcirc \overline{\D}= \overline{\iota\otimes\D}\mcirc \overline{\D}:\
    M(A)\rarr M(A\otimes A\otimes A), \] 
which is the proper formulation of coassociativity for $(A,\D)$.

As observed earlier, the $*$-homomorphism $\ve: A\rarr\7C$ is non-degenerate, thus 
$\ve\otimes\iota,\, \iota\otimes\ve:\ A\otimes A\rarr A$ are non-degenerate too and have unique
extensions $\overline{\ve\otimes\iota},\, \overline{\iota\otimes\ve}$ to $M(A\otimes A)$. Again, it
is clear that 
$\overline{\ve\otimes\iota}=\psi_1\mcirc (\widetilde{\ve\otimes\iota})\mcirc \psi_2^{-1}$ and
similarly for $\overline{\iota\otimes\ve}$. Thus with Proposition \ref{p-ve} we conclude
\[ \overline{\ve\otimes\iota}\mcirc \overline{\D}=\overline{\iota\otimes\ve}\mcirc
   \overline{\D}=\id: \ M(A)\rarr M(A). \] 
We conclude that $(A,\D)$ is a discrete quantum semigroup.
\qed

\btheor \label{t-mon-equiv}
Let $\2C$ be a semisimple tensor $*$-category and $E:\2C\rarr\2H$ an embedding functor. 
Let $(A,\D)$ be the corresponding discrete quantum semigroup. Then the functor 
$F: \2C\rarr\Rep_f(A,\D)$ defined by 
\[ F(X)=(E(X),\pi_X), \quad \pi_X(a)=a_X \ \forall a\in A, \quad\quad F(s)=E(s) \]
provides an equivalence of tensor $*$-categories such that $K\mcirc  F= E$.
\etheor
\prf With $A$ and $F: \2C\rarr\Rep_f A$ as defined in Theorem \ref{t-equiv} we have an equivalence
of categories. It remains to show that $F: \2C\rarr\Rep_f(A,\D)$ is monoidal for the monoidal
structure of $\Rep_f(A,\D)$ in Definition \ref{d-repf}. Thus we must exhibit isomorphisms
$d^F_{X,Y}: F(X)\otimes F(Y)\rarr F(X\otimes Y)$ satisfying the conditions of Definition
\ref{d-mon-func}. In view of $F(X)=(E(X),\pi_X)$, an obvious candidate is the isomorphism 
$d^E_{X,Y}: E(X)\otimes E(Y)\rarr E(X\otimes Y)$ which we have by virtue of the assumed monoidality
of $E$. If we set $d^F_{X,Y}=d^E_{X,Y}$ for all $X,Y\in \2C$, the equations in Definition
\ref{d-mon-func} are clearly satisfied. It only remains to prove that
$d^F_{X,Y}$ is in $\Mor_{\Rep_f(A,\D)}(F(X)\otimes F(Y), F(X\otimes Y))$, i.e.
\[ d^F_{X,Y} \mcirc  (\pi_X\times\pi_Y)(a)=\pi_{X\otimes Y}(a)\mcirc  d^F_{X,Y} \quad\forall
   a\in A. \] 

Keeping in mind the various identifications we have made so far, we get
\[d^F_{X,Y} \mcirc  (\pi_X\times\pi_Y)(a)=d^F_{X,Y} \mcirc  \overline{(\pi_X\otimes\pi_Y)}\D (a)
=d^F_{X,Y} \mcirc  \overline{(\pi_X\otimes\pi_Y)}\overline{\D }(a)\]
\[=d^F_{X,Y} \mcirc  \overline{(\pi_X\otimes\pi_Y)}\psi_2\tilde{\D}(a)
=\pi_{X,Y}\tilde{\D}(a)\mcirc  d^F_{X,Y}
=a_{X\otimes Y}\mcirc  d^F_{X,Y} =\pi_{X\otimes Y}(a)\mcirc  d^F_{X,Y} \]
for all $a\in A$. (The reader is invited to check the required identity explicitly for 
$X=X_i, Y=X_j$, as suffices by naturality of the isomorphisms $d^F$.)  
\qed


\subsection{Conjugates and Antipodes: Discrete Quantum Groups}
We have already seen that the concrete semisimple tensor $*$-categories
$\Rep_f (A,\D )$ and $\Corep_f (\hat{A},\hat{\D})$ come with conjugates
whenever $(A,\D )$ is a discrete aqg. In this section we prove the converse.
We start off with a couple of preparatory results.

\blemma
\label{l-inclusions}
Let $\2C$ be a semisimple tensor $*$-category with conjugates and let $E :\2C\rightarrow\2H$ be an
embedding functor. Then the coproduct $\D :A\rarr M(A\otimes A)$ for the discrete quantum semigroup
$(A,\D )$ given by Proposition \ref{p-cdisc} has the following inclusion properties
\[ (1\otimes A)\D (A)\subset A\otimes A, \quad\quad (A\otimes 1)\D (A)\subset A\otimes A. \]
\elemma

\prf 
We prove only the inclusion $(A\otimes 1)\D (A)\subset A\otimes A$, the proof of the inclusion
$(1\otimes A)\D (A)\subset A\otimes A$ being completely analogous. We need to show that 
$(a\otimes 1)\D(b)\in A\otimes A$ for all $a,b\in A$. It is clearly sufficient to show this 
for $b=(b_i)$ where all components except $i=m$ vanish and for $a=I_n$.
By the calculation done in the proof of Lemma \ref{l-nondeg} we have
\[ \D(b)_{ij}=b_{X_i\otimes X_j} = \sum_{k\in I_\2C} \sum_{\alpha=1}^{N_{ij}^k} 
   E(V_{ij}^{k\alpha})\mcirc  b_k\mcirc E(V_{ij}^{k\alpha*})\ \in\End\,E(X_i\otimes X_j). \]
In view of $a=I_n$ and $b_i=0$ for $i\ne m$ we obtain
\[ ((a\otimes 1)\D(b))_{ij}=\delta_{i,n} \sum_{\alpha=1}^{N_{nj}^m} E(V_{nj}^{m\alpha})\mcirc     
 b_m\mcirc E(V_{nj}^{m\alpha*}). \]
Now, since $\2C$ has conjugates we have isomorphisms 
$\Mor(X_n\otimes X_j,X_m)\cong\Mor(X_j,\overline{X_n}\otimes X_m)$. By semisimplicity, for fixed
$n,m$ the latter space is non-zero only for finitely many $j\in I_\2C$. Thus 
$((a\otimes 1)\D(b))_{ij}$ is non-zero only for $i=n$ and finitely many $j$ and therefore 
$(a\otimes 1)\D(b)\in A\otimes A$. 
\qed

\bprop
Let $\2C$ be a semisimple tensor $*$-category with conjugates and let $E :\2C\rightarrow\2H$ be an
embedding functor. Let $b=(b_{X,Y} )\in \Nat\,E\otimes$. Using the identifications 
$B((E(X\otimes Y))=B(E(X))\otimes B(E(Y))$, for any objects $X$ and $Y$ in $\2C$, we may write
$b_{X,X}=\sum_k a^k_X\otimes b^k_X\in B(E(X))\otimes B(E(X))$, and then define
$\tilde{m}(b)_X\in B(E(X))$ by
\[\tilde{m}(b)_X =\sum_k a^k_X\mcirc  b^k_X, \]
for any object $X\in\2C$. This gives a unital linear map $\tilde{m}:\Nat\,E\otimes \rarr\Nat\,E$
which sends $b\in Nat\,E\otimes$ to $\tilde{m}(b)$. 
Furthermore, it restricts to the linearized multiplication $m :A\otimes A\rightarrow A$ on 
$A\subset\Nat\,E$. 
\eprop

\prf The bilinear multiplication map from $B(E(X))\otimes B(E(X))$ to $B(E(X))$ exists since $E(X)$
is finite dimensional. The rest is obvious.
\qed

\bprop \label{p-coinv}
Let $\2C$ be a semisimple tensor $*$-category with conjugates and let $E :\2C\rightarrow\2H$ be an
embedding functor. Then the discrete quantum semigroup $(A,\D )$ given by Proposition \ref{p-cdisc}
is a discrete aqg. The coinverse $S:A\rarr A$ is given by
\[ S(a)_X =E(\id_X\otimes r^*_X )\mcirc \id_{E(X)}\otimes a_{\overline{X}}\otimes \id_{E(X)}
   \mcirc  E(\overline{r}_X\otimes \id_X ) \]
for all $a\in A\subset\Nat\,E$ and objects $X$ in $\2C$, and where
$(r_X, \overline{r}_X, \overline{X})$ is any solution of the conjugate equations associated to $X$. 
\eprop 

\prf
In view of Proposition \ref{p-converse} it is sufficient to show that 
$(A\otimes 1)\D (A)\subset A\otimes A$ and $(1\otimes A)\D (A)\subset A\otimes A$ and to produce an
invertible  coinverse $S$. We thus get an aqg since Theorem \ref{thm-haar} provides a (non-zero)
left-invariant positive functional for discrete multiplier Hopf $*$-algebras. Now the inclusions
$(A\otimes 1)\D (A)\subset A\otimes A$ and $(1\otimes A)\D (A)\subset A\otimes A$ are the content of
Lemma \ref{l-inclusions}.

We show that the map $S:A\rightarrow A$ given by the formula in the theorem is well-defined.
We first show that for $a\in A$ and $X\in\2C$, the formula for
$S(a)_X \in\End\,E(X)$ does not depend on the choice of 
solution $(r_X, \overline{r}_X, \overline{X})$ of the conjugate equations associated to $X$.   
Suppose therefore that $(r'_X, \overline{r}'_X, \overline{X}')$ is another solution
of the conjugate equations associated to $X$. By \cite{LR} there exists a unique invertible 
$t\in\Hom(\overline{X},\overline{X}')$ such that
$r'_X =t\otimes \id_X\mcirc  r_X$ and $\overline{r}'_X =\id_X\otimes t^{*-1}\mcirc  \overline{r}_X$.
Thus
\begin{eqnarray*} \lefteqn{
 E(\id_X\otimes {r'}^*_X )\mcirc  \id_{E(X)}\otimes a_{\overline{X}'}\otimes \id_{E(X)}
      \mcirc  E(\overline{r}'_X\otimes \id_X ) } \\
 && =E(\id_X\otimes (r^*_X \mcirc  t^*\otimes \id_X ))\mcirc  \id_{E(X)}\otimes
      a_{\overline{X}'}\otimes \id_{E(X)} 
    \mcirc  E((\id_X\otimes t^{*-1}\mcirc \overline{r}_X)\otimes \id_X ) \\
 &&=E(\id_X\otimes r^*_X )\mcirc  \id_{E(t)}\otimes (E(t^* )\mcirc  a_{\overline{X}'}\mcirc 
      E(t^{*-1})) \otimes \id_{E(X)}\mcirc  E(\overline{r}_X\otimes \id_X ) \\
 &&=E(\id_X\otimes r^*_X )\mcirc  \id_{E(X)}\otimes (E(t^*\mcirc t^{*-1}\mcirc  a_{\overline{X}})
    \otimes \id_{E(X)})\mcirc  E(\overline{r}_X\otimes \id_X ) \\
 &&=E(\id_X\otimes r^*_X )\mcirc  (\id_{E(X)}\otimes a_{\overline{X}}
    \otimes \id_{E(X)})\mcirc  E(\overline{r}_X\otimes \id_X ),
\end{eqnarray*}
and $S(a)_X \in\End\,E(X)$ does not depend on the choice of 
solution $(r_X, \overline{r}_X, \overline{X})$ of the conjugate equations associated with $X$.

Next we show that $S(a) =(S(a)_X )\in\Nat\,E$ for $a\in A$, i.e., we must show that 
$E(t)\mcirc  S(a)_X =S(a)_Y\mcirc  E(t)$, for all $a\in A$ and $t\in\Hom(X,Y)$.
By the existence of the transpose $t^{\bullet *}\in\Hom(Y,X)$ given in \cite{LR}, we get
\begin{eqnarray*}  
E(t)\mcirc  S(a)_X 
&=& E(\id_Y\otimes r_{X}^* )\mcirc  E(t)\otimes a_{\overline{X}}
  \otimes \id_{E(X )}\mcirc  E(\overline{r}_X\otimes \id_X) \\
&=& E(\id_Y\otimes r_{X}^* )\mcirc  \id_{E(Y )}\otimes 
  a_{\overline{X}}\otimes \id_{E(X )}\mcirc  
  E((t\otimes\id_{\overline{X}}\mcirc \overline{r}_X)\otimes \id_X) \\
&=& E(\id_Y\otimes r_X^* )\mcirc \id_{E(Y)}\otimes a_{\overline{X}}\otimes \id_{E(X )}\mcirc  
  E((\id_Y\otimes t^{\bullet *}\mcirc \overline{r}_Y)\otimes \id_X) \\
&=& E(\id_Y\otimes r_X^* )\mcirc \id_{E(Y)}\otimes (a_{\overline{X}}\circ E(t^{\bullet *}))
   \otimes \id_{E(X )}\mcirc E(\overline{r}_{Y}\otimes \id_X) \\
&=& E(\id_{Y}\otimes r_X^* )\mcirc \id_{E(Y )}\otimes (E(t^{\bullet *})\circ a_{\overline{Y}})
   \otimes \id_{E(X )} \mcirc  E(\overline{r}_Y\otimes \id_X) \\
&=& E(\id_Y\otimes (r_{X}^*\mcirc t^{\bullet *}\otimes \id_X ))
  \mcirc \id_{E(Y )}\otimes a_{\overline{Y}}\otimes \id_{E(X )}
  \mcirc  E(\overline{r}_Y\otimes \id_X) \\
&=& E(\id_Y\otimes (r_Y^*\mcirc \id_{\overline{Y}}\otimes t))
  \mcirc \id_{E(Y )}\otimes a_{\overline{Y}}\otimes \id_{E(X )}
  \mcirc  E(\overline{r}_Y\otimes \id_X) \\
&=& E(\id_Y\otimes r_Y^* )\mcirc  \id_{E(Y )}\otimes
  a_{\overline{Y}}\otimes E(t)\mcirc  E(\overline{r}_Y \otimes \id_Y ) \\
&=& S (a)_Y \mcirc  E(t),
\end{eqnarray*}
as desired.

Clearly, we have $S(a)\in A$ for all $a\in A$ and it is also clear that the map $S$ is linear.
We must prove the two identities
\[m(\iota\otimes S)((a\otimes 1)\D (b))=\varepsilon (b)a,\ \ \ \ \ \ \ 
m(S\otimes\iota )((1\otimes a)\D (b))=\varepsilon (b)a,\]
for all $a,b\in A$. We prove only the first one. Thus let $a,b\in A$ and write
$(a\otimes 1)\D (b)=\sum_k a^k\otimes b^k$ for some $a^k, b^k\in A$. 
This means that
\[ (a_X\otimes 1_Y )\mcirc  b_{X\otimes Y}=\sum_k a^k_X \otimes b^k_Y \]
for all objects $X$ and $Y$ of $\2C$.
Thus we get
\begin{eqnarray*} \lefteqn{
  m (\iota\otimes S)((a\otimes 1)\D (b))_X  \ \  =\ \  \sum_k m(a^k\otimes S(b^k ))_X \ \ 
  =\ \ \sum_k a^k_X\mcirc  S(b^k )_X } \\
&& =\sum_k a^k_X\mcirc  E(\id_X\otimes r^*_X )\mcirc \id_X\otimes b^k_{\overline{X}}\otimes\id_{E(X)}
  \mcirc  E(\overline{r}_X\otimes\id_X ) \\
&&=E(\id_X\otimes r^*_X )\mcirc  (\sum_k a^k_X\otimes b^k_{\overline{X}} )\otimes \id_{E(X)}
  \mcirc  E(\overline{r}_X\otimes \id_X ) \\
&& =E(\id_X\otimes r^*_X )\mcirc  (a_X\otimes \id_{E(\overline{X})}\mcirc  b_{X\otimes\overline{X}})
  \otimes\id_{E(X)}\mcirc  E(\overline{r}_X\otimes \id_X ) \\
&& =a_X\mcirc  E(\id_X\otimes r^*_X )\mcirc  b_{X\otimes\overline{X}}
  \otimes\id_{E(X)}\mcirc  E(\overline{r}_X\otimes \id_X ) \\
&& =a_X\mcirc  E(\id_X\otimes r^*_X )\mcirc  E(\overline{r}_X\otimes \id_X )b_{\11} \quad=\quad
  a_X\mcirc  E(\id_X\otimes r^*_X \mcirc \overline{r}_X\otimes \id_X )b_{\11} \\
&& =a_X\mcirc  E(\id_X )b_{\11} \quad=\quad a_X b_{\11} \quad=\quad a_X \varepsilon (b) \quad=\quad
(\varepsilon (b)a)_X
\end{eqnarray*}
for all objects $X$ in $\2C$, so $m(\iota\otimes S)((a\otimes 1)\D (b))=\varepsilon (b)a$ as natural 
transformations.

It remains to show that $S$ is invertible. We do that by showing $S(S(a)^* )=a^*$ for all $a\in A$.
Now with $X\in\2C$, we calculate
\begin{multline*} S (S (a)^* )_X=E(\id_{X}\otimes r_{X}^* )\mcirc  \id_{E(X )}
   \otimes S (a)^*_{\overline{X}}\otimes \id_{E(X )}
   \mcirc  E(\overline{r}_{X}\otimes \id_{X} ) \\
 =E(\id_{X}\otimes r_{X}^* )\mcirc  \id_{E(X )}\otimes
  [E(\id_{\overline{X}}\otimes \overline{r}_{X}^* )\mcirc 
  \id_{E(\overline{X})}\otimes a_X \otimes \id_{E(\overline{X})}
  \mcirc  E(r_{X}\otimes \id_{\overline{X}})]^*
  \otimes \id_{E(X )}\mcirc  E(\overline{r}_{X}\otimes \id_{X}), 
\end{multline*}
where we used $(\overline{r}_{X}, r_{X},X )$ as a solution of the conjugate equations for
$\overline{X}\in\2C$. Thus
\begin{eqnarray*}
\lefteqn{ S (S (a)^* )_X } \\
&=& E(\id_{X}\otimes r_{X}^* )
 \mcirc  \id_{E(X )}\otimes [E(r_{X}^*\otimes \id_{\overline{X}})
 \mcirc  \id_{E(\overline{X} )}\otimes a_X^*
 \otimes \id_{E(\overline{X} )}
 \mcirc  E(\id_{\overline{X}}\otimes \overline{r}_{X})]
 \otimes \id_{E(X )}\mcirc  E(\overline{r}_{X}\otimes \id_{X}) \\
&=& E(\id_{X}\otimes r_{X}^*\mcirc  \id_{X}\otimes r_{X}^*
 \otimes \id_{\overline{X}\otimes X})\mcirc  \id_{E(X\otimes\overline{X})}
 \otimes a_X ^*\otimes \id_{E(\overline{X}\otimes X)}\mcirc 
 E(\id_{X\otimes\overline{X}}\otimes \overline{r}_{X}\otimes \id_{X}
 \mcirc \overline{r}_{X}\otimes \id_{X}) \\
&=& E(\id_{X}\otimes r_{X}^*\otimes r_{X}^* )\mcirc 
 \id_{E(X\otimes\overline{X} )}\otimes a_X^*\otimes 
 \id_{E(\overline{X}\otimes X )}\mcirc  E(\overline{r}_{X}\otimes
 \overline{r}_{X}\otimes \id_{X}) \\
&=& E(\id_{X}\otimes r_{X}^*\otimes r_{X}^*\mcirc  
 \overline{r}_{X}\otimes \id_{X\otimes\overline{X}\otimes X})
 \mcirc  a_X^*\otimes \id_{E(\overline{X}\otimes X )}\mcirc 
 E(\overline{r}_{X}\otimes \id_{X}) \\
&=& E((\id_{X}\otimes r_{X}^* )\mcirc  (\overline{r}_{X}
 \otimes \id_{X})\otimes \id_{\overline{X}\otimes X}\mcirc  \id_{X}
 \otimes r_{X}^* )\mcirc  a_X^*\otimes 
 \id_{E(\overline{X }\otimes X )}\mcirc  E(\overline{r}_{X}\otimes \id_{X}) \\
&=& E(\id_{X}\otimes r_{X}^* )\mcirc  a_X^*\otimes 
 \id_{E(\overline{X}\otimes X )}\mcirc  E(\overline{r}_{X}\otimes \id_{X})
 =a_X^*\mcirc  E(\id_{X}\otimes r_{X }^*\mcirc \overline{r}_{X}
 \otimes \id_{X}) \\
&=& a_X^*\mcirc  \id_{E(X )} \ =\ a_X^* \ = \ (a^* )_X.
\end{eqnarray*}
\qed

After these preparatory computations, the following Main Result is essentially a restatement of  
Theorem \ref{t-mon-equiv} and Proposition \ref{p-coinv}.

\btheor \label{t-main-equiv}
Let $\2C$ be a semisimple tensor $*$-category with conjugates 
and $E:\2C\rarr\2H$ an embedding functor. Let $(A,\D)$ be the corresponding discrete aqg. 
Then the functor $F: \2C\rarr\Rep_f(A,\D)$ defined by
\[ F(X)=(E(X),\pi_X), \quad \pi_X(a)=a_X \ \forall a\in A,\quad F(s)=E(s) \]
provides an equivalence of tensor $*$-categories such that $K\mcirc  F= E$.
\etheor

\brem
1. In combination with Theorem \ref{t-dualcor} this gives Woronowicz' generalized Tannaka-Krein
result. Woronowicz constructs $\Corep_f (\hat{A},\hat{\D})$ for a compact aqg $(\hat{A},\hat{\D})$ 
more directly and uses the universal corepresentation $U$ of $(A,\D)$ implicitly. In his proof
Woronowicz constructs a compact matrix pseudogroup from a category with a generator. Our result is
more general in that this assumption is redundant.

2. In view of Theorem \ref{t-dual-comod}, an alternative reconstruction theorem for a concrete
semisimple tensor $*$-category $\2C$ with conjugates would be that
\[ \2C \stackrel{\otimes}{\cong}{\Comod}_f (\hat{A},\hat{\D}^{op}), \]
where $(\hat{A},\hat{\D}^{op})$ is a Hopf $*$-algebra with a positive invariant functional. We refer
to \cite{ulb, JS1} for similar constructions, where concrete non-semisimple tensor categories are
shown to be equivalent to tensor categories of finite dimensional comodules over (infinite
dimensional) Hopf algebras (with no $*$-operation).
\erem

We now discuss how the aqg $(A,\D)$ depends on the choice of the embedding functor $E$. This
discussion is only preliminary and will be taken up again in part II of this series. We need the
following 

\bdefin Let $\2C,\2D$ be strict tensor categories and $E,E':\2C\rarr\2D$ be tensor functors (with
structural isomorphisms $e, e', d_{X,Y}, d'_{X,Y}$). Then $E,E'$ are said to be isomorphic as tensor
functors, $E\stackrel{\otimes}{\cong}E'$, if the diagram
\begin{equation} \label{e-nat}
\begin{diagram} E(X)\otimes E(Y) & \rTo^{u_X\otimes u_Y} & E'(X)\otimes E'(Y) \\
\dTo_{d_{X,Y}} & & \dTo^{d'_{X,Y}} \\
E(X\otimes Y) & \rTo_{u_{X\otimes Y}} & E'(X\otimes Y)
\end{diagram} 
\end{equation}
commutes for all $X,Y\in\2C$, where $u: E\rarr E'$ is a natural transformation whose components
$u_X: E(X)\rarr E'(X)$ are isomorphisms.
A similar and obvious commutative diagram involves the structural isomorphisms $e,e'$. If $\2C,\2D$
are tensor $*$-categories and $E,E'$ are $*$-preserving then all $u_X$ must be unitary. 
\edefin

\bprop
Let $\2C$ be a semisimple tensor $*$-category $\2C$ with conjugates.
Let $E,E': \2C\rarr\2H$ be embedding functors, and let $(A,\D), (A',\D')$ be the 
corresponding discrete \aqgs. An isomorphism $u: E\stackrel{\otimes}{\cong} E'$ gives rise
to an isomorphism $\theta: (A,\D )\cong (A', \D ')$ with $\pi'_X\theta=\pi_X$ for all $X\in\2C$.
\eprop

\prf Let $u: E\rarr E'$ be a unitary equivalence of tensor $*$-functors. Then
$\theta :\Nat E\rightarrow\Nat E'$ given by $\theta (a)_X =u_Xa_Xu_X^*$ for $a\in\Nat E$ and
$X\in\2C$ restricts to a unital $*$-homomorphism from $A$ to $A'$. Suppressing as usual the natural
transformations $d$ and $d'$, 
we clearly have $\D'(\theta (a))_{(X,Y)}=u_{X\otimes Y}a_{X\otimes Y}u_{X\otimes Y}^*$, whereas
$(\theta\otimes\theta )\D(a)_{(X,Y)}=(u_X\otimes u_Y )a_{X\otimes Y}(u_X^*\otimes u_Y^* )$, for all
$a\in\Nat E$ and $X,Y\in\2C$. Since $u_{X\otimes Y}=u_X\otimes u_Y$, for all $X,Y\in\2C$, we get
$(\theta\otimes\theta )\D =\D'\theta$, and $(A,\D )\cong (A',\D' )$.
\qed

The following definition will be useful in our discussion of embedded symmetric tensor categories.

\bdefin \label{d-intr}
Let $\2C,\2D$ be tensor categories and $E: \2C\rarr\2D$ a tensor functor. Then the automorphism
group $\Aut^\otimes E$ consists of the natural isomorphisms of $E$ for which the diagram
(\ref{e-nat}) commutes (with $E'=E$), together with the obvious componentwise group structure.
If $\2C$ is a tensor $*$-category and $E$ is a $*$-preserving embedding functor then the intrinsic
group $G_E$ is the subgroup of $\Aut^\otimes E$ where all $u_X$ are unitaries.
\edefin

\subsection{The Non-$*$ Case}
Our starting point in this paper is a tensor $*$-category with direct sums, subobjects, conjugates 
and irreducible unit. From it we (re)construct a multiplier Hopf $*$-algebra with positive left
invariant functional. 
Such algebras were called algebraic quantum groups (\aqgs) by Van Daele, who now prefers to speak of
`multiplier Hopf $*$-algebras with integrals'. We stick to `algebraic quantum groups' mainly for the
sake of definiteness and brevity. When $A$ is no longer a $*$-algebra and the invariant functional
is required to be faithful rather than positive ($\varphi(ab)=0$ for all $a\in A$ implies $b=0$) we
arrive at the notion of a regular multiplier Hopf algebras with left invariant functional. 
Essentially the entire theory works just as well in this setting, provided one adopts a modified
notion of dual. We say a tensor category is rigid if every object $X$ comes with a chosen left dual
$\overline{X}$ and morphisms $e_X: \overline{X}\otimes X\rarr\11$ and 
$d_X: \11\rarr X\otimes\overline{X}$ satisfying the triangular equations instead of the conjugate
equations. In a semisimple category, left duals are automatically two-sided duals since 
dual morphisms $e'_X: \11\rarr\overline{X}\otimes X$ and $d'_X: X\otimes\overline{X}\rarr\11$ exist.
Adopting the standard definition of semisimple categories as abelian categories with all exact
sequences splitting, one can prove the following:

\btheor \label{t-nonstar}
Let $\7F$ be an algebraically closed field (of arbitrary characteristic). Let $\2C$ be a
$\7F$-linear semisimple rigid tensor category with $\End\,\11\cong\7F$ and let 
$E:\2C\rarr\Vect_\7F$ be a faithful tensor functor. Then there is a discrete regular 
multiplier Hopf algebra $(A,\D)$ with faithful left invariant functional $\varphi$ and an
equivalence $F: \2C\rarr\Rep(A,\D)$ such that $K\circ F=E$. 
\etheor

The proof proceeds exactly as for $*$-categories, ignoring all references to the $*$-preserving
property of $E$ and the $*$-involution of $A$. Again, existence, uniqueness and faithfulness of the
invariant functional are proven using results due to Van Daele. Only the construction of the
coinverse requires careful examination since the notion of dual is different.

We close with some remarks on the two different notions of discreteness applying to multiplier Hopf
algebras. In analogy to the $*$-case, a regular multiplier Hopf algebra $(A,\D)$ with left invariant
functional is called {\it discrete} if $A$ is a direct sum of finite dimensional matrix algebras
over $\7F$. In the non-$*$ case Proposition \ref{p-disc} fails, and $(A,\D)$ is said to
be of {\it discrete type} iff there exists $0\ne h\in A$ such that $ah=\ve(a) h$ for all $a\in A$. 
As in the $*$-case, $(A,\D)$ is of {\it compact type} or, equivalently, {\it compact} if $A$ has a
unit. It is well known \cite{VD} that $(A,\D)$ is of discrete type iff $(\hat{A},\hat{\D})$ is of
compact type. 

Every discrete multiplier Hopf algebra has a copy of $\7F$ as a
direct summand, corresponding to the one-dimensional representation $\ve$.  The unit $I_0$ of this
algebra has the property $aI_0=I_0a=\ve(a)I_0$, thus a discrete multiplier Hopf algebra is of
discrete type. The following characterizes discrete multiplier Hopf algebras, generalizing a well 
known result for finite dimensional Hopf algebras (which automatically possess integrals).

\bprop \label{p-disc2}
Let $(A,\D)$ and $(\hat{A},\hat{\D})$ be mutually dual regular multiplier Hopf algebras with
faithful left invariant functionals. Then the following are equivalent: 
\begin{enumerate} 
\item $(A,\D)$ is discrete.
\item $(A,\D)$ is of discrete type and $h$ can be chosen to be idempotent. 
\item $(\hat{A},\hat{\D})$ is of compact type and $\hat{\varphi}(1)\ne 0$.
\end{enumerate}
\eprop

\prf The implication 1.$\impl$2. has been discussed before. Assuming 2., we have $h=hh=\ve(h)h$,
thus $\ve(h)=1$. Conversely, if $\ve(h)\ne 0$ then $h/\ve(h)$ satisfies 2. Now, 2. implies that
$(\hat{A},\hat{\D})$ is compact, and compactness implies the existence of $0\ne h\in A$ satisfying
$xh=\ve(x)h$ for all $x\in A$.  By definition of the dual $(\hat{A},\hat{\D})$ and assuming the
normalization $\varphi(h)=1$, we have $\hat{\varphi}(1)=\ve(h)$, yielding the equivalence
2.$\Leftrightarrow$3. Now assume 3. holds. Then $(\hat{A},\hat{\D})$ is a Hopf algebra with left
invariant functional $\hat{\varphi}$ satisfying $\hat{\varphi}(1)\ne 0$. By classical results, cf.,
e.g., \cite[14.0.3]{swe} or \cite[Theorem 3.3.2]{abe}, it follows that $(\hat{A},\hat{\D})$ is
cosemisimple.  This means that the $\7F$-coalgebra $(\hat{A},\hat{\D},\hat{\ve})$ is a direct sum of
finite dimensional coalgebras. Since $\7F$ is algebraically closed, every finite dimensional
$\7F$-coalgebra is a matrix coalgebra, i.e.\ it has a basis $\{a_{ij}, 1\le i,j\le n\}$ such that 
$\D(a_{ij})=\sum_l a_{il}\otimes a_{lj}$. Since 
$A=\{ \hat{\varphi}(\cdot a),\ a\in \hat{A}\}\subset\hat{A}^*$ and
$\hat{\varphi}$ is faithful, $A$ has a basis 
$\{ b^\alpha_{ij}=\hat{\varphi}(\cdot a^\alpha_{ij}), \ \alpha\in I, 1\le i,j\le n_\alpha\in\7N\}$. 
These elements multiply according to
$b^\alpha_{ij}b^{\alpha'}_{i'j'}=\delta_{\alpha\alpha'}\delta_{ji'}b^\alpha_{ij'}$, thus 
\[ A=\bigoplus_{\alpha\in I} M_{n_\alpha}(\7F), \] as desired.  
\qed

\brem 1. In the $*$-case, $\hat{\varphi}$ is positive and $1>0$, thus discrete type implies
discreteness. 

2. We indicate how Proposition \ref{p-disc2} can be used to prove Theorem \ref{t-nonstar}. Given a 
semisimple tensor category $\2C$ with embedding functor $E$, by \cite{ulb} there exists a Hopf
algebra $H$ such that $\2C\simeq\Comod\,H$. By \cite[14.0.3]{swe}, semisimplicity of $\2C$ implies
the existence of a left invariant functional $\hat{\varphi}$ on $H$, satisfying 
$\hat{\varphi}(1)\ne 0$. Thus $(H,\hat{\Delta})$ is a compact multiplier Hopf algebra with left
invariant functional. Its dual $(A,\Delta)$ then is a discrete multiplier Hopf algebra with left
invariant functional, and by Theorem \ref{t-dualcor} we have $\2C\simeq\Rep_f(A,\Delta)$.
However, the direct proof of Theorem \ref{t-nonstar} analogous to the proof of
Theorem \ref{t-main-equiv} is more direct and instructive.
\erem


\section{Concrete Braided and Symmetric Categories} \label{s-braid}
In this section we prove that a braiding for the concrete tensor $*$-category $\2C$ gives rise to an
$R$-matrix for the associated aqg $(A,\D)$. If $\2C$ is symmetric and the embedding functor $E$ is 
{\it symmetric} monoidal then we may conclude that $(A,\D)$ is cocommutative.

The flip automorphism $\sigma$ of $A\otimes A$ defined by $\sigma(a\otimes b)=b\otimes a$ is a
non-degenerate $*$-homomorphism. Thus there exists a unique extension to an involutive
$*$-automorphism of $M(A\otimes A)$ which we denote by the same symbol.

\bdefin \label{d-quasit}
An aqg $(A,\D)$ is said to be quasitriangular w.r.t. $R$ iff $R$ is an invertible element of
$M(A\otimes A)$ and satisfies
\begin{enumerate}
\item $(\D\otimes\iota)(R)=R_{13}R_{23}$,
\item $(\iota\otimes\D)(R)=R_{13}R_{12}$,
\item $\D^{op} =R\D(\cdot )R^{-1}$,
\end{enumerate} 
where $\D^\op\equiv \sigma\D$. 
We sometimes refer to such an element $R$ as an $R$-matrix. 
If, in addition, $\sigma(R)=R^{-1}$ we call $(A,\D)$ triangular. 
\edefin

The first equation means that the invertible element $R\in M(A\otimes A)$ is a non-degenerate
corepresentation of $(A,\D)$ on $A$.
It is easy to see that $(A,\D )$ is quasitriangular w.r.t. $R$ iff $(A, \D^{op})$ is quasitriangular 
w.r.t. $\sigma (R)$.

\bdefin \label{d-braiding}
A braiding for a (strict) tensor category $\2C$ is a family of isomorphisms 
$c_{X,Y}: X\otimes Y\rarr Y\otimes X$, for all $X,Y\in\2C$ natural in both arguments
\[ t\otimes s\mcirc c_{X,Y}=c_{X',Y'}\mcirc s\otimes t \quad \forall s: X\rarr X', t: Y\rarr Y', \]
satisfying
\begin{eqnarray*} 
   c_{X,Y\otimes Z} &=& \id_Y\otimes c_{X,Z}\mcirc c_{X,Y}\otimes\id_Z, \\
   c_{X\otimes Y,Z} &=& c_{X,Z}\otimes \id_Y\mcirc \id_X\otimes c_{Y,Z}, \\
   c_{X,\11} &=& c_{\11,X}=\id_X
\end{eqnarray*} 
for all $X,Y,Z\in \2C$. A braiding is called a symmetry if, in addition, it satisfies
\[ c_{X,Y}\mcirc c_{Y,X}=\id_{Y\otimes X} \quad\forall X,Y\in\2C. \]
A tensor category with a (chosen) braiding or symmetry is called a braided or symmetric tensor
category, respectively. A braided functor $F: \2C\rarr\2D$ is a monoidal functor between strict
braided tensor categories making
\begin{equation} \label{e-braidfunc}
\begin{diagram}
F(X)\otimes F(Y) & \rTo^{d^F_{X,Y}} & F(X\otimes Y) \\
\dTo^{c_{F(X),F(Y)}} & & \dTo_{F(c_{X,Y})} \\
F(Y)\otimes F(X) & \rTo_{d^F_{Y,X}} & F(Y\otimes X)
\end{diagram}
\end{equation}
commutative for all $X, Y\in\2C$.
\edefin

The following result establishes the relation between $R$-matrices for $(A,\D)$ and braidings for
$\Rep_f(A,\D)$.

\bprop 
\label{p-braidR}
Let $(A,\D)$ be a discrete aqg. Let $R\in M(A\otimes A)$ and consider
\[  c_{(H,\pi),(H',\pi')} \equiv \Sigma_{H,H'} \, (\overline{\pi\otimes\pi'})(R), \quad (H,\pi),
    (H',\pi')\in\Rep_f(A,\D). \] 
Here $\Sigma_{H,K}$ is the canonical unitary flip from $H\otimes K\rarr K\otimes H$.
Then the family $(c_{(H,\pi),(H',\pi')})$ is a braiding for $\Rep_f(A,\D)$
iff $R$ is an $R$-matrix. Also $c_{(H,\pi),(H',\pi')}$ is unitary for all $(H,\pi),(H',\pi')$ iff
$R$ is unitary. If these equivalent conditions are satisfied then $(A,\D)$ is triangular iff the 
corresponding braiding for $\Rep_f(A,\D)$ is symmetric.
\eprop
\prf
Throughout the proof we write $c_{\pi, \pi'}$ for $c_{(H,\pi),(H',\pi')}$ and note that 
$c_{\pi, \pi'}: H\otimes H'\rightarrow H'\otimes H$, 
for all $(H,\pi ), (H' \pi ')\in \Rep_f(A,\D)$.
Also note that 
\[\overline{(\pi\otimes\pi')}\ \sigma (a) 
=\Sigma_{H',H}\ \overline{(\pi'\otimes\pi )}\ (a)\ \Sigma_{H, H'} 
\quad \forall\  (H,\pi), (H',\pi')\in\Rep_f(A,\D)\quad \forall\ a\in M(A\otimes A).\]
Now suppose $R\in M(A\otimes A)$ satisfies $\D^{op}(\cdot )R =R\D (\cdot)$. 
Then 
\begin{eqnarray*}  c_{\pi, \pi'}\ (\pi\times \pi' )(a) &=& \Sigma_{H,H'} \,
   (\overline{\pi\otimes\pi'})(R)\ (\overline{\pi\otimes\pi'})\D (a) 
   =\Sigma_{H,H'} \, (\overline{\pi\otimes\pi'}) (R\D (a))
   =\Sigma_{H,H'} \, (\overline{\pi\otimes\pi'}) (\D^{op} (a) R)  \\
 &=& \Sigma_{H,H'} \, (\overline{\pi\otimes\pi'}) \D^{op} (a)\ (\overline{\pi\otimes\pi'})(R)
   =\Sigma_{H,H'} \, (\overline{\pi\otimes\pi'}) (\sigma\D (a))\ \Sigma_{H',H}\ \Sigma_{H,H'}\
    (\overline{\pi\otimes\pi'})(R) \\ 
 &=& (\overline{\pi '\otimes\pi})\D (a)\ \Sigma_{H,H'}\ (\overline{\pi\otimes\pi'})(R)
    =(\pi'\times \pi )(a)\ c_{\pi, \pi'},
\end{eqnarray*}
so we have seen that $c_{\pi, \pi'} :\pi\times \pi '\rarr \pi'\times \pi$ for all 
$\pi, \pi '\in\Rep_f (A,\D )$. Conversely, by reversing the above calculations, we see that
$\D^{op}(\cdot ) R=R\D (\cdot )$ if $c_{\pi, \pi'} :\pi\times \pi '\rarr \pi'\times \pi$ for all
$\pi, \pi '\in\Rep_f (A,\D )$.

It is straightforward to see that $R$ is invertible iff 
$c_{\pi, \pi'}: \pi\times\pi'\rarr\pi'\times\pi$ is an isomorphism for all 
$\pi,\pi'\in\Rep_f(A,\D )$, that $R$ is unitary iff  
$c_{\pi, \pi'}$ is unitary for all $\pi, \pi '\in\Rep_f (A,\D )$, and that $R^{-1}=\sigma (R)$ iff 
$c_{\pi', \pi}\ c_{\pi, \pi '} =\id_{\pi, \pi'}$ for all $\pi, \pi '\in\Rep_f (A,\D )$.

Pick $\pi,\pi', \pi''\in\Rep_f (A,\D )$. Then 
\begin{eqnarray*} \lefteqn{  \overline{(\pi\otimes \pi'\otimes\pi'' )}(\iota\otimes\D)(R)
  =\overline{(\pi\otimes (\pi'\times\pi ''))}(R) } \\
 && =\Sigma_{H'\otimes H'', H}\ \Sigma_{H, H'\otimes H''}\ \overline{(\pi\otimes (\pi'\times\pi ''))}(R)
   =\Sigma_{H'\otimes H'', H}\ c_{\pi, \pi'\times \pi'' }, 
\end{eqnarray*}
whereas
\begin{eqnarray*} \lefteqn{ \overline{(\pi\otimes \pi'\otimes\pi'' )}(R_{13} R_{12})
  =(\Sigma_{H'', H}\ \Sigma_{H,H''}\ \overline{(\pi\otimes \pi'')} (R))_{13}\ 
  (\Sigma_{H', H}\ \Sigma_{H,H'}\ \overline{(\pi\otimes \pi')} (R))_{12} } \\
 && =(\Sigma_{H'', H}\ c_{\pi, \pi''} )_{13}\ (\Sigma_{H', H}\ c_{\pi, \pi'} )_{12}
  =(\Sigma_{H'', H})_{13}\ (c_{\pi, \pi''} )_{13}\ (\Sigma_{H', H})_{12}\ (c_{\pi, \pi'} )_{12} \\
 && =\Sigma_{H'\otimes H'', H}\ (c_{\pi, \pi''} )_{23}\ \ (c_{\pi, \pi'} )_{12}. 
\end{eqnarray*}
Hence $(\iota\otimes\D )(R)=R_{13}R_{12}$ iff
\[c_{\pi, \pi'\times \pi'' } =(c_{\pi, \pi''} )_{23}\ \ (c_{\pi, \pi'} )_{12}
=\id_{\pi '}\otimes c_{\pi,\pi ''}\mcirc  c_{\pi, \pi '}\otimes\id_{\pi ''}, \]
for all $\pi,\pi', \pi''\in\Rep_f (A,\D )$. Similarly, one shows that
$(\D\otimes\iota )R=R_{13}R_{23}$ iff
\[c_{\pi\times \pi', \pi'' } =(c_{\pi, \pi''} )_{12}\ \ (c_{\pi', \pi''} )_{23}
=c_{\pi,\pi ''}\otimes\id_{\pi'}\mcirc  \id_{\pi }\otimes c_{\pi', \pi ''}, \]
for all $\pi,\pi', \pi''\in\Rep_f (A,\D )$.
\qed

\btheor \label{t-braid1}
Let $\2C$ be a semisimple braided tensor $*$-category with conjugates and $E: \2C\rarr\2H$ a
faithful $*$-preserving tensor functor. Let $(A,\D)$ and $F: \2C\rarr\Rep_f(A,\D)$ be as constructed
in Theorem \ref{t-main-equiv}. Then $(A,\D)$ is quasitriangular w.r.t.\ a unique 
$R\in M(A\otimes A)$ such that the functor $F: \2C\rarr\Rep_f(A,\D)$ is an equivalence of braided
tensor categories w.r.t.\ the braiding for $\Rep_f(A,\D)$ provided by $R$.
$R$ is unitary iff $c_{X,Y}$ is unitary for all $X,Y$. $(A,\D)$ is triangular iff the 
corresponding braiding for $\Rep_f(A,\D)$ is symmetric.
\etheor
\prf For $i,j\in I_\2C$ we have 
\[ F(c_{X_i,X_j})\in\Mor_\2H(E(X_i)\otimes E(X_j), E(X_j)\otimes E(X_i))=B(E(X_i)\otimes E(X_j),
    E(X_j)\otimes E(X_i)). \] 
Thus $\Sigma_{E(X_j),E(X_i)}F(c_{X_i,X_j})\in B(E(X_i)\otimes E(X_j))$ and we can define
$R\in M(A\otimes A)=\prod_{ij}B(E(X_i)\otimes E(X_j))$ by
$R_{ij}=\Sigma_{E(X_j),E(X_i)}F(c_{X_i,X_j})$. 
Naturality of the braiding and functoriality of $F$ now imply that we have
\[ F(c_{X,Y})=\Sigma_{E(X),E(Y)} \, (\overline{\pi_X\otimes\pi_Y})(R) \quad \forall X,Y\in \2C. \]

Since $F$ is an equivalence, we see that the family $F(c_{X,Y})$ is a braiding for $\Rep_f (A,\D)$.
Denote this family by $(c_{(E(X),\pi_X),(E(Y),\pi_Y )})$, so 
\[c_{(E(X),\pi_X),(E(Y),\pi_Y )}=\Sigma_{E(X),E(Y)} \, (\overline{\pi_X\otimes\pi_Y})(R) \quad
    \forall X,Y\in \2C. \] 
By Proposition \ref{p-braidR}, the element $R\in M(A\otimes A)$ is therefore an $R$-matrix which is
unitary iff $c_{X,Y}$ is unitary for all $X,Y$ and satisfies $\sigma (R)=R^{-1}$ iff 
$c_{Y,X}\mcirc  c_{X,Y} =\id_{X\otimes Y}$ for all $X,Y$.    
With the braiding on $\Rep_f (A,\D)$ given by $R$ according to Proposition \ref{p-braidR}, we get by
definition $F(c_{X,Y})=c_{F(X),F(Y)}$, for all $X,Y$, as $F(X)=(E(X),\pi_X )$ for all $X$. Thus $F$
is a braided functor.
\qed

\brem It seems worthwhile comparing our result with \cite[Proposition XIII.1.4]{kas}, where it is
shown that the tensor category of left modules over a bialgebra $(H,m,\eta,\Delta,\ve)$ is braided
iff $H$ is quasitriangular w.r.t.\ some $R$-matrix. The proof uses the left regular representation
of $H$ and therefore would not apply to the category of finite dimensional $H$-modules if $H$ is
infinite dimensional. Our setting differs too in that $A$ is non-unital whenever it is infinite
dimensional. 
\erem

\brem As is well known, the category $\2H$ of Hilbert spaces has a unique braiding, the flip
$\Sigma_{H,K}: H\otimes K\rarr K\otimes H$. It is thus natural to ask when the functors $E$
respects this braiding. If so, then  
\[ E(c_{X,Y}\mcirc  c_{Y,X}) = E(c_{X,Y})\mcirc  E(c_{Y,X}) = \Sigma_{E(X),E(Y)}\mcirc \Sigma_{E(Y),E(X)}
   =\id_{E(Y)\otimes E(X)}= \id_{E(Y\otimes X)}=E(\id_{Y\otimes X}), \]
where we suppressed the isomorphisms $d^E_{X,Y}$ as is our policy throughout. If $E$ is faithful
this implies $c_{X,Y}\mcirc  c_{Y,X}=\id_{Y\otimes X}$, thus $\2C$ is symmetric. In other words, 
for non-symmetric categories $\2C$ there is no embedding functor satisfying
$E(c_{X,Y})=\Sigma_{E(X),E(Y)}$. For symmetric $\2C$ we have the following easy yet important
result.
\erem

\bcoro \label{c-symm}
Let $\2C$ be a semisimple symmetric tensor $*$-category and let $E: \2C\rarr\2H$ be a symmetric
embedding functor, i.e., $E(c_{X,Y})=\Sigma_{E(X),E(Y)}$ for all $X,Y\in\2C$. Then the corresponding
discrete aqg $(A,\D)$ is cocommutative and $F$ is a $*$-preserving symmetric monoidal equivalence.
\ecoro
\prf It is clear from the definition of the $R$-matrix $R$ in the proof of Theorem \ref{t-braid1}, 
that $R=1$ in this case, so $\D^{op}=\D$.
\qed

We are now in a position to re-prove the classical Tannaka-Krein duality theorem for compact groups.

\btheor 
Let $\2C$ be a semisimple tensor $*$-category with conjugates and unitary symmetry, and let 
$E: \2C\rarr\2H$ be a $*$-preserving faithful symmetric tensor functor. Then there exists a compact
group $G$ and an equivalence $F': \2C\rarr\Rep_fG$ of symmetric tensor $*$-categories such that
$K'\circ F'=E$, where $K': \Rep_fG\rarr\2H$ is the forgetful functor. The group $G$ is canonically
isomorphic to the intrinsic group $G_E$ of the embedding functor.
\etheor
\prf Under the given assumptions, Theorem \ref{t-braid1} provides us with a discrete aqg $(A,\D)$
and an equivalence $F: \2C\rarr\Rep_f(A,\D)$ of tensor $*$-categories such that $K\circ F=E$. By
assumption, $E$ preserves the symmetries, i.e.\ $E(c_{X,Y})=\Sigma_{E(X),E(Y)}$ for all
$X,Y\in\2C$. Then by Corollary \ref{c-symm} $(A,\D)$ is cocommutative. Let $G$ be the intrinsic
group of $(A,\D)$ as defined in Definition \ref{d-intrins}. By Theorem \ref{t-cocomm2} we have a
(concretely given) equivalence $D: \Rep_f(A,\D)\rarr\Rep_f G$ of tensor $*$-categories. Thus the
composition $F'=D\circ F: \2C\rarr\Rep_f G$ is the desired equivalence of tensor $*$-categories. It
satisfies $K'\circ F'=K'\circ D\circ F=K\circ F=E$. The final claim is obvious, since the unitary
grouplike elements of $(A,\D)$ are by definition of $(A,\D)$ precisely the unitary natural
isomorphisms of $E$ that are compatible with the tensor structures, in the sense of (\ref{e-nat}). 
\qed

We end this section with some general facts on $R$-matrices for quasitriangular \aqgs.

\bprop
Let $(A,\D)$ be a quasitriangular aqg w.r.t. $R$. Then the following statements hold:
\begin{enumerate}
\item $R_{12}R_{13}R_{23}=R_{23}R_{13}R_{12}$.
\item $(\varepsilon\otimes\iota )R=1=(\iota\otimes\varepsilon )R$.
\item $(S\otimes\iota )R=R^{-1} =(\iota\otimes S^{-1} )(R)$.
\item $(S\otimes S)R=R$.
\end{enumerate}
The first equation is called the `quantum Yang-Baxter equation'. 
\eprop

\prf $R$ is a non-degenerate corepresentation of $(A,\D)$ on $A$, so 2. and the first equation in
3. follow immediately. The second equation in 3. thus follows by noting that $(A, \D^{op})$ is
quasitriangular w.r.t. $\sigma (R)$ with coinverse $S^{-1}$. Claim 4 follows from 3. by
$(S\otimes S)R=(\iota\otimes S)(S\otimes\iota )R=(\iota\otimes S)(R^{-1}) 
=(\iota\otimes S)(\iota\otimes S^{-1})R =R$.
Finally, to prove 1. first note that both the maps $\D^{op}\otimes\iota $ and
$R_{12}(\D\otimes\iota )(\cdot )R_{12}^{-1}$ are non-degenerate homomorphisms from $A\otimes A$
to $M(A\otimes A\otimes A )$, and thus both have unique extensions to $M(A\otimes A)$. 
Because $\D^{op}=R\D (\cdot )R^{-1}$ these maps clearly coincide on $A\otimes A$, and thus on
$M(A\otimes A)$, extensions being unique. Hence
\[ R_{12}R_{13}R_{23}R_{12}^{-1}=R_{12}(\D\otimes\iota)(R)R_{12}^{-1}=(\D^{op}\otimes\iota )(R)
=(\sigma\otimes\iota )(\D\otimes\iota)(R)=(\sigma\otimes\iota )(R_{13}R_{23})=R_{23}R_{13}, \]
as desired.
\qed

\bprop Suppose $(A,\D)$ is a quasitriangular aqg w.r.t. $R$. 
Let $U$ denote the universal corepresentation of $(A,\D)$ as defined in Theorem \ref{t-univ-corep}.
Define $\pi_R :\hat{A}\rightarrow M(A)$
as in Proposition \ref{p-functP}. Then the following hold:
\begin{enumerate}
\item $\pi_R$ is a non-degenerate homomorphism.
\item $(\pi_R\otimes\pi_R )\hat{\D}=\D^{op}\pi_R$.
\item $\D^{op} (\cdot)(\iota\otimes\pi_R )U =(\iota\otimes \pi_R )U\D (\cdot )$.
\end{enumerate}
Conversely, any map $\pi: \hat{A}\rightarrow M(A)$ satisfying these three properties arises from a 
unique $R\in M(A\otimes A)$ with $\pi=\pi_R$ making $(A,\D)$ into a quasitriangular aqg. Moreover,
the map $\pi_R$ is $*$-preserving iff $R$ is unitary.
\eprop

\prf
It suffices to show that statement 2 is equivalent to $(\iota\otimes\D)(R)=R_{13}R_{12}$, 
this being the only non-trivial step not covered by Proposition \ref{p-functP} and Theorem
\ref{t-univ-corr}. But by Theorem \ref{t-univ-corep}, we get
\begin{multline*} 
 (\pi_R\otimes\pi_R )\hat{\D}\omega =(\pi_R\otimes\pi_R )\hat{\D}(\omega\otimes\iota )U
   =(\pi_R\otimes\pi_R )(\omega\otimes\iota\otimes\iota )(\iota\otimes \hat{\D})U
 =(\pi_R\otimes\pi_R )(\omega\otimes\iota\otimes\iota )(U_{12}U_{13}) \\
   =(\omega\otimes\iota\otimes\iota )(\iota\otimes\pi_R\otimes\pi_R )(U_{12}U_{13})
=(\omega\otimes\iota\otimes\iota )[((\iota\otimes\pi_R )U)_{12}((\iota\otimes\pi_R )U)_{13}]
=(\omega\otimes\iota\otimes\iota )(R_{12}R_{13}) 
\end{multline*}
for all $\omega\in\hat{A}$, whereas
\[ \D^{op}\pi_R (\omega )=\D^{op}(\omega\otimes\iota )R 
   =(\omega\otimes\iota\otimes\iota )(\iota\otimes\D^{op})R \]
for all $\omega\in\hat{A}$.
\qed

Note that 1. and 2. just mean that $\pi_R$ is a morphism of multiplier Hopf algebras from
$(\hat{A},\hat{\D })$ to $(A,\D^{op} )$. It is even an aqg morphism iff $R$ is unitary by the last
statement of the proposition.

In Definition \ref{d-corepf} we introduced the tensor product 
$V\times V' =V_{12}V'_{13}\in M(A\otimes B\otimes B')$  
of two unitary corepresentations $V$ and $V'$ on $B$ and $B'$, respectively. 
Another choice of a tensor product which works equally well, is given by
$V\times_{op} V' =V'_{13}V_{12}\in M(A\otimes B\otimes B' )$.
In general there is no relation between these two tensor products. 
In view of the relation between representations of an aqg and corepresentations of its dual aqg
(and the relation between the identity representation and the universal corepresentation), the
following result should not come as a big surprise.

\bprop
Let $(A,\D)$ be an aqg with universal corepresentation $U$, so $hat{U}=\sigma (U)$ is the universal  
corepresentation of $(\hat{A},\hat{\D})$ on $A$. Let $R$ be any element of $M(A\otimes A)$. Then
$\D^{op}(\cdot )R=R\D (\cdot )$ iff 
\[(\hat{U}\times_{op}\hat{U})R_{23} =R_{23}(\hat{U}\times\hat{U}).\]  
\eprop

\prf
By Theorem \ref{t-univ-corep}, we get
\[ \D^{op}(a )R =(\D^{op}(\iota\otimes a)U )R =(\sigma\D (\iota\otimes a)U)R\] 
\[ =((\sigma\otimes a)(\D\otimes\iota )U)R 
   =((\iota\otimes\iota\otimes a )(\sigma\otimes\iota )(U_{13}U_{23}))R
   =(\iota\otimes\iota\otimes a )(U_{23}U_{13}R_{12}) \]
for all $a\in A$, whereas
\[ R\D(a ) =R\D (\iota\otimes a)U  =R(\iota\otimes\iota\otimes a)(\D\otimes\iota )U 
   =R(\iota\otimes\iota\otimes a )(U_{13}U_{23})
   =(\iota\otimes\iota\otimes a )(R_{12}U_{13}U_{23}) \]
for all $a\in A$. Hence $\D^{op}(\cdot )R=R\D (\cdot )$
iff $U_{23}U_{13}R_{12}=R_{12}U_{13}U_{23}$.
Permuting indices, we see that the latter equation is equivalent to
$U_{31}U_{21}R_{23}=R_{23}U_{21}U_{31}$, which again is equivalent to
$\hat{U}_{13}\hat{U}_{12}R_{23}=R_{23}\hat{U}_{12}\hat{U}_{13}$
as desired.
\qed

So given any unitary $R$-matrix $R$, the map $f\in \Hom(A\otimes A,A\otimes A)$ given 
by $f(b)=RbR^*$, for all $b\in A\otimes A$, is an arrow from $\Delta$ to $\Delta^{op}$ and from
$\hat{U}\times\hat{U}$ to $\hat{U}\times_{op}\hat{U}$ in the sense of Remark \ref{r-cat}.
Also $\delta_R (a)=f(1\otimes a)$ for $a\in A$, is a coaction in the sense of Remark \ref{r-actcoact}.
Similarly, the linear map $\delta_R ' :A\rarr M(A\otimes A)$ given by $\delta_R ' (a)=f(a\otimes 1)$, 
for all $a\in A$, is a non-degenerate $*$-homomorphism such that 
$(\delta_R '\otimes\iota )\delta_R ' =(\iota\otimes\D^{op} )\delta_R '$ and 
$(\iota\otimes\varepsilon )\delta_R ' =\iota$.


\bigskip
{\parindent=0pt \parskip=4pt Addresses of the authors:

\smallskip Michael M\"uger, Korteweg-de Vries Instituut, Universiteit van Amsterdam, Amsterdam, Netherlands. \\
E-mail: mmueger@science.uva.nl. \\ 

John E. Roberts, Dipartimento di Matematica, Universit\`{a} di Roma II `Tor Vergata', Roma, Italy. \\
E-mail: roberts@mat.uniroma2.it. \\

Lars Tuset, Faculty of Engineering, Oslo University College, Oslo, Norway. \\
E-mail: Lars.Tuset@iu.hio.no.}


\end{document}